\documentstyle[leqno]{article} \makeatletter
\ifx\@auxdefsloaded\relax\endinput
\else\let\@auxdefsloaded\relax\fi

\def\newenvironment{%
   \@ifnextchar *{\@@newenv{\global\@ignoretrue}}{\@@newenv{}*}}

\def\@@newenv#1*#2{%
   \@ifnextchar [{\@newenv{#1}{#2}}{\@newenv{#1}{#2}[0]}}

\long\def\@newenv#1#2[#3]#4#5{%
   \expandafter\newcommand\csname#2\endcsname[#3]{#4}%
   \expandafter\long\expandafter\def\csname end#2\endcsname{#5#1}}

\def\renewenvironment{%
   \@ifnextchar *{\@@renewenv{\global\@ignoretrue}}{\@@renewenv{}*}}

\def\@@renewenv#1*#2{%
   \@ifnextchar [{\@renewenv{#1}{#2}}{\@renewenv{#1}{#2}[0]}}

\long\def\@renewenv#1#2[#3]#4#5{%
   \expandafter\renewcommand\csname#2\endcsname[#3]{#4}%
   \expandafter\long\expandafter\def\csname end#2\endcsname{#5#1}}

\def\newoptcommand#1#2{%
   \@ifnextchar [{\@optargdef#1#2}{\@optargdef#1#2[1]}}

\def\renewoptcommand#1#2{%
   \edef\@tempa{\expandafter\@cdr\string#1\@nil}%
   \@ifundefined{\@tempa}{%
      \@latexerr{\string#1\space undefined}\@ehc}{}%
   \@ifnextchar [{\@reoptargdef#1#2}{\@reoptargdef#1#2[1]}}

\long\def\@optargdef#1#2[#3]#4{%
   \@ifdefinable #1{\@reoptargdef#1#2[#3]{#4}}}

\long\def\@reoptargdef#1#2[#3]#4{%
   \@tempcnta#3\relax \@tempcntb \@ne
   \let#1\relax \let\@tempa\relax
   \edef\@tempb{\long\def\csname\string#1\endcsname
      [\@tempa\the\@tempcntb]}%
   \advance\@tempcntb \@ne \advance\@tempcnta \m@ne
   \@whilenum\@tempcnta>0\do{%
      \edef\@tempb{\@tempb\@tempa\the\@tempcntb}%
      \advance\@tempcntb \@ne \advance\@tempcnta \m@ne}%
   \let\@tempa=##\@tempb{#4}%
   \def#1{\@ifnextchar [{\csname\string#1\endcsname}{%
      \csname\string#1\endcsname[#2]}}}

\def\newoptenvironment{%
   \@ifnextchar *{\@@newoptenv{\global\@ignoretrue}}{%
      \@@newoptenv{}*}}

\def\@@newoptenv#1*#2#3{%
   \@ifnextchar [{\@newoptenv{#1}{#2}{#3}}{%
      \@newoptenv{#1}{#2}{#3}[0]}}

\long\def\@newoptenv#1#2#3[#4]#5#6{%
   \expandafter\newoptcommand\csname#2\endcsname{#3}[#4]{#5}%
   \expandafter\long\expandafter\def\csname end#2\endcsname{#6#1}}

\def\renewoptenvironment{%
   \@ifnextchar *{\@@renewoptenv{\global\@ignoretrue}}{%
      \@@renewoptenv{}*}}

\def\@@renewoptenv#1*#2#3{%
   \@ifnextchar [{\@renewoptenv{#1}{#2}{#3}}{%
      \@renewoptenv{#1}{#2}{#3}[0]}}

\long\def\@renewoptenv#1#2#3[#4]#5#6{%
   \expandafter\renewoptcommand\csname#2\endcsname{#3}[#4]{#5}%
   \expandafter\long\expandafter\def\csname end#2\endcsname{#6#1}}

\newcounter{keepoptional}
\newcounter{optnestctr}

\newoptenvironment*{optional}{1}[1]{%
   \ifnum#1>\value{keepoptional}\relax
      \setcounter{optnestctr}{0}\@powerup\expandafter\@endopt\fi
   \ignorespaces}{}

\def\@powerup{\catcode`\{=12 \catcode`\}=12 \catcode`\\=12 \relax}
\def\@powerdown{\catcode`\{=1 \catcode`\}=2 \catcode`\\=0 \relax}
\begingroup \catcode`|=0 \catcode`[=1 \catcode`]=2 \@powerup
   |long|gdef|@endopt#1\end{optional}[%
      |@beginopt#1\begin{optional}*]%
   |long|gdef|@beginopt#1\begin{optional}[%
      |@ifstar[%
         |ifnum|value[optnestctr]>0|relax
            |addtocounter[optnestctr][-1]|let|@zzzap=|@endopt
         |else
            |@powerdown|end[optional]|let|@zzzap=|relax|fi
	 |@zzzap]%
      [
         |addtocounter[optnestctr][1]|@beginopt]]%
|endgroup

\makeatother
\makeatletter
\ifx\@auxdefsloaded\relax
\else \input{auxdefs.sty}\fi
\newskip\dgARROWLENGTH  \dgARROWLENGTH=2.5em\relax
\newskip\dgHORIZPAD     \dgHORIZPAD=1em\relax
\newskip\dgVERTPAD      \dgVERTPAD=2ex\relax
\newskip\dgLABELOFFSET  \dgLABELOFFSET=.7ex\relax
\newcount\dgARROWPARTS  \dgARROWPARTS=4\relax
\newcommand{\dgeverynode}{\displaystyle}
\newcommand{\dgeverylabel}{\scriptstyle}
\newskip\dgDOTSPACING   \dgDOTSPACING=0.35em
\newskip\dgDOTSIZE      \dgDOTSIZE=1.5\fontdimen8\tenln
\newcount\dgMAXSQUARE   \dgMAXSQUARE=4\relax
\newcount\dgMINDBLSQ    \dgMINDBLSQ=10\relax
\newskip\dgCOLUMNWIDTH  \dgCOLUMNWIDTH=2em\relax
\chardef\f@ur=4
\@namedef{dgo@}{\let\dg@VECTOR=\vector
   \dg@LBLPOS=\dgARROWPARTS \divide\dg@LBLPOS\tw@}%
\@namedef{dgo@-}{\let\dg@VECTOR=\line}%
\@namedef{dgo@!}{\let\dg@VECTOR=\dg@novector}%
\@namedef{dgo@1}{\dg@LBLPOS=\@ne\relax}%
\@namedef{dgo@2}{\dg@LBLPOS=\tw@\relax}%
\@namedef{dgo@3}{\dg@LBLPOS=\thr@@\relax}%
\@namedef{dgo@4}{\dg@LBLPOS=\f@ur\relax}%
\@namedef{dgo@5}{\dg@LBLPOS=5\relax}%
\@namedef{dgo@6}{\dg@LBLPOS=6\relax}%
\@namedef{dgo@7}{\dg@LBLPOS=7\relax}%
\@namedef{dgo@8}{\dg@LBLPOS=8\relax}%
\@namedef{dgo@9}{\dg@LBLPOS=9\relax}%
\def\dgt@e{\dg@DX=\@ne \dg@DY=\z@ \dg@SIZE=\@ne}%
\def\dgt@w{\dg@DX=\m@ne \dg@DY=\z@ \dg@SIZE=\@ne}%
\def\dgt@n{\dg@DX=\z@ \dg@DY=\@ne \dg@SIZE=\@ne}%
\def\dgt@s{\dg@DX=\z@ \dg@DY=\m@ne \dg@SIZE=\@ne}%
\def\dgt@ne{\dg@DX=\@ne \dg@DY=\@ne \dg@SIZE=\@ne}%
\def\dgt@se{\dg@DX=\@ne \dg@DY=\m@ne \dg@SIZE=\@ne}%
\def\dgt@nw{\dg@DX=\m@ne \dg@DY=\@ne \dg@SIZE=\@ne}%
\def\dgt@sw{\dg@DX=\m@ne \dg@DY=\m@ne \dg@SIZE=\@ne}%
\def\dgt@nne{\dg@DX=\@ne \dg@DY=\tw@ \dg@SIZE=\@ne}%
\def\dgt@nnw{\dg@DX=\m@ne \dg@DY=\tw@ \dg@SIZE=\@ne}%
\def\dgt@sse{\dg@DX=\@ne \dg@DY=-\tw@ \dg@SIZE=\@ne}%
\def\dgt@ssw{\dg@DX=\m@ne \dg@DY=-\tw@ \dg@SIZE=\@ne}%
\def\dgt@ene{\dg@DX=\tw@ \dg@DY=\@ne \dg@SIZE=\tw@}%
\def\dgt@ese{\dg@DX=\tw@ \dg@DY=\m@ne \dg@SIZE=\tw@}%
\def\dgt@wnw{\dg@DX=-\tw@ \dg@DY=\@ne \dg@SIZE=\tw@}%
\def\dgt@wsw{\dg@DX=-\tw@ \dg@DY=\m@ne \dg@SIZE=\tw@}%
\def\dggeometry{
   \dg@ZTEMP=\dg@XGRID \multiply\dg@ZTEMP\tw@
   \ifnum\dg@YGRID=\z@ \dg@ZTEMP=\tw@
   \else \divide\dg@ZTEMP\dg@YGRID \fi
   \ifnum\dg@ZTEMP>\f@ur \dg@ZTEMP=\f@ur \fi
   \ifnum\dg@ZTEMP<\@ne \dg@ZTEMP=\@ne \fi
   \unitlength=2sp\relax
   \ifnum\dg@ZTEMP<\tw@
      \advance\dg@ZTEMP\@ne
      \multiply\unitlength\dg@YGRID
   \else
      \multiply\unitlength\dg@XGRID \divide\unitlength\dg@ZTEMP
   \fi
   \dg@XGRID=\dg@ZTEMP \dg@YGRID=\tw@
   \dg@rmcommondiv\tw@\dg@XGRID\dg@YGRID
   \divide\unitlength\dg@YGRID \divide\unitlength\@m\relax}
\@namedef{dgo@..}{\let\dg@VECTOR=\dg@dotvector}%

\def\dg@dotvector(#1,#2)#3{%
   \begingroup
   \dg@XTEMP=#1\relax \dg@YTEMP=#2\relax
   \let\dg@NDOTS=\dg@XEND \let\dg@DOTDIAM=\dg@WEND
   \dg@NDOTS=\unitlength \multiply\dg@NDOTS #3\relax
   \dg@ZTEMP=\dg@YTEMP \dg@changesign\dg@YTEMP\dg@ZTEMP
   \ifnum\dg@XTEMP>\z@
      \ifnum\dg@YTEMP>\dg@XTEMP
         \multiply\dg@NDOTS\dg@YTEMP \divide\dg@NDOTS\dg@XTEMP \fi
   \else\ifnum\dg@XTEMP<\z@
      \ifnum\dg@YTEMP>-\dg@XTEMP
         \multiply\dg@NDOTS\dg@YTEMP \divide\dg@NDOTS-\dg@XTEMP \fi
   \fi\fi
   \dg@YTEMP=\dg@ZTEMP
   \divide\dg@NDOTS\dgDOTSPACING
   \ifnum\dg@NDOTS>\z@\else \dg@NDOTS=\@ne \fi
   \dg@ZTEMP=\unitlength \multiply\dg@ZTEMP #3\relax
   \divide\dg@ZTEMP\dg@NDOTS
   \ifnum\dg@XTEMP=\z@
      \dg@changesign\dg@ZTEMP\dg@YTEMP \dg@YTEMP=\dg@ZTEMP
   \else
      \dg@changesign\dg@ZTEMP\dg@XTEMP
      \multiply\dg@YTEMP\dg@ZTEMP \divide\dg@YTEMP\dg@XTEMP
      \dg@XTEMP=\dg@ZTEMP
   \fi
   \divide\dg@XTEMP\unitlength \divide\dg@YTEMP\unitlength
   \begin{picture}(0,0)
      \dg@DOTDIAM=\dgDOTSIZE \divide\dg@DOTDIAM\unitlength
      \multiput(0,0)(\dg@XTEMP,\dg@YTEMP){\dg@NDOTS}{%
         \circle*{\dg@DOTDIAM}}%
      \multiply\dg@XTEMP\dg@NDOTS \multiply\dg@YTEMP\dg@NDOTS
      \put(\dg@XTEMP,\dg@YTEMP){\vector(#1,#2){0}}%
   \end{picture}%
   \endgroup}%
\newif\ifdg@LATEXGEOM \dg@LATEXGEOMfalse

\@ifundefined{lamsvector}{}{%
   \@namedef{dgo@}{%
      \let\dg@VECTOR=\lamsvector
      \lamsshaft{?}\lamssource{?}\lamstarget{?}%
      \dg@LBLPOS=\dgARROWPARTS \divide\dg@LBLPOS\tw@\relax}%
   \@namedef{dgo@..}{\lamsshaft{.}}%
   \@namedef{dgo@=}{\lamsshaft{=}}%
   \@namedef{dgo@-}{\lamstarget{0}}%
   \@namedef{dgo@'}{\lamstarget{'}}%
   \@namedef{dgo@`}{\lamstarget{`}}%
   \@namedef{dgo@A}{\lamstarget{A}}%
   \@namedef{dgo@V}{\lamssource{V}}%
   \@namedef{dgo@J}{\lamssource{J}}%
   \@namedef{dgo@L}{\lamssource{L}}%
   \@namedef{dgo@S}{\lamssource{S}}%
   \def\dggeometry{
      \dg@ZTEMP=\dg@XGRID \multiply\dg@ZTEMP\tw@
      \ifnum\dg@YGRID=\z@ \dg@ZTEMP=\tw@
      \else \divide\dg@ZTEMP\dg@YGRID \fi
      \ifnum\dg@ZTEMP>6\relax \dg@ZTEMP=6\relax \fi
      \ifdg@LATEXGEOM\ifnum\dg@ZTEMP>\f@ur \dg@ZTEMP=\f@ur \fi\fi
      \ifnum\dg@ZTEMP<\@ne \dg@ZTEMP=\@ne \fi
      \unitlength=2sp\relax
      \ifnum\dg@ZTEMP<\tw@
         \advance\dg@ZTEMP\@ne
         \multiply\unitlength\dg@YGRID
      \else
         \multiply\unitlength\dg@XGRID \divide\unitlength\dg@ZTEMP
      \fi
      \dg@XGRID=\dg@ZTEMP \dg@YGRID=\tw@
      \dg@rmcommondiv\tw@\dg@XGRID\dg@YGRID
      \divide\unitlength\dg@YGRID \divide\unitlength\@m
      \dg@LATEXGEOMfalse}
   \def\dgt@nee{\dg@DX=\tw@ \dg@DY=\@ne \dg@SIZE=\tw@}%
   \def\dgt@see{\dg@DX=\tw@ \dg@DY=\m@ne \dg@SIZE=\tw@}%
   \def\dgt@nww{\dg@DX=-\tw@ \dg@DY=\@ne \dg@SIZE=\tw@}%
   \def\dgt@sww{\dg@DX=-\tw@ \dg@DY=\m@ne \dg@SIZE=\tw@}%
   \def\dgt@nnne{\dg@DX=\@ne \dg@DY=\thr@@ \dg@SIZE=\@ne}%
   \def\dgt@nnnw{\dg@DX=\m@ne \dg@DY=\thr@@ \dg@SIZE=\@ne}%
   \def\dgt@sssw{\dg@DX=\m@ne \dg@DY=-\thr@@ \dg@SIZE=\@ne}%
   \def\dgt@ssse{\dg@DX=\@ne \dg@DY=-\thr@@ \dg@SIZE=\@ne}%
   \def\dgt@nnnee{\dg@DX=\tw@ \dg@DY=\thr@@ \dg@SIZE=\tw@}%
   \def\dgt@nnnww{\dg@DX=-\tw@ \dg@DY=\thr@@ \dg@SIZE=\tw@}%
   \def\dgt@sssww{\dg@DX=-\tw@ \dg@DY=-\thr@@ \dg@SIZE=\tw@}%
   \def\dgt@sssee{\dg@DX=\tw@ \dg@DY=-\thr@@ \dg@SIZE=\tw@}%
   \def\dgt@nneee{\dg@DX=\thr@@ \dg@DY=\tw@ \dg@SIZE=\thr@@}%
   \def\dgt@nnwww{\dg@DX=-\thr@@ \dg@DY=\tw@ \dg@SIZE=\thr@@}%
   \def\dgt@sswww{\dg@DX=-\thr@@ \dg@DY=-\tw@ \dg@SIZE=\thr@@}%
   \def\dgt@sseee{\dg@DX=\thr@@ \dg@DY=-\tw@ \dg@SIZE=\thr@@}%
   \def\dgt@neee{\dg@DX=\thr@@ \dg@DY=\@ne \dg@SIZE=\thr@@
      \global\dg@LATEXGEOMtrue}%
   \def\dgt@nwww{\dg@DX=-\thr@@ \dg@DY=\@ne \dg@SIZE=\thr@@
      \global\dg@LATEXGEOMtrue}%
   \def\dgt@swww{\dg@DX=-\thr@@ \dg@DY=\m@ne \dg@SIZE=\thr@@
      \global\dg@LATEXGEOMtrue}%
   \def\dgt@seee{\dg@DX=\thr@@ \dg@DY=\m@ne \dg@SIZE=\thr@@
      \global\dg@LATEXGEOMtrue}%
}
\newcount\dg@HORIZ      \newcount\dg@VERT
\newcount\dg@XLPAD      \newcount\dg@YBPAD
\newcount\dg@XRPAD      \newcount\dg@YTPAD
\newcount\dg@X          \newcount\dg@Y
\newcount\dg@XGRID      \newcount\dg@YGRID
\newcount\dg@SIZE
\newcount\dg@USERSIZE
\newcount\dg@DX         \newcount\dg@DY
\newcount\dg@XLBL       \newcount\dg@YLBL
\newcount\dg@XOFFSET    \newcount\dg@YOFFSET
\newcount\dg@LBLOFF
\newcount\dg@XLBLOFF    \newcount\dg@YLBLOFF
\newcount\dg@LBLPOS
\newcount\dg@XTEMP      \newcount\dg@YTEMP
\newcount\dg@ZTEMP
\newcount\dg@XNODE      \newcount\dg@YNODE
\newcount\dg@XEND       \newcount\dg@YEND
\newcount\dg@WEND       \newcount\dg@HEND
\newcount\dg@COUNT
\newbox\dg@NODEBOX
\newoptenvironment*{diagram}{}[1]{%
   \global\advance\dg@COUNT\@ne \typeout{[diagram \the\dg@COUNT:}%
   \let\node=\dg@node \let\\=\dg@cr \let\arrow=\dg@arrow
   \def\dg@BIGNODE{#1}%
   \ifx\dg@BIGNODE\@empty\else
      \setbox\dg@NODEBOX=\hbox{$\dgeverynode{#1}$}%
      \dg@XTEMP=\wd\dg@NODEBOX
      \dg@YTEMP=\ht\dg@NODEBOX \advance\dg@YTEMP\dp\dg@NODEBOX
      \ifnum\dg@YTEMP=\z@ \dg@ZTEMP=\z@
      \else \dg@ZTEMP=\dg@XTEMP \divide\dg@ZTEMP\dg@YTEMP\fi
      \ifnum\dg@ZTEMP>\dgMAXSQUARE
         \ifnum\dg@ZTEMP<\dgMINDBLSQ \dg@XGRID=\thr@@ \dg@YGRID=\tw@
         \else \dg@XGRID=\tw@ \dg@YGRID=\@ne \fi
      \else \dg@XGRID=\@ne \dg@YGRID=\@ne \fi
      \advance\dg@XTEMP\dgHORIZPAD \advance\dg@YTEMP\dgVERTPAD
      \dg@XLPAD=\dg@XTEMP \divide\dg@XLPAD\tw@
      \dg@XRPAD=\dg@XTEMP \divide\dg@XRPAD\tw@
      \dg@YBPAD=\dg@YTEMP \divide\dg@YBPAD\tw@
      \dg@YTPAD=\dg@YTEMP \divide\dg@YTPAD\tw@
      \advance\dg@XTEMP\dgARROWLENGTH
      \divide\dg@XTEMP\dg@XGRID \divide\dg@XTEMP\@m
      \unitlength=1sp\relax \multiply\unitlength\dg@XTEMP
   \fi
   \typeout{saving...}%
   \dg@X=-\@ne \dg@Y=\z@ \dg@HORIZ=\z@\relax
   \let\dg@SLIST=\@empty
   \let\dg@NLIST=\@empty \let\dg@ALIST=\@empty
   \let\dg@PASS=\dg@savepass
}{
   \dg@VERT=-\dg@Y\relax
   \typeout{calculating...}%
   \ifx\dg@BIGNODE\@empty
      \dg@XGRID=\z@ \dg@YGRID=\z@
      \dg@XLPAD=\z@ \dg@XRPAD=\z@ \dg@YBPAD=\z@ \dg@YTPAD=\z@
      \let\dg@PASS=\dg@geompass
      \dg@NLIST \dg@ALIST
      \ifnum\dg@XGRID>\z@\else \dg@XGRID=\quad\relax \fi
      \ifnum\dg@YGRID>\z@\else \dg@YGRID=\quad\relax \fi
      \dggeometry
      \ifdim\unitlength>\z@\else \unitlength=\quad \fi
      \ifnum\dg@XGRID>\z@\else \dg@XGRID=\@ne \fi
      \ifnum\dg@YGRID>\z@\else \dg@YGRID=\@ne \fi
   \fi
   \dg@LBLOFF=\dgLABELOFFSET \divide\dg@LBLOFF\unitlength
   \multiply\dg@HORIZ\@m \multiply\dg@HORIZ\dg@XGRID
   \multiply\dg@VERT\@m \multiply\dg@VERT\dg@YGRID
   \divide\dg@XLPAD\unitlength \divide\dg@XRPAD\unitlength
   \divide\dg@YBPAD\unitlength \divide\dg@YTPAD\unitlength
   \typeout{drawing...}%
   \let\dg@PASS=\dg@drawpass
   \hspace*{\dg@XLPAD\unitlength}%
   \vcenter{%
      \hsize=0pt\relax\parindent=0pt\relax
      \parskip=0pt\relax\baselineskip=0pt\relax
      \vspace*{\dg@YTPAD\unitlength}%
      \begin{picture}(0,0)(0,0)\dg@NLIST\dg@ALIST\end{picture}%
      \raisebox{\z@}[\z@][\dg@VERT\unitlength]{}%
      \vspace*{\dg@YBPAD\unitlength}%
   }
   \hspace*{\dg@HORIZ\unitlength}%
   \hspace*{\dg@XRPAD\unitlength}%
   \typeout{done]}}%

\def\dg@savepass{s}
\def\dg@geompass{g}
\def\dg@drawpass{d}

\newoptcommand{\dg@node}{\@ne}[2]{%
   \ifx\dg@PASS\dg@savepass
      \dg@XTEMP=#1\relax
      \ifnum\dg@XTEMP<\@ne \dg@XTEMP=\@ne\fi
      \advance\dg@X\dg@XTEMP
      \ifnum\dg@HORIZ<\dg@X \dg@HORIZ=\dg@X \fi
      \setbox\dg@NODEBOX=\hbox{$\dgeverynode{#2}$}%
      \dg@XTEMP=\wd\dg@NODEBOX \advance\dg@XTEMP\dgHORIZPAD
      \dg@YTEMP=\ht\dg@NODEBOX \advance\dg@YTEMP\dp\dg@NODEBOX
      \advance\dg@YTEMP\dgVERTPAD
      \toks\z@=\expandafter{\dg@SLIST}%
      \edef\dg@SLIST{\the\toks\z@
         ,\noexpand\dg@XNODE=\number\dg@X\noexpand\relax
         \noexpand\dg@YNODE=\number\dg@Y\noexpand\relax
         \noexpand\dg@XTEMP=\number\dg@XTEMP\noexpand\relax
         \noexpand\dg@YTEMP=\number\dg@YTEMP\noexpand\relax}%
      \toks\z@=\expandafter{\dg@NLIST}%
      \toks\tw@={\dg@node{#2}}%
      \edef\dg@NLIST{\the\toks\z@
         \noexpand\dg@X=\number\dg@X\noexpand\relax
         \noexpand\dg@Y=\number\dg@Y\noexpand\relax
         \the\toks\tw@}%
   \else\ifx\dg@PASS\dg@geompass
      \ifnum\dg@X=\z@
         \dg@getnodesize
            {\dg@SLIST}{\dg@X}{\dg@Y}{\dg@WEND}{\dg@HEND}%
         \divide\dg@WEND\tw@
         \ifnum\dg@XLPAD<\dg@WEND \dg@XLPAD=\dg@WEND \fi\fi
      \ifnum\dg@X=\dg@HORIZ
         \dg@getnodesize
            {\dg@SLIST}{\dg@X}{\dg@Y}{\dg@WEND}{\dg@HEND}%
         \divide\dg@WEND\tw@
         \ifnum\dg@XRPAD<\dg@WEND \dg@XRPAD=\dg@WEND \fi\fi
      \ifnum\dg@Y=\z@
         \dg@getnodesize
            {\dg@SLIST}{\dg@X}{\dg@Y}{\dg@WEND}{\dg@HEND}%
         \divide\dg@HEND\tw@
         \ifnum\dg@YTPAD<\dg@HEND \dg@YTPAD=\dg@HEND \fi\fi
      \ifnum\dg@Y=-\dg@VERT\relax
         \dg@getnodesize
            {\dg@SLIST}{\dg@X}{\dg@Y}{\dg@WEND}{\dg@HEND}%
         \divide\dg@HEND\tw@
         \ifnum\dg@YBPAD<\dg@HEND \dg@YBPAD=\dg@HEND \fi\fi
   \else\ifx\dg@PASS\dg@drawpass
      \dg@XNODE=\dg@X \multiply\dg@XNODE\@m
      \multiply\dg@XNODE\dg@XGRID
      \dg@YNODE=\dg@Y \multiply\dg@YNODE\@m
      \multiply\dg@YNODE\dg@YGRID
      \setbox\dg@NODEBOX=\hbox{$\dgeverynode{#2}$}%
      \put(\dg@XNODE,\dg@YNODE){\dg@makebox{\box\dg@NODEBOX}}%
   \fi\fi\fi}%

\newoptcommand{\dg@cr}{\@ne}[1]{%
   \ifx\dg@PASS\dg@savepass
      \dg@YTEMP=#1\relax
      \ifnum\dg@YTEMP<\@ne \dg@YTEMP=\@ne \fi
      \advance\dg@Y -\dg@YTEMP\relax
      \dg@X=-\@ne\relax\fi}%
\newoptcommand{\dg@arrow}{\@ne}[2]{%
   \begingroup
   \dg@USERSIZE=#1\relax
   \ifnum\dg@USERSIZE<\@ne \dg@USERSIZE=\@ne \fi
   \dg@parse{#2}%
   \ifx\dg@PASS\dg@savepass
      \ifx\dg@label\dgl@b \let\dg@label=\dgl@t \fi
      \ifx\dg@label\dgl@r \let\dg@label=\dgl@l \fi
      \let\dg@process=\dg@save
   \else\ifx\dg@PASS\dg@geompass
      \let\dg@process=\dg@ignore
      \dg@geomcalc
   \else\ifx\dg@PASS\dg@drawpass
      \let\dg@process=\dg@draw
      \dg@drawcalc
   \fi\fi\fi
   \dg@label{\dg@process{#1}{#2}}}%

\newoptcommand{\arrow}{\@ne}[2]{%
   \dg@parse{#2}%
   \ifx\dg@label\dgl@b \let\dg@label=\dgl@t \fi
   \ifx\dg@label\dgl@r \let\dg@label=\dgl@l \fi
   \dg@label{\dg@textarrow{#1}{#2}}}%

\def\dg@textarrow#1#2#3#4{%
   \mathop{{\dgHORIZPAD=0pt\relax\dgVERTPAD=0pt\relax
      \begin{diagram}
         \node{}\arrow[#1]{#2}{#3}{#4}\node{}
      \end{diagram}}}}

\def\dg@parse#1{%
   \let\dg@label=\dgl@ \dgo@
   \let\dg@type=\@empty \let\dg@lbltype=\@empty
   \@for\dg@list:=#1\do{%
      \ifx\dg@type\@empty \let\dg@type=\dg@list
      \else\ifx\dg@lbltype\@empty \let\dg@lbltype=\dg@list
         \@ifundefined{dgo@\dg@list}{}{\@nameuse{dgo@\dg@list}}%
      \else
         \@ifundefined{dgo@\dg@list}{}{\@nameuse{dgo@\dg@list}}%
      \fi\fi}%
   \@ifundefined{dgt@\dg@type}{\dgt@e}{\@nameuse{dgt@\dg@type}}%
   \@ifundefined{dgl@\dg@lbltype}{}{%
      \dg@letname\dg@label{dgl@\dg@lbltype}}}

\def\dg@draw#1#2#3#4{%
   \put(\dg@X,\dg@Y){\dg@makebox{%
      \begin{picture}(0,0)%
         \thinlines
         \put(\dg@XOFFSET,\dg@YOFFSET){%
            \dg@VECTOR(\dg@DX,\dg@DY){\dg@SIZE}}%
         \put(\dg@XLBL,\dg@YLBL){\dg@makebox{%
            \begin{picture}(0,0)%
               \put(\dg@XLBLOFF,\dg@YLBLOFF){%
                  \dg@makebox[\dg@LBLONE]{$\dgeverylabel{#3}$}}%
               \put(-\dg@XLBLOFF,-\dg@YLBLOFF){%
                  \dg@makebox[\dg@LBLTWO]{$\dgeverylabel{#4}$}}%
            \end{picture}}}%
      \end{picture}}}%
   \endgroup}%

\def\dg@save#1#2#3#4{%
   \endgroup 
   \toks\z@=\expandafter{\dg@ALIST}%
   \toks\tw@={\dg@arrow[#1]{#2}{#3}{#4}}%
   \edef\dg@ALIST{\the\toks\z@%
      \noexpand\dg@X=\number\dg@X\noexpand\relax
      \noexpand\dg@Y=\number\dg@Y\noexpand\relax
      \the\toks\tw@}}%

\def\dg@ignore#1#2#3#4{\endgroup}

\def\dg@geomcalc{%
   \dg@XEND=\dg@SIZE \multiply\dg@XEND\dg@USERSIZE
   \ifnum\dg@DX=\z@
      \dg@YEND=\dg@XEND \dg@XEND=\z@
      \dg@changesign\dg@YEND\dg@DY
   \else
      \dg@changesign\dg@XEND\dg@DX \dg@YEND=\dg@XEND
      \multiply\dg@YEND\dg@DY \divide\dg@YEND\dg@DX
   \fi
   \advance\dg@XEND\dg@X \advance\dg@YEND\dg@Y
   \dg@getnodesize
      {\dg@SLIST}{\dg@XEND}{\dg@YEND}{\dg@WEND}{\dg@HEND}%
   \dg@XOFFSET=\dg@WEND \dg@YOFFSET=\dg@HEND
   \dg@getnodesize
      {\dg@SLIST}{\dg@X}{\dg@Y}{\dg@WEND}{\dg@HEND}%
   \advance\dg@XOFFSET\dg@WEND \divide\dg@XOFFSET\tw@
   \advance\dg@YOFFSET\dg@HEND \divide\dg@YOFFSET\tw@
   \dg@XTEMP=\dgARROWLENGTH \dg@YTEMP=\dgARROWLENGTH
   \ifnum\dg@DX>\z@
      \dg@ZTEMP=\dg@DX \multiply\dg@XTEMP\dg@DX
   \else \dg@ZTEMP=-\dg@DX \multiply\dg@XTEMP -\dg@DX \fi
   \ifnum\dg@DY>\z@
      \advance\dg@ZTEMP\dg@DY \multiply\dg@YTEMP\dg@DY
   \else \advance\dg@ZTEMP -\dg@DY \multiply\dg@YTEMP -\dg@DY\fi
   \ifnum\dg@ZTEMP=\z@\else
      \divide\dg@XTEMP\dg@ZTEMP \divide\dg@YTEMP\dg@ZTEMP
      \advance\dg@XOFFSET\dg@XTEMP \advance\dg@YOFFSET\dg@YTEMP
   \fi
   \divide\dg@XOFFSET\dg@SIZE \divide\dg@XOFFSET\dg@USERSIZE
   \divide\dg@YOFFSET\dg@SIZE \divide\dg@YOFFSET\dg@USERSIZE
   \ifnum\dg@DX=\z@ \dg@XOFFSET=\z@ \fi
   \ifnum\dg@DY=\z@ \dg@YOFFSET=\z@ \fi
   \ifnum\dg@XGRID<\dg@XOFFSET \global\dg@XGRID=\dg@XOFFSET\fi
   \ifnum\dg@YGRID<\dg@YOFFSET \global\dg@YGRID=\dg@YOFFSET\fi
   \relax}

\def\dg@drawcalc{%
   \dg@XEND=\dg@SIZE \multiply\dg@XEND\dg@USERSIZE
   \ifnum\dg@DX=\z@
      \dg@YEND=\dg@XEND \dg@XEND=\z@
      \dg@changesign\dg@YEND\dg@DY
   \else
      \dg@changesign\dg@XEND\dg@DX \dg@YEND=\dg@XEND
      \multiply\dg@YEND\dg@DY \divide\dg@YEND\dg@DX
   \fi
   \advance\dg@XEND\dg@X \advance\dg@YEND\dg@Y
   \dg@getnodesize
      {\dg@SLIST}{\dg@XEND}{\dg@YEND}{\dg@WEND}{\dg@HEND}%
   \divide\dg@WEND\unitlength \divide\dg@HEND\unitlength
   \multiply\dg@DX\dg@XGRID \multiply\dg@DY\dg@YGRID
   \dg@rmcommondiv\tw@\dg@DX\dg@DY
   \dg@rmcommondiv\tw@\dg@DX\dg@DY 
   \dg@rmcommondiv\thr@@\dg@DX\dg@DY
   \multiply\dg@SIZE\dg@USERSIZE \multiply\dg@SIZE\@m
   \ifnum\dg@DX=\z@
      \multiply\dg@SIZE\dg@YGRID
      \divide\dg@HEND\tw@ \advance\dg@SIZE -\dg@HEND
      \dg@getnodesize
         {\dg@SLIST}{\dg@X}{\dg@Y}{\dg@WEND}{\dg@YOFFSET}%
      \divide\dg@YOFFSET\unitlength \divide\dg@YOFFSET\tw@
      \advance\dg@SIZE -\dg@YOFFSET
      \dg@XOFFSET=\z@
      \def\dg@LBLONE{r}\def\dg@LBLTWO{l}%
      \dg@XLBL=\z@ \dg@YLBL=\dg@SIZE
      \multiply\dg@YLBL\dg@LBLPOS
      \divide\dg@YLBL\dgARROWPARTS\relax
      \advance\dg@YLBL\dg@YOFFSET
      \dg@changesign\dg@YLBL\dg@DY
      \dg@changesign\dg@YOFFSET\dg@DY
   \else
      \multiply\dg@SIZE\dg@XGRID
      \ifnum\dg@DY=\z@
         \divide\dg@WEND\tw@ \advance\dg@SIZE -\dg@WEND
         \dg@getnodesize
            {\dg@SLIST}{\dg@X}{\dg@Y}{\dg@XOFFSET}{\dg@HEND}%
         \divide\dg@XOFFSET\unitlength \divide\dg@XOFFSET\tw@
         \advance\dg@SIZE -\dg@XOFFSET
         \dg@YOFFSET=\z@
         \def\dg@LBLONE{b}\def\dg@LBLTWO{t}%
         \dg@YLBL=\z@ \dg@XLBL=\dg@SIZE
         \multiply\dg@XLBL\dg@LBLPOS
         \divide\dg@XLBL\dgARROWPARTS\relax
         \advance\dg@XLBL\dg@XOFFSET
         \dg@changesign\dg@XLBL\dg@DX
         \dg@changesign\dg@XOFFSET\dg@DX
      \else
         \divide\dg@WEND\tw@ \divide\dg@HEND\tw@
         \multiply\dg@HEND\dg@DX \divide\dg@HEND\dg@DY
         \ifnum\dg@HEND<\z@ \multiply\dg@HEND\m@ne \fi
         \ifnum\dg@WEND<\dg@HEND \advance\dg@SIZE -\dg@WEND
         \else \advance\dg@SIZE -\dg@HEND \fi
         \dg@getnodesize
            {\dg@SLIST}{\dg@X}{\dg@Y}{\dg@WEND}{\dg@HEND}%
         \divide\dg@WEND\unitlength \divide\dg@WEND\tw@
         \divide\dg@HEND\unitlength \divide\dg@HEND\tw@
         \multiply\dg@HEND\dg@DX \divide\dg@HEND\dg@DY
         \ifnum\dg@HEND<\z@ \multiply\dg@HEND\m@ne \fi
         \ifnum\dg@WEND<\dg@HEND \dg@XOFFSET=\dg@WEND
         \else \dg@XOFFSET=\dg@HEND \fi
         \advance\dg@SIZE -\dg@XOFFSET
         \dg@changesign\dg@XOFFSET\dg@DX
         \dg@YOFFSET=\dg@XOFFSET
         \multiply\dg@YOFFSET\dg@DY \divide\dg@YOFFSET\dg@DX
         \def\dg@LBLONE{br}\def\dg@LBLTWO{tl}%
         \ifnum\dg@DX<\z@ \ifnum\dg@DY>\z@
            \def\dg@LBLONE{bl}\def\dg@LBLTWO{tr}\fi\fi
         \ifnum\dg@DX>\z@ \ifnum\dg@DY<\z@
            \def\dg@LBLONE{bl}\def\dg@LBLTWO{tr}\fi\fi
         \dg@XLBL=\dg@SIZE
         \multiply\dg@XLBL\dg@LBLPOS
         \divide\dg@XLBL\dgARROWPARTS\relax
         \dg@changesign\dg@XLBL\dg@DX
         \dg@YLBL=\dg@XLBL
         \multiply\dg@YLBL\dg@DY \divide\dg@YLBL\dg@DX
         \advance\dg@XLBL\dg@XOFFSET
         \advance\dg@YLBL\dg@YOFFSET
      \fi
   \fi

   \dg@XLBLOFF=-\dg@DY \dg@changesign\dg@XLBLOFF\dg@DX
   \dg@YLBLOFF=\dg@DX \dg@changesign\dg@YLBLOFF\dg@DX
   \ifnum\dg@DX=\z@ \dg@XLBLOFF=\m@ne \fi
   \dg@XTEMP=\dg@DX \dg@changesign\dg@XTEMP\dg@DX
   \dg@YTEMP=\dg@DY \dg@changesign\dg@YTEMP\dg@DY
   \ifnum\dg@YTEMP>\dg@XTEMP \dg@XTEMP=\dg@YTEMP \fi
   \ifnum\dg@XTEMP=\z@ \dg@XTEMP=\@ne \fi
   \multiply\dg@XLBLOFF\dg@LBLOFF \divide\dg@XLBLOFF\dg@XTEMP
   \multiply\dg@YLBLOFF\dg@LBLOFF \divide\dg@YLBLOFF\dg@XTEMP
   \multiply\dg@X\@m \multiply\dg@X\dg@XGRID
   \multiply\dg@Y\@m \multiply\dg@Y\dg@YGRID
   \relax}%

\def\dg@rmcommondiv#1#2#3{%
   \dg@XTEMP=#2\relax
   \divide\dg@XTEMP #1\relax \multiply\dg@XTEMP #1\relax
   \dg@YTEMP=#3\relax
   \divide\dg@YTEMP #1\relax \multiply\dg@YTEMP #1\relax
   \ifnum\dg@XTEMP=#2\relax \ifnum\dg@YTEMP=#3\relax
      \divide#2#1\relax \divide#3#1\relax \fi\fi}%

\def\dg@changesign#1#2{%
   \ifnum #2<\z@ \multiply#1\m@ne
   \else\ifnum #2=\z@ #1=\z@ \fi\fi}%

\def\dg@getnodesize#1#2#3#4#5{%
   #4=\z@\relax #5=\z@\relax
   \expandafter\@for\expandafter\dg@trynode
   \expandafter:\expandafter=#1\do{%
      \dg@XNODE=\m@ne 
      \dg@trynode
      \ifnum #2=\dg@XNODE \ifnum #3=\dg@YNODE
         #4=\dg@XTEMP\relax #5=\dg@YTEMP\relax\fi\fi}}%

\newoptcommand{\dg@makebox}{}[2]{%
   \expandafter\makebox\expandafter(\expandafter
      0\expandafter,\expandafter0\expandafter)\expandafter
      [#1]{#2}}%

\def\dg@novector(#1,#2)#3{}%

\def\dg@letname#1#2{%
   \relax\expandafter
   \let\expandafter #1\csname #2\endcsname\relax}%

\def\dgl@#1{#1{}{}}%
\def\dgl@t#1#2{#1{#2}{}}%
\def\dgl@b#1#2{#1{}{#2}}%
\def\dgl@tb#1#2#3{#1{#2}{#3}}%
\def\dgl@l#1#2{#1{#2}{}}%
\def\dgl@r#1#2{#1{}{#2}}%
\def\dgl@lr#1#2#3{#1{#2}{#3}}%

\makeatother
\newcounter {bean}

\newcounter{butter}
\newcounter{au}
\def\description label#1{\hfil\bf[#1]\hfil} \parskip 5pt plus 1pt
\topmargin
-.5in
\newcommand{\Z}{\bf{Z}\rm} \newcommand{\C}{\cal C}
\newcommand{\R}{\cal R}

\newcommand{\B}{\cal B}

\newcommand{\A}{\cal A}

\newcommand{\M}{\cal M}

\newcommand{\J}{\cal J}
\newcommand{\p}{\cal P}
\newcommand{\Q}{\cal Q}
\newcommand{\doublebar}[1]{\overline{\overline{#1}}}
\newcommand{\doubletilde}[1]{\stackrel{\approx}{#1}}
\begin{document}
\large

\centerline{\bf The Fibers of the Prym Map}

\medskip

\centerline{Ron Donagi}

\bigskip
\bigskip

\centerline{Department of Mathematics}
\centerline{University of Pennsylvania}
\centerline{Philadelphia, PA 19104}

\bigskip

Introduction

\S1.  Pryms

\hspace*{.3in} \S1.1  Pryms and parity

\hspace*{.3in} \S1.2  Bilinear and quadratic forms

\hspace*{.3in} \S1.3  The Prym maps

\S2.  Polygonal constructions

\hspace*{.3in} \S2.1  The $n$-gonal constructions

\hspace*{.3in} \S2.2  Orientation

\hspace*{.3in} \S2.3  The bigonal construction

\hspace*{.3in} \S2.4  The trigonal construction

\hspace*{.3in} \S2.5  The tetragonal construction

\S3.  Biellpitic Pryms

\S4.  Fibers of ${\p}:{\R}_6\to{\A}_5$

\hspace*{.3in} \S4.1  The structure

\hspace*{.3in} \S4.2  Special fibers

\hspace*{.7in}(4.6)  Cubic threefolds

\hspace*{.7in}(4.7)  Jacobians

\hspace*{.7in}(4.8)  Quartic double solids and the branch locus of
${\p}$

\hspace*{.7in}(4.9)  Boundary behavior

\S5.  Fibers of ${\p}:{\R}_5\to{\A}_4$

\hspace*{.3in} \S5.1  The general fiber

\hspace*{.3in} \S5.2  Isotropic subgroups

\hspace*{.3in} \S5.3  Proofs

\hspace*{.3in} \S5.4  Special fibers

\hspace*{.7in}(5.11)  Nodal cubic threefolds
\footnote{

\noindent
Research partially supported by NSF grant No.
DMS90-08113. This paper is in final form and no version of it  will
be
submitted
for publication elsewhere.}

\hspace*{.7in}(5.12)  Degenerations in ${\R\C}^+$

\hspace*{.7in}(5.13)  Degenerations in ${\R}^2{\Q}^+$

\hspace*{.7in}(5.14)  Jacobians

\hspace*{.7in}(5.15)  The boundary

\hspace*{.7in}(5.16)  Thetanulls

\hspace*{.7in}(5.17)  Pentagons and wheels

\S6.  Other genera
\bigskip
\bigskip

\noindent{\Large\bf Introduction}

The Prym map $${\p}:{\R}_g\to{\A}_{g-1}$$ sends a pair
$(C,\widetilde{C})\in{\R}_g$, consisting of a curve $C\in{\M}_g$
and an unramified double cover $\widetilde{C}$, to its Prym variety
$$P={\p}(C,\widetilde{C}):=\ker^0({\rm Nm}:J(\widetilde{C})\to
J(C)).$$Prym
varieties and the Prym map are central to several
approaches to the Schottky problem, e.g. [B1], [D3-D5], [Deb1],
[vG], [vGvdG], [I], [M2], [W].  The purpose of this work is to
describe the fibers of the Prym map.  When $g=5$ or 6, these fibers
turn out to have some beautiful, and perhaps unexpected, structure.
We spend much of our effort in \S\S4, 5 on analyzing the picture in
these two cases, both generically and over some of the natural
special loci in ${\A}_4$ and ${\A}_5$.  In \S6 we summarize some of
what is known in other genera.

Here are some of the results.  When $g=6$, the map is generically
finite of degree 27 [DS].  We show that its monodromy group equals
the Weyl group $WE_6$, and that the general fiber has on it a
structure which is equivalent to the incidence correspondence on
the 27 lines on a non-singular cubic surface (Theorem (4.2)).  The
map fails to be finite over some of the interesting loci in
${\A}_5$, such as ${\J}_5$ (5-dimensional Jacobians) and ${\C}$
(intermediate Jacobians of cubic threefolds).  Finiteness is
restored when ${\p}$ is compactified (\S1.3) and blown up (\S4.2);
the  \linebreak resulting finite fibers can be described
very explicitly ((4.6), (4.7)).  A similarly explicit description
of the fibers is available over the locus of intermediate Jacobians
of Clemens' quartic double solid ((4.8), following [C1], [DS]).The
latter is
the branch locus of ${\p}:{\R}_6\to{\A}_5$ [D6].

When $g=5$, we show (Theorems (5.1)-(5.3)) that the fiber ${\p}^{-
1}(A)$, over generic $A\in{\A}_4$, is a double cover
$\widetilde{F(X)}$ of the Fano surface $F(X)$ of lines on a cubic
threefold $X$.  The correspondence between $A\in{\A}_4$ and the
pair $(X,\delta)\in{\R\C}^+$ consisting of the cubic threefold $X$
and the non zero, ``even", point $\delta$ of order 2 in $J(X)$,
is a birational equivalence of the moduli spaces.  It also turns
out that ${\R}_5$ has an involution $\lambda$ which commutes with
${\p}$, inducing the sheet interchange on the double cover
$\widetilde{F(X)}$.  This $\lambda$ is quite exotic; for example,
it interchanges double covers of trigonal curves with ``Wirtinger"
double covers of nodal curves (5.14).  Again, we can describe the
fiber in more detail over the three distinguished divisors in
$\bar{\A}_4$:  Jacobians, the boundary \linebreak(= degenerate
abelian
varieties), and the locus $\theta_{null}$ of abelian varieties
with a vanishing thetanull:  in all three cases, the cubic
threefold becomes nodal, and the covers $(C,\widetilde{C})$ in the
fiber can be described.  A particularly pretty picture arises for
$\bar{\p}^{-1}(A)$, where $A\in{\A}_4$  is the unique
4-dimensional, non-hyperelliptic PPAV with 10 vanishing
thetanulls.  Varley [V] showed that all Humbert curves (with their
natural double covers) are in this fiber.  We observe that the
corresponding cubic threefold is Segre's 10-nodal cubic (4.8); this
leads quickly to a complete description of the whole fiber, (5.17).

For other values of $g$, the picture does not seem to be quite as
rich.  For $g\le 4$, one can give a rather elementary description
of the fibers using Masiewicki's criterion [Ma] and Recillas'
trigonal construction [R].  When $g\ge 7$, the map is generically
injective ([FS], [K], [W]), but we show that it is never injective
(\S6).

The main tool used to analyze the Prym map is the tetragonal
construction (\S2.5), a triality on the locus of curves with a
$g^1_4$ in ${\R}_g$, which commutes with ${\p}$.  We exploit it
consistently, together with standard facts [ACGH] on the existence
of $g^1_4$'s on curves of low genus, to establish the various
structures on the fibers of ${\p}$.  In genus 5 this fits into a
larger symmetry, indexed by the finite projective plane ${\bf
P}^{2}({\bf F}_2)$, which we describe in \S5.2 and use to find the
cubic threefold.

Almost all the results in this work were announced in [D1].  Since
then, several preliminary manuscripts have circulated, but most of
these results have not been published before.  Several interesting
recent developments concerning closely related questions,
especially Clemens' notes [C2] and Izadi's thesis [I], convinced me
that these ideas may still be useful, and should be published.  The
present work, then, provides the details for almsot everything in
[D1].  The main exception are the results on the Andreotti-Mayer
locus, which have since appeared (in a corrected form) in [Deb1]
and in [D5].  I include here only the underlying idea, which is the
systematic application of the tetragonal construction to double
covers of bielliptic curves (\S3).

As mentioned above, several beautiful extensions of our results
have recently been obtained by Clemens [C2] and Izadi [I].  Their
basic idea is that the cubic threefold $X$ associated to an abelian
variety $A\in{\A}_4$ can be realized concretely inside the van
Geemen-van der Geer linear system $\Gamma_{00}$ on $A$, through
use of Clemens' quartic double solids.  The period map ${\J}$ for
these is analyzed in [D6], and the fiber ${\J}^{-1}(A)$ turns out
to be a certain cover of the cubic threefold $X$.  Clemens
constructs a map
$$c:{\J}^{-1}(A)\to\Gamma_{00}$$
whose image is $X$.  He
conjectures, and Izadi proves, that the projective dual $X^*$ of
$X$ can be recovered as an irreducible component of the branch locus
of the rational map from $A$ to $\Gamma_{00}^{*}$
determined by the linear system $\Gamma_{00}$.  This concrete
model of $X$ leads to several interesting
applictions:

 \noindent $\bullet$  Over $A\in\partial{\A}_5$, which is a ${\bf
C}^*$-extension of $A_0\in{\A}_4$, Izadi obtains the cubic surface
(of Theorem (4.2)) as hyperplane section of the cubic threefold $X$
of $A_0$. (cf. (4.9) for some more details.)

\noindent $\bullet$  The Abel-Prym models of the six genus-5 curves
making up a ${\bf P}^2({\bf F}_2)$-diagram (\S5.2) can be realized
as the intersection $\Theta_a\cap\Theta_{-a}\cap H$ of two
\mbox{theta-translates} with a divisor in $\Gamma_{00}$

\noindent $\bullet$  Izadi is able to describe precisely where our
birational map ${\A}_4\sim{\R\C}^+$ fails to be an isomorphism.

\noindent $\bullet$  She is also able to verify some of the [vGvdG]
conjectures in genus 4.

A second area of current activity is conjecture (6.5.1), which says
that all non-injectivity of the Prym maps is due to the tetragonal
construction.  For non-hyperelliptic, non-trigonal and
non-bielliptic curves of genus $\ge 13$, this was proved in [Deb2].
The generic bielliptic case, $g\ge 10$, is in [N].  Radionov
[Ra] has recently proved that for $g\ge 7$ the graph of the
tetragonal construction provides at least an irreducible component
of the non-injectivity locus of ${\p}$.

Some of the results of the present work were used in [D3] and [D4]
to study the Schottky-Jung loci.  This leads to a proof, which I
hope to publish in the near future, of the Schottky-Jung conjecture
in genus 5, i.e. that the Schottky-Jung equations in genus 5
characterize Jacobians.  An exciting new idea in [vGP] is the
interpretation of Schottky-Jung and tetragonal-type identities
via rank-2 vector bundles; we wonder whether the results on the
geometry of the Prym map will also admit interpretations in terms
of the geometry of the moduli space of vector bundles.

It is a pleasure to acknowledge many beneficial conversations on
the subject of Prym varieties which I have had over the years with
Arnaud Beauville, Roy Smith, Robert Varley, and especially Herb
Clemens, who introduced me to Prym geometry and to his double
solids, and who has had a profound motivating effect on my
thinking.

\bigskip

\noindent{\Large\bf Notation}

\noindent{\bf Moduli Spaces}:

\begin{tabbing}
xxxxxxxxxxx \=\kill
${\M}_g$:  \>curves of genus $g$. \\
$\bar{\M}_g$:\>  the Deligne-Mumford compactification. \\
${\A}_g$:  \>$g$-dimensional principally polarized abelian
varieties (PPAV). \\
${\J}_g$: \> the closure in ${\A}_g$ of the locus of Jacobians. \\
${\Q}\subset{\M}_6$: \> plane quintic curves. \\
${\C}\subset{\A}_5$:  \>(Intermediate Jacobians of) cubic
threefolds. \\
${\R\A}_g$: \> pairs $(A,\mu), \ \ A\in{\A}_g, \ \ \mu\in A_2$ a
non-zero point of order 2. \\
${\R}_g,{\R\Q}, {\R\C}$: \> the pullback of the cover
${\R\A}_g\to{\A}_g$ to ${\M}_g, {\Q,\C}$ respectively. \\
$\doublebar{\R}_g$: \> the Deligne-Mumford compactification of
${\R}_g$.
\\
$\bar{\R}_g$:  \> the open subset of $\doublebar{\R}_g$
ofBeauville-allowable
double covers \\
\>(\S1.3).
\end{tabbing}

\medskip

\noindent{\bf Maps}

\begin{tabbing}
xxxxxxxxxxxxxx \=\kill
${\p} : {\R}_g\to{\A}_{g-1}$:  \>the Prym map. \\
$\bar{\p}:\bar{\R}_g\to{\A}_{g-1}$: \> Beauville's proper version
of ${\p}$. \\ $\doublebar{\p}:\doublebar{\R}_g\to\bar{\A}_{g-1}$:\>
a
compactification of ${\p}$, where $\bar{\A}_{g-1}$ denotes
(Satake's \\ \> compactification, or) an appropriate toroidal
compactifica-\\
\> tion. \\
$\Phi, \Psi, \varphi,\psi$: \> canonical,
Prym canonical, Abel-Jacobi and Abel-Prym \\ \> maps of a curve.
\end{tabbing}

We work throughout over the complex number field ${\bf C}$.

\section{Pryms.}

\subsection{Pryms and parity.}

\ \ \ \ Let $$\pi:\widetilde{C}\rightarrow C$$ be an unramified,
irreducible double cover of a curve $C\in{\cal M}_g$.  The genus of
$\widetilde{C}$ is then $2g-1$, and we have the Jacobians $$J:=J(C),
\
\ \ \ \ \ \ \widetilde{J}:=J(\widetilde{C})$$ of dimensions $g,\,
2g-1$
respectively, and the norm homomorphism
$${\rm Nm}:\widetilde{J}\longrightarrow J.$$ Mumford shows [M2] that
$${\rm Ker}({\rm Nm})=P\cup P^-$$ where $P={\cal P}(C,
\widetilde{C})$ is an
abelian subvariety of $\widetilde{J}$, called the Prym variety, and
$P^-$ is its translate by a point of order 2 in $\widetilde{J}$.  The
principal polarization on $\widetilde{J}$ induces twice a principal
polarization on the Prym.  This appears most naturally when we
consider instead the norm map on line bundles of degree $2g-2$,
$${\rm Nm}:{\rm Pic}^{2g-2}(\widetilde{C})\rightarrow {\rm
Pic}^{2g-2}(C).$$
Let $\omega_{C}\in {\rm Pic}^{2g-2}(C)$ be the canonical bundle of
$C$.
\bigskip

\noindent {\bf Theorem 1.1}  (Mumford [M1], [M2])
\begin{list}{{\bf(\arabic{bean})}}{\usecounter{bean}}
\item The two components $P_0, \,P_1$ of ${\rm Nm}^{-1}(\omega_C)$
can be
distinguished by their parity:  $$P_i=\{L\in
{\rm Nm}^{-1}(\omega_C) \ \ \ | \ \ \ h^0(L)\equiv i \ \ \ \ \ {\rm
mod.
\ 2}\}, \ \ \ \ \ \ \ \ \ \ i=0,1.$$
\item Riemann's theta divisor $\widetilde{\Theta}'\subset {\rm
Pic}^{2g-2}(\widetilde{C})$ satisfies $$\widetilde{\Theta}'\supset
P_1$$
and $$\widetilde{\Theta}'\cap P_0=2\Xi'$$ where $\Xi'\subset P_0$ is
a
divisor in the principal polarization on $P_0$.
\end{list}

\subsection{Bilinear and quadratic forms.}

\ \ \ \ Let $X\in{\cal A}_g$ be a {\rm PPAV}, and let $Y$ be a torser
(=principal homogeneous space) over $X$.  By theta divisor in $Y$
we mean an effective divisor whose translates in $X$ are in the
principal polarization.  $X$ acts by translation on the variety
$Y'$ of theta divisors in $Y$, making $Y'$ also into an $X$-torser.
In $X'$ there is a distinguished divisor $$\Theta':=\{\Theta\subset
X|\Theta\ni 0\}\subset X'$$ which turns out to be a theta divisor,
$\Theta'\in X''$.  In particular, we have a natural identification
$X''\approx X$ sending $\Theta'$ to $0$. Let $X_2$ be the subgroup
of points of order $2$ in $X$.  Inversion on $X$ induces an
involution on $X'$; the invariant subset $X'_2$, consisting of
symmetric theta divisors in $X$, is an $X_2$-torser.  Let
$\langle\,  ,\, \rangle$
denote the natural ${\bf F}_2$-valued (Weil) pairing on $X_2$.  On
$X'_2$ we have an ${\bf F}_2$-valued function
$$q=q_X:X'_2\rightarrow{\bf F}_2$$ sending $\Theta\in X'$ to its
multiplicity at $0\in X$, taken mod. 2.

\bigskip

\noindent {\bf Theorem 1.2}  [M1]  The function $q_X$ is quadratic.
Its associated bilinear form, on $X_2$, is $\langle \, ,\, \rangle$.
When
$(X,\Theta)$ vary in a family, $q_X(\Theta)$ is locally constant.

When $X$ is a Jacobian $J=J(C)$, these objects have the following
interpretations:

\noindent
\begin{tabbing}
$q(L)$ \= $\approx$ \=  $\{ L \in {\rm Pic}^{0}(C) = J \ | \ L^{2} \
\approx
{\cal O}_{C} \} \ \approx \ H^{1}(C,{\bf F}_{2})$ \= \ (semi periods)
\kill
$J'$ \> $\approx$ \> ${\rm Pic}^{g-1}(C) \,$ \ \ \ \ \ \ \ \ \ \ \ \
\ \ \ \ \
\ \ \ \ \ \ \ \ \ \ \ \ \
(use Riemann's theta divisor) \\
$J_{2}$ \> $\approx$ \>  $\{ L \in {\rm Pic}^{0}(C) = J \ | \ L^{2} \
\approx
{\cal O}_{C} \} \ \approx \ H^{1}(C,{\bf F}_{2})$ \=  \ \ \ \  (semi
periods)
\\
$J'_{2}$ \> $\approx$ \>  $\{ L \in {\rm Pic}^{g-1}(C) \ \ \ \ | \
L^{2} \
\approx
{\cal \omega}_{C} \ \} \,$ \ \ \ \ \ \ \ \ \ \ \ \ (theta
characteristics) \\
$q(L)$ \> $\equiv$ \> $h^{0}(C,L)$ mod. 2 \ \ \ \ \ \ \ \ \ \ \ \ \ \
\ \ \ \ \
\ \ \ \ \ \ \ \ \ \ \ \ \ \
(by Riemann-Kempf)
\end{tabbing}

Explicitly, the theorem says in this case that for $\nu, \sigma\in
J_2$ and $L\in J'_2$:

\bigskip

\noindent ({\bf 1.3}) $\langle \nu,\sigma \rangle \ \ \ \equiv \ \ \
h^0(L)+h^0(L\otimes\nu)+h^0(L\otimes\sigma)+h^0(L\otimes\nu
\otimes\sigma)$ \begin{flushright} mod. 2. \end{flushright}

We note that non-zero elements $\mu\in J_2$ correspond exactly to
irreducible double covers $\pi:\widetilde{C}\rightarrow C$.  Let $X$
be
the Prym $P={\cal P}(C,\widetilde{C})$, which we also denote
$P(C,\mu)$, $P(C,\widetilde{C})$, \ $P(\widetilde{C}/C)$ etc.  Now
the
divisor $\Xi'\subset P_0$ of Theorem 1.1 gives
a natural identification $$P'\approx P_0\subset\widetilde{J}'.$$  The
pullback $$\pi^*:J\longrightarrow\widetilde{J}$$ sends $J_2$ to
$\widetilde{J}_2$.  Since ${\rm Nm}\circ\pi^*=2$, we see that
$$\pi^*(J_2)\subset P_2\cup P^-_2.$$  Let $(\mu)^\perp$ denote the
subgroup of $J_2$ perpendicular to $\mu$ with respect to $\langle \,
, \,
\rangle$.

\bigskip

\noindent {\bf Theorem 1.4}
[M2]\begin{list}{{\bf(\arabic{bean})}}{\usecounter{bean}}
\item For $\tau\in J_2$, $\pi^*\tau\in P_2$ iff
$\tau\in(\mu)^{\perp}$.
\item This gives an exact sequence
$$0\rightarrow(\mu)\rightarrow(\mu)^\perp
\stackrel{\pi^*}{\rightarrow}P_2\rightarrow 0.$$
\item In (2), $\pi^*$ is symplectic, i.e. $$\langle \nu,
\sigma \rangle_J \, = \, \langle \pi^*\nu, \pi^*\sigma\rangle _P, \ \
\ \ \ \ \
\ \
\ \nu, \sigma\in(\mu)^\perp\subset J_2.$$
\end{list}

This equality of bilinear forms can be refined to an equality of
quadratic functions.  The identifications $$J'\approx {\rm Pic}^{g-
1}(C), \ \ \ \ \ \ \ \widetilde{J}'\approx {\rm Pic} ^{2g-2}
(\widetilde{C})$$ convert the pullback $$\pi^* :{\rm
Pic}^{g-1}(C)\rightarrow{\rm Pic}^{2g-2}(\widetilde{C})$$ into a map
of
torsers $${\pi^{*}}' : J' \rightarrow \widetilde{J}'$$ over the group
homomorphism $$\pi^{*}:J\rightarrow\widetilde{J}.$$ Let
$${(\mu)^{\perp
}}':=({\pi^{*}}')^{-1}(P'_2).$$ the refinement is:

\bigskip

\noindent {\bf Theorem 1.5}
[D4]\begin{list}{{\bf(\arabic{bean})}}{\usecounter{bean}}
\item ${(\mu)^{\perp}}'$ is contained in $J'_2$ and is a
$(\mu)^\perp$-coset there.
\item ${\pi^{*}}':{(\mu)^{\perp}}'\rightarrow P'_2$ is a map of
torsers over $\pi^*:(\mu)^\perp\rightarrow P_2$.
\item In (2), ${\pi^{*}}'$ is orthogonal, i.e.
$$q_J(\nu)=q_P({\pi^{*}}'\nu), \ \ \ \ \ \ \ \ \ \ \ \ \ \ \ \ \
\nu\in{(\mu)^{\perp}}'.$$
\end{list}

\subsection{The Prym Maps.}

\ \ \ \ Let ${\cal R}_g$ be the moduli space of irreducible double
covers $\pi:\widetilde{C}\rightarrow C$ of non-singular curves
$C\in{\cal M}_g$.  Equivalently, ${\cal R}_g$ parametrizes pairs
$(C, \mu)$ with $\mu\in J_2(C)\backslash(0)$, a semiperiod on
$C$.  The assignment of the Prym variety to a double cover gives a
morphism $${\cal P}:{\cal R}_g\rightarrow{\cal A}_{g-1}.$$
Let $\iota$ be the involution on $\widetilde{C}$ over $C$.  The Abel-
Jacobi map $$\varphi:\widetilde{C}\rightarrow J(\widetilde{C})$$
induces
the Abel-Prym map $$\psi:\widetilde{C}\rightarrow{\rm Ker}({\rm
Nm})$$
$$ \ \ \ \ \ \ \ \ \ \ \ x\longmapsto\varphi(x)-\varphi(\iota x).$$
The image actually lands in the wrong component, $P^-$, but
at least $\psi$ is well-defined up to translation (by a point of
order 2).  In particular, its derivative is well-defined; it
factors through $C$, yielding the Prym-canonical map
$$\Psi:C\rightarrow{\bf P}^{g-2}$$ given by the complete linear
system $|\omega_C\otimes\mu|$. Beauville computed the
codifferential of the Prym map:

\bigskip

\noindent {\bf Theorem 1.6} [B1]  The codifferential $$d{\cal
P}:T^*_P{\cal A}_{g-1}\rightarrow T^*_{(C,\mu)}{\cal R}_g$$ can be
naturally identified with restriction $$\Psi^*:S^2 H^0(\omega_C
\otimes \mu) \rightarrow H^0(\omega^2_C).$$ In particular,
${\rm Ker}(d{\cal P})$ is given by quadrics through the
Prym-canonical
curve $\Psi(C)\subset{\bf P}^{g-2}$.

Let $\bar{\cal A}_g$ denote a toroidal compactification of ${\cal
A}_g$.  Its boundary $\partial\bar{\cal A}_g$ maps to $\bar{\cal
A}_{g-1}$, and the fiber over generic $A\in{\cal
A}_{g-1}\subset\bar{\cal A}_{g-1}$ is the Kummer variety
$K(A):=A/(\pm 1)$.  In codimension 1, this picture is independent
of the toroidal compactification used.

Let ${\cal RA}_g$ denote the level moduli space parametrizing
pairs $(A,\mu)$ with $A\in{\cal A}_g \ , \ \mu\in
A_2\backslash(0)$, and let $\overline{\cal RA}_g$ be a toroidal
compactification.  In [D3] we noted that its boundary has 3
irreducible components, distinguished by the relation of the
vanishing cycle (mod. 2), $\lambda$, to the semiperiod $\mu$:

\noindent{\bf (1.7)} \ \ \ \ \ \ \ \ \ \ \ \ \ \ \ \ \ \
$\begin{array}{llll} \partial^{\rm I} \ \ : \!\!\! & \!\!\!\lambda
\!\!\! & \!\! = \!\! & \!\!\! \mu \\ \partial^{\rm II} \ : \!\!\! &
\!\!\! \lambda \!\!\! & \!\! \neq \!\! & \!\!\! \mu, \ { \ }
\!\langle
\lambda,\mu\rangle =0\in{\bf F}_2 \\ \partial^{\rm III}:\langle
\!\!\! &
\!\!\!\lambda \!\!\! & \! , \!\! & \!\!\! \mu\rangle \neq0.
\end{array} $

Let $\bar{\cal M}_g$, $\doublebar{\cal R}_g$ denote the
Deligne-Mumford
stable-curve compactifications of ${\cal M}_g$ and ${\cal R}_g$.At
least in
codimension one, the Jacobi map extends: $$\bar{\cal
M}_g\rightarrow\bar{\cal A}_g \ , \ \doublebar{\cal R}_g\rightarrow
\overline{\cal RA}_g.$$  We use $\partial\bar{\cal M}_g$ \ , \
$\partial^i\doublebar{\cal R}_g$ \ \ \ $(i = {\rm I, \,  II,  \,
III})$ to
denote
the intersections of $\bar{\cal M}_g$, $\doublebar{\cal R}_g$ with
the
corresponding boundary divisors in $\bar{\cal A}_g$,
$\overline{\cal RA}_g$.

In [B1], Beauville introduced the notion of an allowable double
cover.  This leads to the construction ([DS] I, 1.1) of a proper
version of the Prym map, $$\bar{\cal P}:\bar{\cal R}_g\rightarrow
{\cal A}_{g-1}.$$ Roughly, one extends ${\cal P}$ to
$$\doublebar{\cal P}:\doublebar{\cal R}_g\rightarrow\bar{\cal
A}_{g-1},$$ then restricts to the open subset $\bar{\cal
R}_g\subset\doublebar{\cal R}_g$ of covers which are allowable, in
the sense that their Prym is in ${\cal A}_{g-1}$.This condition can
be made
more explicit:

\bigskip

\noindent {\bf Theorem 1.8} [B1] A stable curve $\widetilde{C}$ with
involution $\iota$, quotient $C$, is allowable if and only if all
the fixed points of $\iota$ are nodes of $\widetilde{C}$ where the
branches are not exchanged, and the number of nodes exchanged under
$\iota$ equals the number of irreducible components exchanged under
$\iota$.

We illustrate the possibilities in codimension 1:
\bigskip

\noindent {\bf Examples 1.9}
\begin{list}{{\rm(\Roman{butter})}}{\usecounter{butter}}
\item $X\in{\cal M}_{g-1}, \ p, q\in X, \ \ p\neq q$; let $X_0,
X_1$ be isomorphic copies of $X$.  Then $C:=X/(p\sim q)$ is a point
of $\partial\bar{\cal M}_g$.  The Wirtinger cover
$$\widetilde{C}:=(X_0\amalg X_1)/(p_0\sim q_1, p_1\sim q_0)$$ gives a
point $$(C,\widetilde{C})\in\partial^{\rm I}\bar{\cal R}_g$$ which is
allowable.  The Prym is $${\cal P}(C,\widetilde{C})\approx
J(X)\in{\cal
A}_{g-1}.$$
\item Start with $(\widetilde{X}\rightarrow X)\in{\cal R}_{g-1}$,
choose distinct points $p, q\in X$, let $p_i, q_i(i=0, 1)$ be their
inverse images in $\widetilde{X}$, and set $$C:=X/(p\sim q), \ \ \
\widetilde{C}:=\widetilde{X}/(p_0\sim q_0, p_1\sim q_1).$$ Then
$$(C,\widetilde{C})\in\partial^{\rm II}\doublebar{\cal R}_g$$ is an
unallowable
cover.  Its Prym is a ${\bf C}^*$-extension of ${\cal
P}(X,\widetilde{X})$; the extension datum defining this extension is
given by $$\psi(p_0)-\psi(q_0)\in{\cal P}(X,\widetilde{X}),$$ which
is
well defined modulo $\pm 1$ (i.e. in the Kummer), as it should be.
\item $X, p, q$ as before, but now $\widetilde{X}\rightarrow X$ is a
double cover branched at $p, q$; consider Beauville's cover
$$C:=X/(p\sim q), \ \ \
\\widetilde{C}:=\widetilde{X}/(\widetilde{p}\sim\widetilde{q})$$
where
$\widetilde{p},
\widetilde{q}$ are the ramification points in $\widetilde{X}$ above
$p, q$.
Then $(C, \widetilde{C})\in\partial^{\rm III}\bar{\cal R}_g$ is
allowable.
\end{list}

In [M1], Mumford lists all covers $(C, \widetilde{C})\in{\cal R}_g$
whose Pryms are in the Andreotti-Mayer locus (i.e. have theta
divisors singular in codimension 4).  A major result in [B1]
(Theorem (4.10)) is the extension of this list to allowable covers
in $\bar{\cal R}_g$.  We do not copy Beauville's list here, but we
will refer to it when needed.

\section{Polygonal constructions}

\subsection{The $n$-gonal constructions}

\ \ \ \ Let $$f:C\rightarrow K$$ be a map of non singular algebraic
curves, of
degree $n$, and $$\pi:\widetilde{C}\rightarrow C$$ a branched double
cover.  These two determine a $2^n$-sheeted branched cover of $K$,
$$f_*\widetilde{C}\rightarrow K,$$ whose fiber over a
general point $k\in K$
consists of the $2^n$ sections $s$ of $\pi$ over $k$: $$s:f^{-
1}(k)\rightarrow\pi^{-1}f^{-1}(k), \ \ \ \pi\circ s=id.$$ The curve
$f_*\widetilde{C}$ can be realized, for instance, as sitting in ${\rm
Pic}^n(\widetilde{C})$ or $S^n\widetilde{C}$:

\noindent           {\bf (2.1)} \ \ \ $f_*\widetilde{C}=\{D\in
S^n\widetilde{C} \ \ | \ \ {\rm Nm}(D)=f^{-1}(k), \ {\rm some} \ k\in
K\}.$

\noindent
(If we think of $\widetilde{C}$ as a local system on an open subset
of  $C$,
this
is just the direct image local system on $K$, hence our notation
$f_*\widetilde{C}$.) On $f_*\widetilde{C}$ we have two structures:
an
involution $$\iota:f_*\widetilde{C}\rightarrow f_*\widetilde{C}$$
obtained
by changing all $n$ choices in the section $s$ via the involution
(also denoted $\iota$) of $\widetilde{C}$, and an equivalence
relation
$$f_*\widetilde{C}\rightarrow\widetilde{K}\rightarrow K$$ where
$\widetilde{K}$
is a branched double cover of $K$: two sections $$s_1, s_2:f^{-
1}(k)\rightarrow\pi^{-1}f^{-1}(k)$$ are equivalent if they differ
by an even number of changes.

For even $n$, the involution $\iota$ respects equivalence, so we
have a sequence of maps

 \noindent{\bf (2.1.1)} \ \ \ \ \ \ \ \ \ \ \ \ \ \ \ \ \ \ \
$f_*\widetilde{C}\rightarrow
f_*\widetilde{C}/\iota\rightarrow\widetilde{K}\rightarrow K$

\noindent of degrees $2, 2^{n-2}, 2$ respectively.  For odd $n$ the
equivalence classes are exchanged by $\iota$, so we have instead a
Cartesian diagram:
\begin{equation}
\renewcommand{\theequation}{\bf
{\arabic{section}}.{\arabic{subsection}}.{\arabic{equation}}}
\setcounter{equation}{2}
\begin{diagram}[f ]
\node[2]{f_{*}\widetilde{C}} \arrow{sw} \arrow{se} \\
\node{f_{*}\widetilde{C}/\iota} \arrow{se} \node[2]{\widetilde{K}}
\arrow{sw}
\\
\node[2]{K}
\end{diagram}
\end{equation}
{\bf Remark 2.1.3} In prctice we will often want to allow $C$ to
acquire some
nodes, over which $\pi$ may be etale (as in (1.9 II)) or ramified (as
in
\linebreak (1.9 III)).
We will always consider this as a limiting case of the non-singular
situation,
and interpret the $n$-gonal construction in the limit so as to make
it depend
continuously on the parameters, whenever possible. We will see
various examples
of this below.

\subsection{Orientation}

\ \ \ \ We observe that the branched cover $\widetilde{K}\rightarrow
K$
depends on $f\circ\pi:\widetilde{C}\rightarrow K$, but not on $f,
\pi$
or $C$ directly.  More generally, to an $m$-sheeted branched cover
$$g:M\rightarrow K$$ we can associate an $m!$-sheeted branched
cover (the Galois closure of $M$) $$g!:M!\rightarrow K,$$ with an
action of the symmetric group $S_m$; the quotient by the
alternating group $A_m$ gives a branched double cover
$$O(g):O(M)\rightarrow K$$ which we call the orientation cover of
$M$.  We say $M$ is orientable (over $K$) if the double cover
$O(M)$ is trivial.  One verifies easily that the double cover
$\widetilde{K}\rightarrow K$
(obtained in \S 2.1 from the maps
$\widetilde{C}\stackrel{\pi}{\rightarrow}C\stackrel{f}{\rightarrow}K$
as quotient of $f_*\widetilde{C}$) is the orientation cover
$O(f\circ\pi)$ of $\widetilde{C}$.

\bigskip
\noindent{\bf Corollary 2.2}  If $\widetilde{C}$ is orientable over
$K$
then $f_*\widetilde{C}=\widetilde{C}_0\cup\widetilde{C}_1$ is
reducible:
\begin{list}{{\rm(\roman{butter})}}{\usecounter{butter}}
\item For $n$ even, the involution $\iota$ acts on each
$\widetilde{C}_i$ with quotient $C_i$ of degree $2^{n-2}$ over $K, \
\
\ i=0, 1$.
\item For $n$ odd, $\iota$ exchanges $\widetilde{C}_0,
\widetilde{C}_1$.Each
$\widetilde{C}_i$ has degree $2^{n-1}$ over $K$.
\end{list}

\bigskip

\noindent{\bf Lemma 2.3}  Branch $(\widetilde{K}/K)=f_*({\rm Branch}
\
(\widetilde{C}/C))$.
This means:  if one point of $f^{-1}(k)$ is a branch point of
$\widetilde{C}\rightarrow C$, then $k$ is a branch point of
$\widetilde{K}\rightarrow K$; if two points of $f^{-1}(k)$ are branch
points of $\widetilde{C}\rightarrow C$, then $k$ is not a branch
point
of (the normalization of) $\widetilde{K}\rightarrow K$, but the two
sheets of $\widetilde{K}$ there intersect; etc.  In particular, the
ramification behavior of $f:C\rightarrow K$ does not affect the
ramification of $\widetilde{K}$.

\bigskip

\noindent{\bf Corollary 2.4}  Let $f:C\rightarrow{\bf P}^1$ be a
branched cover, $\pi:\widetilde{C}\rightarrow C$ an (unramified)
double
cover.  Then $\widetilde{C}$ is orientable over ${\bf P}^1$.

(More generally, the conclusion holds whenever $$f_*({\rm Branch}
(\pi))=2D$$ for some divisor $D$ on ${\bf P}^1$, since the
normalization of $O(\widetilde{C})$ is then an unramified double
cover
of the simply connected ${\bf P}^1$, by (2.3).  In this situation
we say that $\pi$ has \underline{cancelling ramification.})
\bigskip

\noindent{\bf Remark 2.5}  Assume $K={\bf P}^1$ and $\pi$
unramified.  The image of $f_*\widetilde{C}$ in ${\rm
Pic}(\widetilde{C})$
is: $$\{L\in{\rm Pic}^n(\widetilde{C}) \ \ | \ \ {\rm Nm}(L)=f^*{\cal
O}_{{\bf P}^1}(1), \ \ \ h^0(L) > 0\}.$$ This is contained in a
translate of $${\rm Nm}^{-1}(\omega_C)=P_0\cup P_1,$$ and the
splitting
(2.2) of $f_*\widetilde{C}$ is ``explained", in this case, by the
splitting (1.1) of ${\rm Ker}({\rm Nm})$, i.e. after translation:
$$\widetilde{C}_i\subset P_i, \ \ \ \ \ i=0, 1,$$ cf. [D1, \S 6],
[B2].

\bigskip
\noindent{\bf Remark 2.6}  The splitting of $f_*\widetilde{C}$ can
also
be explained group theoretically.  Let $WC_n$ be the group of
signed permutations of $n$ letters, i.e. the subgroup of $S_{2n}$
centralizing a fixed-point-free involution of the $2n$ letters.Let
$WD_n$ be
its subgroup of index 2 consisting of even signed
permutations, i.e. permutations of $n$ letters followed by an even
number of sign changes.  (These are the Weyl groups of the Dynkin
diagrams $C_n,D_n$.)  Over an arbitrary space $X$, we have
equivalences:
\pagebreak[4]

$$\{ \ \ \ \ \ n{\rm -sheeted \ cover \ } Y\rightarrow
X \ \ \ \ \ \} \ \longleftrightarrow\{ \ \ \ {\rm Representation \
} \pi_1(X)\rightarrow \ \ \ \ \ \ S_n \ \ \}$$  $$ \left\{
\begin{array}{l} n{\rm -sheeted \ cover \ } Y\rightarrow X \\ {\rm
with \ a \ double \ cover \ } \widetilde{Y}\rightarrow Y \end{array}
\right\} \longleftrightarrow \left\{ \begin{array}{l} {\rm
 \ Representation \ } \pi_1(X)\rightarrow WC_n \end{array} \right\}
$$ $$ \left\{ \begin{array}{l} n{\rm -sheeted \ cover \ }
Y\rightarrow X \\ {\rm with \ an \ orientable \ double \ \ } \\
{\rm cover \ } \widetilde{Y}\rightarrow Y \ \end{array} \right\}
\longleftrightarrow \left\{ \begin{array}{l} {\rm Representation \
} \pi_1(X)\rightarrow WD_n \end{array} \right\} $$

\medskip

The basic construction of $f_*\widetilde{C}$ then corresponds to the
standard representation $$\rho:WC_n\hookrightarrow S_{2^n}.$$ The
existence of the involution $\iota$ on $f_*\widetilde{C}$ corresponds
to the factoring of $\rho$ through $$WC_{2^{n-1}}\subset S_{2^n}.$$
The restriction $\bar{\rho}$ of $\rho$ to $WD_n$ factors through
$$S_{2^{n-1}}\times S_{2^{n-1}},$$ explaining the splitting when
$\widetilde{C}$ is orientable.

\subsection{The bigonal construction}

\ \ \ \ The case $n=2$ of our construction (``bigonal") takes a
pair of maps of degree 2:
$$\widetilde{C}\stackrel{g}{\rightarrow}C\stackrel{f}{\rightarrow}K$$
and produces another such pair
$$f_*\widetilde{C}\stackrel{g'}{\rightarrow}\widetilde{K}
\stackrel{f'}{\rightarrow} K.$$
Above any given point $k \in K$, the possibilities are:
\begin{list}{(\roman{bean})}{\usecounter{bean}}
\item If $f$, $g$ are etale then so are $f'$ and $g'$.
\item If $f$ is etale and $g$  is branched at one of the two points
$f^{-1}(k)$,
then $f'$ is branched at $k$, and $g'$ is etale there.
\item Vise versa, if $f$ is branched and $g$ is etale then $f'$ is
etale and
$g'$
 is branched at one point of $f'^{-1}(k)$.
 \item If both $f$ and $g$ are branched over $k$ then so are $f'$,
$g'$.
 \item If $f$ is etale and $g$ is branched at both points
$f^{-1}(k)$, then
 $\widetilde{K}$ will have a node over $k$, and $g' :
f_{*}\widetilde{C}
 \rightarrow \widetilde{K}$ will be a $\partial^{\rm III}$
degeneration, i.e.
 will look like (1.9 III).
 \item Vice versa, we can extend the bigonal construction by
continuity, as in
 (2.1.3), to allow $g : \widetilde{C} \rightarrow  C$ to degenerate
to a
 $\partial^{\rm III}$-cover. This leads to $f'$ which is etale and
$g'$ which
 is branched at both points of $f'^{-1}(k)$.
\end{list}

The following general properties are immediately verified:

\bigskip

\noindent{\bf Lemma
2.7}\begin{list}{{\bf(\arabic{bean})}}{\usecounter{bean}}
\item The bigonal construction is symmetric, i.e. if it takes
$\widetilde{C}\stackrel{g}{\rightarrow}C\stackrel{f}{\rightarrow}K$
to
$\widetilde{C}'\stackrel{g'}{\rightarrow}C'\stackrel{f'}{\rightarrow}
K$
then it takes $\widetilde{C}'\rightarrow C'\rightarrow K$ to
$\widetilde{C}\rightarrow C\rightarrow K$.
\item The bigonal construction exchanges branch loci: $${\rm
Branch}(g')=f_*({\rm Branch}(g)), \ \ \ \ \ {\rm
Branch}(f)=g'_*({\rm Branch}(f')).$$
\end{list}
(As in lemma (2.3), this requires the following convention in case
(vi)
above: the local contribution to Branch($f$) is $2k$, and the
contribution
to Branch($g$) is 0).

The symmetry group of this situation, $WC_2$, is the
dihedral group
of the square: $$WC_2= \ \ \langle r, f \ \ \ | \ \ \
f^2=r^4=(rf)^2=1\rangle.$$
($r=90^\circ$ rotation, $f=$flip around $x$-axis, in the
2-dimensional representation.)  It has a non-trivial outer
automorphism (=conjugation by a $45^\circ$ rotation), which
explains why conjugacy classes of representations (of $\pi_1(X)$)
in $WC_2$ come in (bigonally related) pairs.
We list all conjugacy classes of subgroups of $WC_2$ in the
following diagram ($\sim$ denotes conjugate subgroups):

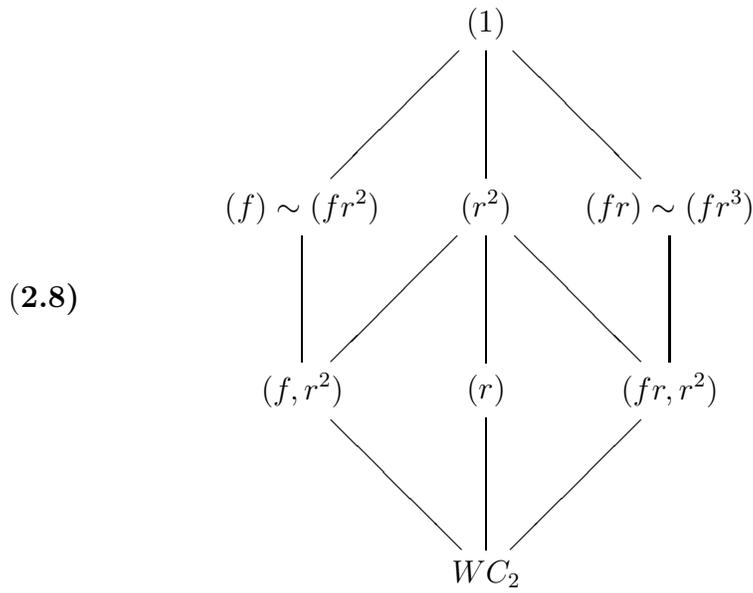
\begin{equation}
\renewcommand{\theequation}{\bf
{\arabic{section}}.{\arabic{equation}}}
\setcounter{equation}{8}
\begin{diagram}[(fr) \sim]
\node[2]{(1)} \arrow{sw,-} \arrow{s,-} \arrow{se,-} \\
\node{(f) \sim (fr^{2})} \arrow{s,-} \node{(r^{2})} \arrow{sw,-}
\arrow{s,-}
\arrow{se,-}
\node{(fr) \sim (fr^{3})} \arrow{s,-} \\
\node{(f,r^{2})} \arrow{se,-} \node{(r)} \arrow{s,-}
\node{(fr,r^{2})}
\arrow{sw,-} \\
\node[2]{WC_{2}}
\end{diagram}
\end{equation}

Correspondingly, we obtain the diagram of curves and maps of degree
2:

\begin{equation}
\renewcommand{\theequation}{\bf
{\arabic{section}}.{\theau}.{\arabic{equation}}}
\setcounter{au}{8}
\setcounter{equation}{1}
\begin{diagram}[\widetilde{C}]
\node[2]{\doubletilde{C}} \arrow{sw} \arrow{s} \arrow{se} \\
\node{\widetilde{C}} \arrow{s} \node{C \times_{K} C'} \arrow{sw}
\arrow{s}
\arrow{se}
\node{\widetilde{C}'} \arrow{s} \\
\node{C} \arrow{se} \node{C''} \arrow{s} \node{C'} \arrow{sw} \\
\node[2]{K}
\end{diagram}
\end{equation}
Here the two sides are bigonally related.

Note that $C'$ is $O(\widetilde{C})$; so if $\widetilde{C}$ is
orientable
(e.g. if $K={\bf P}^1$ and $g$ is unramified) then everything
splits:  $$\widetilde{C}'=C_0\amalg C_1\rightarrow K\amalg K=C',$$
$\widetilde{C}$ is Galois over $K$ with group $({\bf Z}/2{\bf Z})^2$
and quotients
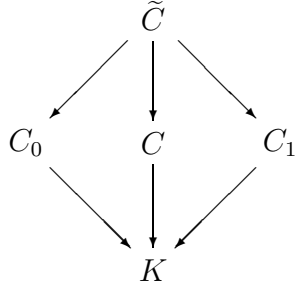
\[
\begin{diagram}[C_{0}]
\node[2]{\widetilde{C}} \arrow{sw} \arrow{s} \arrow{se} \\
\node{C_{0}} \arrow{se} \node{C} \arrow{s} \node{C_{1}} \arrow{sw} \\
\node[2]{K}
\end{diagram}
\]

(cf. [M1]), and (2.8.1) simplifies to:
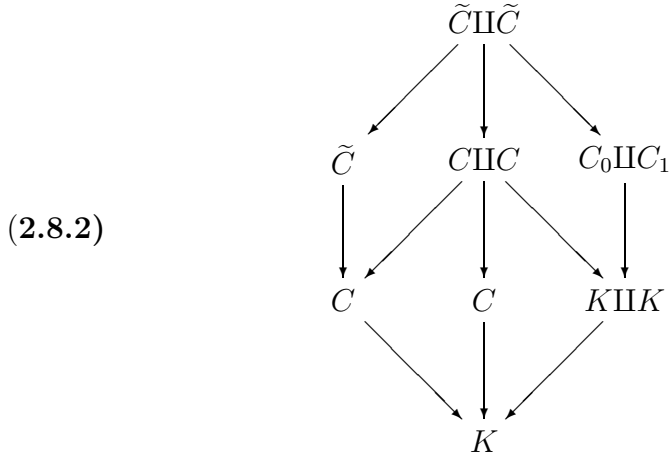
\begin{equation}
\renewcommand{\theequation}{\bf
{\arabic{section}}.{\theau}.{\arabic{equation}}}
\setcounter{au}{8}
\setcounter{equation}{2}
 \begin{diagram}[\widetilde{C} {\textstyle \amalg}]
\node[2]{\widetilde{C} {\textstyle \amalg} \widetilde{C}} \arrow{sw}
\arrow{s}
\arrow{se} \\
\node{\widetilde{C}} \arrow{s} \node{C {\textstyle \amalg} C}
\arrow{sw}
\arrow{s} \arrow{se}
\node{C_{0} {\textstyle \amalg} C_{1}} \arrow{s} \\
\node{C} \arrow{se} \node{C} \arrow{s} \node{K {\textstyle \amalg} K}
\arrow{sw} \\
\node[2]{K}
\end{diagram}
\end{equation}

Given an arbitrary branched double cover $\widetilde{C}\rightarrow
C$,
we form its Prym variety
$$P(\widetilde{C}/C):={\rm Ker}^0({\rm
Nm}:J(\widetilde{C})\rightarrow J(C)).$$
It is
an abelian variety (for $C, \widetilde{C}$ non-singular), but in
general not a principally polarized one.  Nevertheless, there is a
simple relationship between the bigonally-related Pryms
$P(\widetilde{C}/C)$ and $P(\widetilde{C}'/C'):$ in the case $K={\bf
P}^1$,
Pantazis [P] showed that these abelian varieties are dual to each
other.

\subsection{The trigonal construction.}

\ \ \ \ The case $n=3$ of our construction was discovered by
Recillas [R].  Start with a tower
$$\widetilde{C}\stackrel{\pi}{\rightarrow}C\stackrel{f}{\rightarrow}{
\bf P}^1$$ where $f$ has degree 3, and $\widetilde{C}\rightarrow C$
is
an unramified double cover.  By Corollaries (2.4) and (2.2),
$f_*\widetilde{C}$ consists of two copies of a tetragonal curve
$g:X\rightarrow{\bf P}^1$.  Since $f$ and $g$ have the same branch
locus by Lemma (2.3), we find from Hurwitz' formula: $${\rm
genus}(X)={\rm genus}(C)-1.$$  All in all, we have constructed a
map: $$T: \left\{ \begin{array}{c} {\rm trigonal \ curves \ } C
{\rm  \ of \ } \\ {\rm genus \ } g {\rm \ with \ a \ double \ cover
\ } \widetilde{C} \end{array} \right\} \rightarrow \left\{
\begin{array}{c} {\rm tetragonal \ curves \ }  \\ X {\rm \ of \
genus \ } g-1 \end{array} \right\}.$$
We claim that this map is a bijection (except that
sometimes a nonsingular object on one side may correspond to a
singular one on the other):  given $g:X\rightarrow{\bf P}^1$, we
recover $\widetilde{C}$ as the relative second symmetric product of
$X$
over ${\bf P}^1$, $$\widetilde{C}:=S^2_{{\bf P}^1}X\rightarrow{\bf
P}^1,$$ whose fiber over $p\in{\bf P}^1$ consists of all unordered
pairs in $g^{-1}(p)$; this has an involution $\iota$
(=complementation of pairs), giving the quotient
$C:=\widetilde{C}/\iota$.

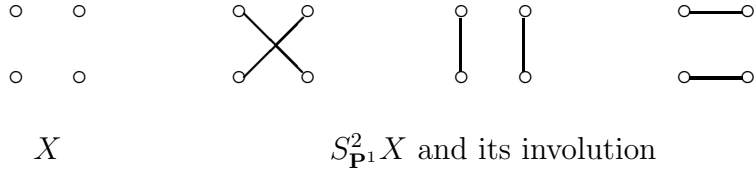
\begin{center}
\begin{tabular}{cccc}
\hspace{1in} & \hspace{1in} & \hspace{1in} & \hspace{1in} \\
\begin{picture}(30,30)(2,1)
\thicklines
\put(2,1){$\circ$} \put(2,25.8){$\circ$}
\put(26,1){$\circ$} \put(26,25.8){${\circ}$}
\end{picture} &
\begin{picture}(30,30)(2,1)
\thicklines
\put(2,1){$\circ$}
\put(7,3.9){\line(1,1){22.4}}
 \put(2,25.8){$\circ$}
 \put(7,27.9){\line(1,-1){22}}
\put(28,1){$\circ$} \put(28,25.8){${\circ}$}
\end{picture} &
\begin{picture}(30,30)(2,1)
\thicklines
\put(2,1){$\circ$} \put(2,25.8){$\circ$}
\put(4.5,6){\line(0,1){20}}
\put(26,1){$\circ$}  \put(28.5,6){\line(0,1){20}}
\put(26,25.8){${\circ}$}
\end{picture} &
 \begin{picture}(30,30)(2,1)
\thicklines
\put(2,1){$\circ$} \put(2,25.8){$\circ$} \put(7,28.3){\line(1,0){20}}
\put(7,3.5){\line(1,0){20}}
\put(26,1){$\circ$}
\put(26,25.8){${\circ}$}
\end{picture}  \\
 &   &  &  \\
$X$ & \multicolumn{3}{c}{$S^{2}_{{\bf P}^{1}}X$ and its involution}
\end{tabular}
\end{center}

In the group-theoretic setup of Remark (2.6), $\bar{\rho}$ induces
an isomorphism $$WD_3\stackrel{\sim}{\rightarrow}S_4.$$ (This is
the standard isomorphism, reflecting the isomorphism of the Dynkin
diagrams $D_3, A_3$.)  Recillas' map $T$ then corresponds to
composition of a representation with this isomorphism.

We list a few of the subgroups of $S_4$:

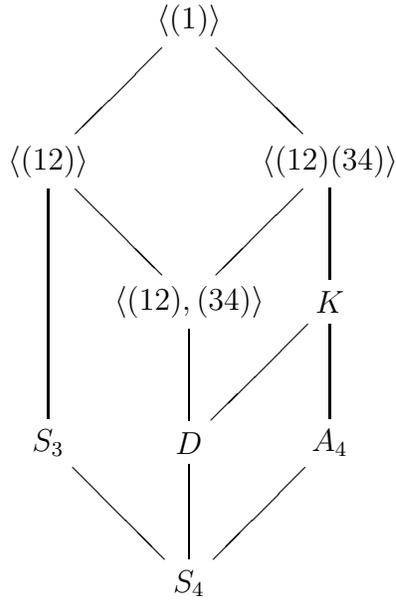
\[
 \begin{diagram}[\widetilde{C} {\textstyle \amalg}]
 \node[2]{\langle (1) \rangle} \arrow{sw,-} \arrow{se,-}  \\
 \node{\langle (12) \rangle} \arrow[2]{s,-} \arrow{se,-}
\node[2]{\langle
(12)(34) \rangle}
 \arrow{sw,-} \arrow{s,-} \\
 \node[2]{\langle (12) , (34) \rangle} \arrow{s,-}  \node{K}
\arrow{s,-}
\arrow{sw,-} \\
 \node{S_{3}} \arrow{se,-} \node{D} \arrow{s,-} \node{A_{4}}
\arrow{sw,-} \\
 \node[2]{S_{4}}
 \end{diagram}
 \]
$D$:  The dihedral group $\langle\!(12), (1324)\!\rangle$ \\
$K=D\cap A_4$:  The Klein group $\langle\!(12)(34),
(13)(24)\!\rangle$.

The corresponding curves are:

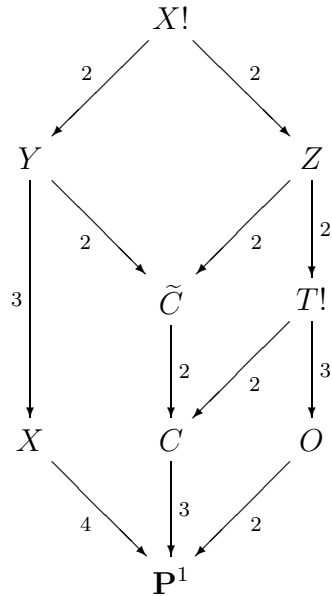
\[
 \begin{diagram}[\widetilde{C} {\textstyle \amalg}]
 \node[2]{X !} \arrow{sw,t}{2} \arrow{se,t}{2}  \\
 \node{Y} \arrow[2]{s,l}{3} \arrow{se,b}{2} \node[2]{Z}
 \arrow{sw,b}{2} \arrow{s,r}{2} \\
 \node[2]{\widetilde{C}} \arrow{s,r}{2}  \node{T !} \arrow{s,r}{3}
\arrow{sw,b}{2} \\
 \node{X} \arrow{se,b}{4} \node{C} \arrow{s,r}{3} \node{O}
\arrow{sw,b}{2} \\
 \node[2]{{\bf P}^{1}}
 \end{diagram}
 \]

\noindent $O\approx O(X)\approx O(C)$:  The orientation \\
$Y\approx(X\times_{{\bf P}^1}X)$ $\backslash$ (diagonal) \\
$Z\approx\widetilde{C}\times_{{\bf P}^1}O$.

Using either of these constructions, we can easily describe the
behavior of $X, C, \widetilde{C}$ around various types of branch
points.  Keeping $X$ non-singular, there are the following five
possible local pictures, cf. [DS, III 1.4].

\begin{center}
 \begin{tabular}{|c|c|c|c|c|c|} \cline{2-6}
 \multicolumn{1}{c|}{ } & (i) & (ii) & (iii) & (iv) & (v) \\
\cline{2-6} \hline
  & \hspace{.76in} & \hspace{.76in} & \hspace{.76in} & \hspace{.76in}
&
 \hspace{.76in} \\
   & \hspace{.76in} & \hspace{.76in} & \hspace{.76in} &
\hspace{.76in} &
 \hspace{.76in} \\
  $X$ & \hspace{.76in} & \hspace{.76in} & \hspace{.76in} &
\hspace{.76in} &
 \hspace{.76in} \\
  & \hspace{.76in} & \hspace{.76in} & \hspace{.76in} & \hspace{.76in}
&
 \hspace{.76in} \\
  & \hspace{.76in} & \hspace{.76in} & \hspace{.76in} & \hspace{.76in}
&
 \hspace{.76in} \\  \hline
 & \hspace{.76in} & \hspace{.76in} & \hspace{.76in} & \hspace{.76in}
&
 \hspace{.76in} \\
& \hspace{.76in} & \hspace{.76in} & \hspace{.76in} & \hspace{.76in} &
 \hspace{.76in} \\
   & \hspace{.76in} & \hspace{.76in} & \hspace{.76in} &
\hspace{.76in} &
 \hspace{.76in} \\
  $\widetilde{C}$ & \hspace{.76in} & \hspace{.76in} & \hspace{.76in}
&
\hspace{.76in} &
 \hspace{.76in} \\
  & \hspace{.76in} & \hspace{.76in} & \hspace{.76in} & \hspace{.76in}
&
 \hspace{.76in} \\
  & \hspace{.76in} & \hspace{.76in} & \hspace{.76in} & \hspace{.76in}
&
 \hspace{.76in} \\
 & \hspace{.76in} & \hspace{.76in} & \hspace{.76in} & \hspace{.76in}
&
 \hspace{.76in} \\   \hline
& \hspace{.76in} & \hspace{.76in} & \hspace{.76in} & \hspace{.76in} &
 \hspace{.76in} \\
   & \hspace{.76in} & \hspace{.76in} & \hspace{.76in} &
\hspace{.76in} &
 \hspace{.76in} \\
  $C$ & \hspace{.76in} & \hspace{.76in} & \hspace{.76in} &
\hspace{.76in} &
 \hspace{.76in} \\
  & \hspace{.76in} & \hspace{.76in} & \hspace{.76in} & \hspace{.76in}
&
 \hspace{.76in} \\
  & \hspace{.76in} & \hspace{.76in} & \hspace{.76in} & \hspace{.76in}
&
 \hspace{.76in} \\  \hline
 \end{tabular}
 \end{center}
 \begin{center}
 \begin{tabular}{cc}
 \multicolumn{2}{c}{\Large\bf Legend} \\
 \begin{tabular}{cl}
 \hspace{.35in} & \hspace{1.2in} \\
 \hspace{.35in} &unramified sheet \\
 \hspace{.35in} & \hspace{1.2in} \\
 \hspace{.35in} & \hspace{1.2in} \\
 \hspace{.35in} & simple ramification \\
 \hspace{.35in} & \hspace{1.2in} \\
 \hspace{.35in} & \hspace{1.2in} \\
 \hspace{.35in} & node (two unramified \\
 \hspace{.35in} & sheets glued together) \\
\hspace{.35in} & \hspace{1.2in} \\
 \hspace{.35in} & two ramified sheets\\
 \hspace{.35in} & glued together  \\
\end{tabular} &
\begin{tabular}{cl}
 \hspace{.35in} & \hspace{1.2in} \\
 \hspace{.35in} & ramification point \\
\hspace{.35in} & of index 2 \\
 \hspace{.35in} & \hspace{1.2in} \\
 \hspace{.35in} & \hspace{1.2in} \\
 \hspace{.35in} & ramification point \\
\hspace{.35in} & of index 3 \\
\hspace{.35in} & \hspace{1.2in} \\
 \hspace{.35in} & \hspace{1.2in} \\
 \hspace{.35in} & glueing of two sheets  \\
\hspace{.35in} & of different \\
 \hspace{.35in} & ramification indices \\
 \end{tabular}
\end{tabular}
\end{center}

\begin{list}{{\rm(\roman{butter})}}{\usecounter{butter}}
\item $f, \pi, g$ are \'{e}tale.\item $f$ and $g$ have
simple ramification points, $\pi$ is \'{e}tale.
\item $f$ and $g$ each have a ramification point of index 2, $\pi$
is \'{e}tale.
\item $g$ has two simple ramification points, $\pi$ is a
Beauville cover: \\ $\bar{f}:N\rightarrow{\bf P}^1$ is trigonal,
with a fiber $\{p, q, r\}; \ \bar{\pi}:\widetilde{N}\rightarrow N$ is
branched at $p, q: \ \bar{\pi}^{-1}(p)=\widetilde{p}, \
\bar{\pi}^{-1}(q)=\widetilde{q}$; \ and we have $C=N/(p\sim q), \
\ \ \widetilde{C}=\widetilde{N}/(\widetilde{p}\sim\widetilde{q})$, \
\ and
$\pi:\widetilde{C}\rightarrow C, \ \ \ f:C\rightarrow{\bf P}^1$ are
induced by $\bar{\pi}, \bar{f}$.
\item $g$ has a ramification point of index 3, $\pi$ is Beauville,
$f$ is ramified at one of the two branches of the node of $C$.
\end{list}

\bigskip

Considering first the first three cases, then all five, we
conclude:

\noindent{\bf Theorem 2.9}  The trigonal construction gives
isomorphisms $$T^0:{\cal R}^{{\rm
Trig}}_g\stackrel{\sim}{\rightarrow}{\cal M}^{{\rm Tet}, 0}_{g-1}$$
and $$T:\bar{\cal R}^{{\rm
Trig}}_g\stackrel{\sim}{\rightarrow}{\cal M}^{{\rm Tet}}_{g-1},$$
where:

${\cal M}^{{\rm Tet}}_{g-1}$ is the moduli space of (non-singular)
curves of genus $g-1$ with a tetragonal line bundle.

${\cal M}^{{\rm Tet},0}_{g-1}$ is the open subset of tetragonal
curves $X$ with the property that above each point of ${\bf P}^1$
there is at least one etale point of $X$.

${\cal R}^{{\rm Trig}}_g$ is the moduli space of etale double
covers of non-singular curves of genus $g$ with a trigonal  bundle.

$\bar{\cal R}^{{\rm Trig}}_g$ is the partial compactification of
${\cal R}^{{\rm Trig}}_g$ using allowable covers in $\bar{\cal
R}_g$ of type $\partial^{\rm III}$ (cf (1.9.III)).

\bigskip

\noindent{\bf Examples 2.10}
\begin{list}{{\rm(\roman{butter})}}{\usecounter{butter}}
\item $\widetilde{C}$ is the trivial cover, $\widetilde{C}=C_0\amalg
C_1$,
iff $X$ is disconnected, \\ $X={\bf P}^1\amalg C$, with
$f=g|_C, \ \ \ id_{{\bf P}^1}=g|_{{\bf P}^1}$.
\item Wirtinger covers $(C_0\amalg C_1) \ / \ (p_0\sim q_1, \
q_0\sim p_1)\rightarrow C/(p\sim q)$, where $\{p, q, r\}$
form a trigonal fiber in $C$, correspond to reducible $X={\bf
P}^1\cup_rC$, the
two components meeting at $r\in C$.
\item $C$ is reducible:  $C=H\cup{\bf P}^1$, with $H$
hyperelliptic, and \\ $\widetilde{C}=\widetilde{H}\cup{\bf P}^1$ with
$\widetilde{H}\rightarrow H$ and $\widetilde{{\bf
P}^1}\rightarrow{\bf P}^1$
branched over \\ $B:=H\cap{\bf P}^1$.  This situation
corresponds to $g:X\rightarrow{\bf P}^1$ factoring through a
hyperelliptic $H'$.  Indeed, such a pair $(C, \widetilde{C})$ is
uniquely determined by the tower $\widetilde{H}\rightarrow
H\rightarrow{\bf P}^1$.  The trigonal construction for $C$ is
reduced to the bigonal construction for $H$, which then gives
$X=\widetilde{H}'\rightarrow H'\rightarrow{\bf P}^1$.  In particular:
\item $C=H\amalg{\bf P}^1$ is disconnected iff $X=H_0\amalg H_1$ is
disconnected with hyperelliptic pieces, and then
$\widetilde{C}=\widetilde{H}\amalg{\bf P}^1\amalg{\bf P}^1$, where
$\widetilde{H}$ is the Cartesian cover:
$$\widetilde{H}=H_0\times_{{\bf
P}^1}H_1.$$
\end{list}

So far, we have only used the fact that $\widetilde{C}$ is an
orientable double cover of a triple cover.  We now use our two
assumptions, that $\pi$ is unramified and that the base $K$ equals
${\bf P}^1$, to obtain an identity of abelian varieties.  Namely,
by Remark 2.5 we have a map, natural up to translation.
$$\alpha:X\rightarrow P(\widetilde{C}/C).$$
The result, due to S. Recillas, is:

\bigskip

\noindent{\bf Theorem 2.11}  [R]  If $X$ is trigonally related to
$(\widetilde{C},C)$, then the above map $\alpha$ induces an
isomorphism
$$\alpha_*:J(X)\stackrel{\sim}{\rightarrow}P(\widetilde{C}/C).$$

\noindent{\bf Proof.}

By naturality of $\alpha$ and irreducibility of ${\cal M}^{{\rm
Tet}}_{g-1}$, it suffices to prove this for any one convenient $X$.
We take $\widetilde{C}\rightarrow C$ to be a Wirtinger cover as in
(2.10)(ii), so $$X={\bf P}^1\cup_rC'.$$ where $p+q+r$ is a trigonal
divisor on $C'$, and $C=C'/(p\sim q)$.  We have natural
identifications:  $$J(X)\approx J(C')\approx P(\widetilde{C}/C),$$ in
terms of which $\alpha$ becomes the Abel-Jacobi map $\varphi$ on
$C'$, and collapses ${\bf P}^1$ to a point.The induced $\alpha_*$ is
therefore
the identity.
\begin{flushright} QED \end{flushright}
\noindent{\bf Corollary 2.12}  All trigonal Pryms are Jacobians,
and all tetragonal Jacobians are Pryms.

\subsection{The tetragonal construction}

 \ \ \ \ Consider now a tower $$\widetilde{C}\rightarrow
C\stackrel{f}{\rightarrow}{\bf P}^1$$ where $f$ has degree 4 and
$\widetilde{C}$ is a double cover (unramified) of $C$.  The general
construction yields a sequence of maps of degrees  2, 4, 2:
$$f_*\widetilde{C}\rightarrow
f_*\widetilde{C}/\iota\rightarrow\widetilde{\bf
P}^1\rightarrow{\bf P}^1.$$By (2.2) and (2.4) again, $\widetilde{\bf
P}^1$ is
unramified, hence we have
splittings:
\[\begin{array}{rcl}
\widetilde{\bf P}^1 & = & {\bf P}_0^1\amalg{\bf P}_1^1 \\
f_*\widetilde{C} & = & \widetilde{C}_0\amalg\widetilde{C}_1 \\
f_*\widetilde{C}/\iota & =  & {C}_0\amalg {C}_1.
\end{array}
\]      The tetragonal
construction thus associates to a tower
$$\widetilde{C}\stackrel{2}{\rightarrow}C\stackrel{4}{\rightarrow}{\b
f
P}^1$$ two other towers $$\widetilde{C}_i\rightarrow
C_i\rightarrow{\bf P}^1, \
\ \ \ \ \ \ \ \ \ \ \ i=0, 1$$ of the same type.

\bigskip
\noindent{\bf Lemma 2.13}  The tetragonal construction is a
triality, i.e. starting with \ \ \ $\widetilde{C}_0\rightarrow
C_0\rightarrow{\bf P}^1$ \ \ \ it returns \ \ \
$\widetilde{C}\rightarrow C\rightarrow{\bf P}^1$ \ \ \ and \\
$\widetilde{C}_1\rightarrow
C_1\rightarrow{\bf P}^1$.
\ \ \ \ On the group level, the point is this:  Our tower
$\widetilde{C}\rightarrow C\rightarrow{\bf P}^1$ corresponds to a
representation (of $\pi_1({\bf P}^1\backslash$ (branch locus)))
in $WD_4$.  Now the Dynkin diagram $D_4$ has an automorphism of
order 3:

\vspace{3in}

\noindent This corresponds to an outer automorphism of $WD_4$, of
order 3.  Hence representations in $WD_4$ come in packets of three.
The various groups involved are described in some detail in the
proof of Lemma (5.5), below.

\bigskip
\pagebreak[4]

\noindent{\bf Local pictures 2.14}  Given the local behavior of $C$
and $\widetilde{C}$ over a point of ${\bf P}^1$, it is quite
straightforward to compute $f_*\widetilde{C}$ and hence
$\widetilde{C}_i,
C_i \ \ (i=0, 1)$ over the same point.  Since these local pictures
are needed quite frequently, we record the simplest ones here.

\begin{list}{{\bf(\arabic{bean})}}{\usecounter{bean}}
\item $C, \widetilde{C}$ unramified $\Rightarrow C_i,
\widetilde{C}_i$ are
also unramified.
\item $C$ has one simple ramification point (and two unramified
sheets), $\widetilde{C}\rightarrow C$ unramified $\Rightarrow C_i,
\widetilde{C}_i$ have the same local picture as $C, \widetilde{C}$
respectively.
\item $C$ has two distinct simple ramification points,
$\widetilde{C}\rightarrow C$ unramified $\Rightarrow$ One pair, say
$C_0, \widetilde{C}_0$, has the same local pictures as $C,
\widetilde{C}$,
while the other is a Beauville degeneration:  $C_1$ is unramified
but two of its four sheets are glued, $\widetilde{C}_1\rightarrow
C_1$
is ramified over these two sheets (and the ramification points are
glued) while the other sheets are unramified.
\item $C$ is unramified but two of its sheets are glued,
$\widetilde{C}\rightarrow C$ is ramified over these two sheets
$\Rightarrow$ $C_i$ has two distinct ramification points,
$\widetilde{C}_i\rightarrow C_i$ is unramified $(i=0, 1)$.  (This is
the same triple as in (3).)
\item $C$ has a simple ramification point and the other two sheets
are glued, $\widetilde{C}$ is ramified over the glued sheets
$\Rightarrow C_i, \widetilde{C}_i$ have the same local pictures as
$C,
\widetilde{C}$.
\item $C$ has a ramification point of index 2 (i.e. 3 of its sheets
are permuted by the local monodromy), $\widetilde{C}\rightarrow C$
unramified $\Rightarrow$ same local picture for
$\widetilde{C}_i\rightarrow C_i$.
\item $C$ has a ramification point of index 3 (all 4 sheets
permuted), $\widetilde{C}\rightarrow C$ unramified $\Rightarrow C_0,
\widetilde{C}_0$ have the same local picture as $C, \widetilde{C}$,
but
$C_1$ has a simple ramification point glued to an unramified point,
so $\widetilde{C}_1$ must be simply ramified over each.  (I.e.
$\widetilde{C}_1$ has a point which is simply ramified over ${\bf
P}^1$, glued to a point which has ramification index 3 over ${\bf
P}^1$!)
\end{list}

We note that in examples (3) and (7), the tetragonal construction
applied to $(\widetilde{C}\rightarrow C)\in{\cal RM}_g$ produces an
(allowable) degenerate cover, $(\widetilde{C}_1\rightarrow C_1)\in
\partial^{\rm III}({\cal RM}_g)$.

\begin{center}
 \begin{tabular}{ccc}
 \begin{tabular}{|c|c|c|} \hline
 \hspace{.32in} & \hspace{.32in} & \hspace{.32in} \\
\hspace{.32in} & \hspace{.32in} & \hspace{.32in} \\
 \hspace{.32in} & \hspace{.32in} & \hspace{.32in} \\
 \hspace{.32in} & \hspace{.32in} & \hspace{.32in} \\
 \hspace{.32in} & \hspace{.32in} & \hspace{.32in} \\
 \hspace{.32in} & \hspace{.32in} & \hspace{.32in} \\ \hline
 \hspace{.32in} & \hspace{.32in} & \hspace{.32in} \\
\hspace{.32in} & \hspace{.32in} & \hspace{.32in} \\
 \hspace{.32in} & \hspace{.32in} & \hspace{.32in} \\
 \hspace{.32in} & \hspace{.32in} & \hspace{.32in} \\ \hline
\hspace{.32in} & \hspace{.32in} & \hspace{.32in} \\
 \hspace{.32in} & \hspace{.32in} & \hspace{.32in} \\ \hline
\end{tabular} &
\begin{tabular}{|c|c|c|} \hline
 \hspace{.32in} & \hspace{.32in} & \hspace{.32in} \\
\hspace{.32in} & \hspace{.32in} & \hspace{.32in} \\
 \hspace{.32in} & \hspace{.32in} & \hspace{.32in} \\
 \hspace{.32in} & \hspace{.32in} & \hspace{.32in} \\
 \hspace{.32in} & \hspace{.32in} & \hspace{.32in} \\
 \hspace{.32in} & \hspace{.32in} & \hspace{.32in} \\ \hline
 \hspace{.32in} & \hspace{.32in} & \hspace{.32in} \\
\hspace{.32in} & \hspace{.32in} & \hspace{.32in} \\
 \hspace{.32in} & \hspace{.32in} & \hspace{.32in} \\
 \hspace{.32in} & \hspace{.32in} & \hspace{.32in} \\ \hline
\hspace{.32in} & \hspace{.32in} & \hspace{.32in} \\
 \hspace{.32in} & \hspace{.32in} & \hspace{.32in} \\ \hline
 \end{tabular} &
\begin{tabular}{|c|c|c|} \hline
 \hspace{.32in} & \hspace{.32in} & \hspace{.32in} \\
\hspace{.32in} & \hspace{.32in} & \hspace{.32in} \\
 \hspace{.32in} & \hspace{.32in} & \hspace{.32in} \\
 \hspace{.32in} & \hspace{.32in} & \hspace{.32in} \\
 \hspace{.32in} & \hspace{.32in} & \hspace{.32in} \\
 \hspace{.32in} & \hspace{.32in} & \hspace{.32in} \\ \hline
 \hspace{.32in} & \hspace{.32in} & \hspace{.32in} \\
\hspace{.32in} & \hspace{.32in} & \hspace{.32in} \\
 \hspace{.32in} & \hspace{.32in} & \hspace{.32in} \\
 \hspace{.32in} & \hspace{.32in} & \hspace{.32in} \\ \hline
\hspace{.32in} & \hspace{.32in} & \hspace{.32in} \\
 \hspace{.32in} & \hspace{.32in} & \hspace{.32in} \\ \hline
 \end{tabular} \\
 (1) & (2) & (3,4)
 \end{tabular}
 \end{center}

 \begin{center}
 \begin{tabular}{ccc}
 \begin{tabular}{|c|c|c|} \hline
  \hspace{.32in} & \hspace{.32in} & \hspace{.32in} \\
 \hspace{.32in} & \hspace{.32in} & \hspace{.32in} \\
 \hspace{.32in} & \hspace{.32in} & \hspace{.32in} \\
 \hspace{.32in} & \hspace{.32in} & \hspace{.32in} \\
 \hspace{.32in} & \hspace{.32in} & \hspace{.32in} \\
 \hspace{.32in} & \hspace{.32in} & \hspace{.32in} \\ \hline
  \hspace{.32in} & \hspace{.32in} & \hspace{.32in} \\
  \hspace{.32in} & \hspace{.32in} & \hspace{.32in} \\
 \hspace{.32in} & \hspace{.32in} & \hspace{.32in} \\
 \hspace{.32in} & \hspace{.32in} & \hspace{.32in} \\ \hline
\hspace{.32in} & \hspace{.32in} & \hspace{.32in} \\
 \hspace{.32in} & \hspace{.32in} & \hspace{.32in} \\ \hline
 \end{tabular} &
 \begin{tabular}{|c|c|c|} \hline
 \hspace{.32in} & \hspace{.32in} & \hspace{.32in} \\
\hspace{.32in} & \hspace{.32in} & \hspace{.32in} \\

 \hspace{.32in} & \hspace{.32in} & \hspace{.32in} \\
 \hspace{.32in} & \hspace{.32in} & \hspace{.32in} \\
 \hspace{.32in} & \hspace{.32in} & \hspace{.32in} \\
 \hspace{.32in} & \hspace{.32in} & \hspace{.32in} \\ \hline
 \hspace{.32in} & \hspace{.32in} & \hspace{.32in} \\
\hspace{.32in} & \hspace{.32in} & \hspace{.32in} \\
 \hspace{.32in} & \hspace{.32in} & \hspace{.32in} \\
 \hspace{.32in} & \hspace{.32in} & \hspace{.32in} \\ \hline
\hspace{.32in} & \hspace{.32in} & \hspace{.32in} \\
 \hspace{.32in} & \hspace{.32in} & \hspace{.32in} \\ \hline
 \end{tabular}  &
  \begin{tabular}{|c|c|c|} \hline
 \hspace{.32in} & \hspace{.32in} & \hspace{.32in} \\
\hspace{.32in} & \hspace{.32in} & \hspace{.32in} \\
 \hspace{.32in} & \hspace{.32in} & \hspace{.32in} \\
 \hspace{.32in} & \hspace{.32in} & \hspace{.32in} \\
 \hspace{.32in} & \hspace{.32in} & \hspace{.32in} \\
 \hspace{.32in} & \hspace{.32in} & \hspace{.32in} \\ \hline
 \hspace{.32in} & \hspace{.32in} & \hspace{.32in} \\
\hspace{.32in} & \hspace{.32in} & \hspace{.32in} \\
 \hspace{.32in} & \hspace{.32in} & \hspace{.32in} \\
 \hspace{.32in} & \hspace{.32in} & \hspace{.32in} \\ \hline
\hspace{.32in} & \hspace{.32in} & \hspace{.32in} \\
 \hspace{.32in} & \hspace{.32in} & \hspace{.32in} \\ \hline
 \end{tabular} \\
(5) & (6) & (7)
\end{tabular}
\end{center}

\begin{center}
\begin{tabular}{c}
 \begin{tabular}{|c|c|c|} \hline
\hspace{.32in} & \hspace{.32in} & \hspace{.32in} \\
\hspace{.32in}& \hspace{.32in} & \hspace{.32in} \\
 \hspace{.32in} & \hspace{.32in} & \hspace{.32in} \\
 $\widetilde{C}$ & $\widetilde{C}_{0}$ & $\widetilde{C}_{1}$ \\
 \hspace{.32in} & \hspace{.32in} & \hspace{.32in} \\
 \hspace{.32in} & \hspace{.32in} & \hspace{.32in} \\
 \hspace{.32in} & \hspace{.32in} & \hspace{.32in} \\ \hline
 \hspace{.32in} & \hspace{.32in} & \hspace{.32in} \\
 $C$ & $C_{0}$ & $C_{1}$ \\
 \hspace{.32in} & \hspace{.32in} & \hspace{.32in} \\ \hline
\hspace{.32in} & \hspace{.32in} & \hspace{.32in} \\
 ${\bf P}^{1}$ & ${\bf P}^{1}$ & ${\bf P}^{1}$ \\ \hline
 \end{tabular} \\
 (pattern)
 \end{tabular}
 \end{center}

\pagebreak[4]

\noindent{\bf Examples 2.15}
\begin{list}{{\bf(\arabic{bean})}}{\usecounter{bean}}
\item It is perhaps not terribly surprising that the trigonal
construction is a degenerate case of the tetragonal
construction.Start with
\linebreak $\widetilde{C}\rightarrow C$ the split double
cover of the curve $C$ with the tetragonal map
$f:C\stackrel{4}{\rightarrow}{\bf P}^1$.  Then
$f_*\widetilde{C}$ splits into 5 components, of degrees 1, 4, 6, 4, 1
respectively over ${\bf P}^1$.  The components of degree 4 make up
$\widetilde{C}_1\rightarrow C_1$, which is isomorphic to
$\widetilde{C}\rightarrow C$.  The remaining components give
$${\bf
P}^1\amalg\widetilde{T}\amalg{\bf P}^1\rightarrow T\amalg{\bf P}^1$$
where $(\widetilde{T}, T)$ is associated to $C$ by the trigonal
construction.  Vice versa, starting with an (unramified) double
cover \\ ${\bf P}^1\amalg\widetilde{T}\amalg{\bf P}^1$ of
$T\amalg{\bf
P}^1$, the tetragonal construction produces \linebreak $C\amalg
C\rightarrow
C$, twice.
\item Let $p+q+r+s$ be a tetragonal divisor on $C$.  Then $C/(p\sim
q)$ is still tetragonal.  Tacking a node onto the previous example,
we see that the Wirtinger cover $$(C'\amalg C'')/(p'\sim q'',
q'\sim p'')\rightarrow C/(p\sim q)$$ is taken by the tetragonal
construction to :

 \noindent $\bullet$ Another Wirtinger Cover, $$(C'\amalg
C'')/(r'\sim s'', s'\sim r'')\rightarrow C/(r\sim s),$$ and to:

\noindent $\bullet$ ${\bf P}^1\cup_{t'}\widetilde{T}\cup_{t''}{\bf
P}^1\rightarrow T\cup_t{\bf P}^1$, where $(\widetilde{T}, T)$ is
associated by the trigonal construction to $C$.  (Each copy of
${\bf P}^1$ meets $\widetilde{T}$ or $T$ in the unique point
indicated.
$t\in T$ corresponds to the partition $\{\{p, q\}, \{r, s\}\}$.)
\item We will see in Lemma (3.5) that if $C\rightarrow{\bf P}^1$
factors through a hyperelliptic curve, so do $C_0, C_1$.  An
interesting subcase occurs when $C=H^0\cup H^1$ has two
hyperelliptic components, cf. Proposition (3.6).
\item Let $X$ be a non-singular cubic hypersurface in ${\bf P}^4$,
$\ell\subset X$ a line, and $\widetilde{X}$ the blowup of $X$ along
$\ell$, with projection from $\ell$:
$$\pi:\widetilde{X}\rightarrow{\bf
P}^2.$$ This is a conic bundle [CG] whose discriminant is a plane
quintic curve $Q\subset{\bf P}^2$.  The set of lines $\ell'\subset
X$ meeting $\ell$ is a double cover $\widetilde{Q}$ of $Q$.  Now
choose
a plane $A\subset{\bf P}^4$ meeting $X$ in 3 lines $\ell, \ell',
\ell''$; we get 3 plane quinties $Q, Q', Q''$, with double covers
$\widetilde{Q}, \widetilde{Q}', \widetilde{Q}''$.  Note that $\ell',
\ell''$
map to a point $p\in Q$, hence determine a tetragonal map
$f:Q\rightarrow{\bf P}^1$, given by ${\cal O}_Q(1)(-p)$, and
similarly for $Q', Q''$.  Our observation is that the 3 objects
$$(\widetilde{Q}, Q, f) \ \ ; \ \ (\widetilde{Q}', Q', f') \ \ ; \ \
(\widetilde{Q}'', Q'', f'')$$ are tetragonally related.  Indeed, the
3
maps can be realized simultaneously via the pencil of hyperplanes
$S\subset{\bf P}^4$ containing $A$.  Such an $S$ meets $X$ in a
(generally non-singular) cubic surface $Y$.  $A$ line in $Y$ (and
not in $A$) which meets $\ell'$, also meets 4 of the 8 lines (in
$Y$, not in $A$) meeting $\ell$, one in each of 4 coplanar pairs.this
gives the
desired injection $\widetilde{Q}'\hookrightarrow
f_*\widetilde{Q}$.
\end{list}
Our main interest in the tetragonal construction stems from:

\bigskip

\noindent{\bf Theorem 2.16}  The tetragonal construction commutes
with the Prym map, $$P(\widetilde{C}/C)\approx
P(\widetilde{C}_0/C_0)\approx P(\widetilde{C}_1/C_1).$$

\bigskip

\noindent{\bf Proof}

\ \ \ \ As in Remark (2.5), we have a map
$$\alpha:\widetilde{C}_i\hookrightarrow
f_*\widetilde{C}\rightarrow{\rm
Pic}(\widetilde{C}), \ \ \ \ \ \ \ \ \ \ i=0, 1.$$ The image sits in
a

translate of $P(\widetilde{C}/C)$, so we get induced maps
$$\alpha_*:J(\widetilde{C}_i)\rightarrow P(\widetilde{C}/ C)$$ and
by restriction $$\beta:P(\widetilde{C}_i/ C_i)\rightarrow
P(\widetilde{C}/
C).$$  By Masiewicki's criterion [Ma], $\beta$ will
be an isomorphism if we can show:
\begin{list}{{\rm(\arabic{bean})}}{\usecounter{bean}}
\item The image $\alpha(\widetilde{C}_i)$ of $\widetilde{C}_i$ in
$P(\widetilde{C}/C)$ is symmetric;
\item The fundamental class in $P(\widetilde{C}/C)$ of
$\alpha(\widetilde{C}_i)$ is twice the minimal class, $\frac{2}{(g-
1)!}[\Theta]^{g-1}$.
\end{list}

Now (1) is clear, since the involution on $\widetilde{C}_i$
commutes with $-1$ in $P(\widetilde{C}/C)$.  The fundamental class in
(2) can be computed directly, as is done in [B2].  Instead, we find
it here by a degeneration argument:  it varies continuously with
$(C, \widetilde{C})\in{\cal RM}^{{\rm Tet}}_g$, which is an
irreducible parameter space, so it suffices to do the computation
for a single $(C, \widetilde{C})$.  We take  $$C=T\cup_t{\bf P}^1, \
\
\ \widetilde{C}={\bf P}^1\cup_{t'}\widetilde{T}\cup_{t''}{\bf P}^1,$$
as in
Example (2.15)(2).  Then $(C_i, \widetilde{C}_i)$ is a Wirtinger
cover,
 $i=0, 1$, and the normalization of $C_i$ is the tetragonal curve
$N$ associated to $(T, \widetilde{T})$ by the trigonal construction.
We have identifications $$J(N)\approx P(\widetilde{T}/T)\approx
P(\widetilde{C}/C)$$ (Theorem (2.11)), in terms of which
$\alpha(\widetilde{C}_i)$ consists of the Abel-Jacobi image
$\varphi(N)\subset J(N)$ and of its image under the involution.Thus
the
fundamental class is twice that of $\varphi(N)$, as
required.

(Note:  since this argurment works for any double cover
$\widetilde{T}\rightarrow T$, and since any semiperiod on a nearby
\mbox{non-singular} $C$ specializes to a semiperiod on $T\cup_t{\bf
P}^1$
which is supported on $T$, we need only the irreducibility of
${\cal M}^{{\rm Tet}}_g$, instead of ${\cal RM}^{{\rm Tet}}_g$.)
\begin{flushright} QED \end{flushright}

\large
\section{Bielliptic Pryms.}

\ \ \ \ As a first application of the tetragonal construction, we
show that some remarkable coincidences occur among the various loci
in Beauville's list [B1].  The central role here is played byPryms of
bielliptic curves.  We see in (3.7), (3.8) that      the bielliptic
loci
can be tetragonally related to various other components in
Beauville's
list, and therefore give the same Pryms.  As suggested
in [D1], this leads to a complete, short list of the irreducible
components of the Andreotti-Mayer locus in genus $\leq 5$, and
of its intersection with the image of the proper Prym map for
arbitrary $g$.  We do not include here the complete analysis of the
Andreotti-Mayer locus itself, since this has already appeared in
[Deb1] and [D5] (together with some corrections to the original
list in [D1]).  Nevertheless, we could not resist describing
explicitly the operation of the tetragonal construction on
Beauville's list, as it is such a pretty and straightforward
application of the results of \S 2.

\ \ \ \ We recall Mumford's results on hyperelliptic Pryms.  Let
$$f_i:C^i\rightarrow K, \ \ \ \ \ \ i= 0, 1$$ be two ramified
double covers of a curve $K$.  The fiber product
$$\widetilde{C}:=C^0\times_KC^1$$ has 3 natural
involutions:$\tau_i(i=0, 1)$,
with quotient $C^i$, and \\
$\tau:=\tau_0\circ\tau_1$,  with a new quotient, $C$.  This all
fits in a Cartesian diagram:
      \[
      \begin{diagram}[C]
\node[2]{\widetilde{C}} \arrow{sw,t}{\pi_{1}} \arrow{s,r}{\pi}
\arrow{se,t}{\pi_{0}} \\
\node{C^{0}} \arrow{se,b}{f_{0}} \node{C} \arrow{s,r}{f} \node{C^{1}}
\arrow{sw,b}{f_{1}} \\
\node[2]{K}
\end{diagram}
\]
If the branch loci of $f_0, f_1$ are disjoint, then
$$\pi:\widetilde{C}\rightarrow C$$ is unramified.  We say that a
double
cover obtained this way is \linebreak \underline{Cartesian}.

\bigskip

\noindent{\bf Lemma 3.1}  Let $f:C\rightarrow K$ be a ramified
double cover.  A double cover $$\pi:\widetilde{C}\rightarrow C,$$
given
by a semiperiod $\eta\in J_2(C)$, is Cartesian if and only if \\
$f_*(\eta)=0\in J_2(K)$.

\bigskip

\noindent{\bf Proof:}  apply the bigonal construction.
\begin{flushright} QED \end{flushright}

\bigskip

\noindent{\bf Proposition 3.2} [M1]
\begin{list}{{\bf(\arabic{bean})}}{\usecounter{bean}}
\item Any double cover $\widetilde{C}$ of a hyperelliptic $C$ is
Cartesian.
\item Any hyperelliptic Prym is a product of 2 hyperelliptic
Jacobians (one of which may vanish):  If $\widetilde{C}$ arises as
$C^0\times_{{\bf P}^1}C^1$ then $$P(\widetilde{C}/C)\approx
J(C^0)\times J(C^1).$$
\end{list}

\bigskip

\noindent{\bf Proof:}  (2) follows from (1), (1) follows from lemma
(3.1) with $K = {\bf P}^{1}$.  \begin{flushright} QED
\end{flushright}

\ \ \ \ A \underline{bielliptic} curve (aka elliptic-hyperelliptic,
superelliptic, ...) is a branched double cover of an elliptic
curve.  In this section we apply the tetragonal construction to
find various identities between bielliptic Pryms and Pryms of
other, usually degenerate, curves.  Some of the results extend
to \underline{bihyperelliptic} curves, i.e. branched double
covers of hyperelliptic curves.  To warm up, we consider
\underline{Jacobians} of bihyperelliptic curves.  Example
(2.10)(iii) can be restated:

\bigskip

{\bf Lemma 3.3} The trigonal construction gives a bijection between:

\noindent
$\bullet \ \ $ Bihyperelliptic, non singular curves $C$:
\[ C \stackrel{f}{\rightarrow} H \stackrel{g}{\rightarrow} {\bf
P}^{1};      \]
$\bullet \ \ $ Reducible trigonal double covers $\widetilde{X}
\rightarrow X$:
\begin{center}
\begin{tabular}{rlccc}
$\begin{diagram}[X]
\node{\widetilde{X}} \arrow{s} \\
\node{X}
\end{diagram}$ &
$\begin{diagram}[X]
\node{=} \\
\node{=}
\end{diagram}$ &
$\begin{diagram}[X]
\node{C'} \arrow{s}  \\
\node{H'}
\end{diagram}$  &
$ \begin{diagram}[X]
\node{\cup}  \\
\node{\cup}
\end{diagram} $&
$\begin{diagram}[X]
\node{H} \arrow{s} \\
\node{{\bf P}^{1}}
\end{diagram}$
\end{tabular}
\end{center}
 where

 $X = H' \cup {\bf P}^{1}$ is reducible

 $\tau : X \rightarrow {\bf P}^{1}$,  the trigonal map, has
 degree 2 on $H'$ and 1 \linebreak on ${\bf P}^{1}$.

 $\tau(H' \cap {\bf P}^{1})$ = Branch($g$)

 $\widetilde{X} \rightarrow X$ is allowable of type $\partial^{\rm
III}$
 at each point of $H' \cap {\bf P}^{1}$.
 \bigskip

 \noindent
 We note that $C' \rightarrow H' \rightarrow {\bf P}^{1}$ is
bigonally
 related to $C \rightarrow H \rightarrow {\bf P}^{1}$.

\bigskip

\noindent{\bf Corollary 3.4}  The Jacobian of a bihyperelliptic
curve $C$,
$$C\stackrel{f}{\rightarrow}H\stackrel{g}{\rightarrow}{\bf P}^1,$$
is isogenous to the product $$J(H)\times P(g_*C, \iota)$$ of a
hyperelliptic Jacobian and a bihyperelliptic (branched) Prym.

\bigskip

\ \ \ \ We move to the Pryms of bihyperelliptic curves.  First we
note that this class is closed under the tetragonal construction:

\bigskip

\noindent{\bf Lemma 3.5} \ \ Let $(\widetilde{C}_i, C_i)$ be
tetragonally related to $(\widetilde{C}, C)$, with $C$ non-singular.
  If $C\rightarrow{\bf
P}^1$ factors through a (possibly reducible) hyperelliptic $H$, so
do the $C_i$:
$$C_i\stackrel{f_i}{\rightarrow}H_i\stackrel{g_i}{\rightarrow}{\bf
P}^1, \ \ \ \ \ \ i=0, 1.$$

\bigskip
\noindent{\bf Proof.}

The bigonal construction applied to
$$\widetilde{C}\stackrel{\pi}{\rightarrow}C\stackrel{f}{\rightarrow}H
$$
yields $$f_*\widetilde{C}\rightarrow\widetilde{H}\rightarrow H,$$ and
when
applied again to $$\widetilde{H}\rightarrow
H\stackrel{g}{\rightarrow}{\bf P}^1$$ yields
$$g_*\widetilde{H}\rightarrow\widetilde{{\bf P}}^1\rightarrow{\bf
P}^1.$$
Since $\pi$ is unramified, so are $\widetilde{H}\rightarrow H$ and
$\widetilde{\bf P}^1\rightarrow{\bf P}^1$.  Hence $\widetilde{\bf
P}^1$
splits:  $$\widetilde{\bf P}^1={\bf P}^1_0\amalg{\bf P}^1_1,$$ and
this
splitting climbs its way up the tower:
\[
\begin{array}{rlc}
(g\circ f)_{*}\widetilde{C} & = & \widetilde{C}_{0} {\textstyle
\amalg}
\widetilde{C}_{1} \\
& & \downarrow \\
  &  & C_{0} {\textstyle \amalg} C_{1} \\
  & & \downarrow \\
  g_{*}\widetilde{H} & = &  {H}_{0} {\textstyle \amalg} {H}_{1} \\
  & & \downarrow \\
  \widetilde{\bf P}^{1} & = & {\bf P}_{0}^{1} {\textstyle \amalg}
{\bf
P}_{1}^{1} \\
  & & \downarrow \\
    &  &  {\bf P}^{1}
  \end{array}
  \]
\begin{flushright} QED \end{flushright}

\bigskip

\noindent
{\bf Remark 3.5.1}  The rational map $f_{i} : C_{i} \rightarrow
H_{i}$ can,
in a couple of cases, fail to be a morphism; this is easily remedied
by
identifying a pair of points in $H_{i}$. Among the local pictures
(2.14), the
ones that
can occur here are (1), (2), (7) and (3) :
\begin{itemize}
\item In cases (1), (2), the hyperelliptic maps $g$,  $g_{0}$,
$g_{1}$ are all
unramified at the relevant point, and the $f_{i}$ are morphisms.
\item In case (7), $g$ and $g_{0}$ are ramified, $g_{1}$ is not, $f$
and
$f_{0}$
are (ramified) morphisms, but $f_{1}$ is not, since $C_{1}$ is
singular above a
point where $H_{1}$, as constructed above, is nonsingular. To make
$f_{1}$ into
a morphism, we must glue the two points of $g_{1}^{-1}(k)$.
\item In case (3) we find two possibilities:
\begin{list}{(3\alph{bean})}{\usecounter{bean}}
\item $g$ is etale, $f$ is ramified at both points of $g^{-1}(k)$;
then
$g_{0}$,
$g_{1}$ are also etale, $f_{0}$ is ramified at both points of
$g_{0}^{-1}(k)$,
$C_{1}$ has a node but $f_{1}$ is still a morphism.
\item $g$ is ramified, $f$ is etale; then $g_{0}$ is ramified,
$f_{0}$ is
etale, $g_{1}$
is etale, but the two branches of the node of $C_{1}$ are sent by
$f_{1}$ to
opposite sheets of $H_{1}$, so $f_{1}$ is again not a morphism.
\end{list}
\end{itemize}

\noindent{\bf Proposition 3.6}  Let $\widetilde{C}\rightarrow C$ be a
Cartesian double cover of a bihyperelliptic $C$:
$$C\stackrel{f}{\rightarrow}H\stackrel{g}{\rightarrow}{\bf P}^1, \;
\widetilde{C} = C^{0}\times_{H} C^{1}.$$
The tetragonal construction applied to $\widetilde{C}\rightarrow
C\rightarrow{\bf P}^1$ yields:

\noindent $\bullet$ A similar Cartesian tower
$\widetilde{C}_0\rightarrow
C_0\stackrel{f_0}{\longrightarrow}H\stackrel{g_0}{\longrightarrow
}{\bf P}^1$, same $H$.

\noindent $\bullet$ A tower $\widetilde{C}_1\rightarrow
C_1\rightarrow{\bf P}^1$ where:
\begin{tabbing} XXXXX \= \kill
\>$C_1$ is reducible, $C_1=H^0\cup H^1$, \\ \> $H^0, H^1$ are
hyperelliptic, \\ \> $H^0\cap H^1$ maps onto $B:={\rm
Branch}(g)\subset{\bf P}^1$, \\ \>
$\widetilde{C}_1=\widetilde{H}^0\cup\widetilde{H}^1$ is allowable
over
$C_1$, \\
\> $C^{i} \rightarrow H \rightarrow {\bf P}^{1}$ is bigonally related
to
$\widetilde{H}^{i} \rightarrow H^{i} \rightarrow {\bf P}^{1}$, $i$ =
1,2.
\end{tabbing}

Vice versa, the tetragonal construction takes any tower  \linebreak
$\widetilde{C}_1\rightarrow C_1\rightarrow{\bf P}^1$ as above to two
Cartesian bihyperelliptic towers $$\widetilde{C}\rightarrow
C\rightarrow H\rightarrow{\bf P}^1 \ \ \ {\rm and} \ \ \
\widetilde{C}_0\rightarrow C_0\rightarrow H\rightarrow{\bf P}^1.$$

The proof is quite straightforward, and we will simply
write down a few of the relationships involved, using the notation of
the
previous
proof:

\noindent $\bullet$  $\widetilde{H}$ splits into two copies of $H$,
by
(3.1).Hence:

\noindent $\bullet$  $g_*\widetilde{H}\approx H\cup {\bf
P}^1\cup {\bf P}^1$, say $H_0\approx H, \ \ H_1 = R^{0} \cup R^{1},
\;
R^{i} \approx {\bf P}^{1}$, $i$ = 0, 1.

\begin{tabbing}
\noindent $\bullet$  \= Let $H^{i}$, $\widetilde{H}^{i}$ be the
inverse image
of
$R^{i}$ in $C_{1}$, $\widetilde{C}_{1}$ respectively. Then \\
\> $\widetilde{H}^{i} \rightarrow H^{i} \rightarrow {\bf P}^{1}$ is
bigonally
related to $C^{i} \rightarrow H \rightarrow {\bf P}^{1}$.
\end{tabbing}
      \begin{tabbing}
\noindent $\bullet$  \=The intersection properties of the $H^i$ (or
$\widetilde{H}^i$) can be read off the \\ \>local pictures (2.14.3).
\end{tabbing}

\begin{tabbing}
\noindent $\bullet$  \=Finally, let $\varepsilon:H\to H$ be the
hyperelliptic involution.  A cover \\ \>$C^1\to H$ determines a
mirror-image $\varepsilon^*C^1$.  The remaining tower \\
\>$\widetilde{C}_0\to C_0\to H\to{\bf P}^1$ is given by the Cartesian
diagram:
\end{tabbing}
 \[
\begin{diagram}[C]
\node[2]{\widetilde{C}_{0}} \arrow{sw} \arrow{s} \arrow{se} \\
\node{C^{0}} \arrow{se} \node{C_{0}} \arrow{s}
\node{\varepsilon^{*}C^{1}}
\arrow{sw} \\
\node[2]{H}
\end{diagram}
\]
     \begin{flushright} QED \end{flushright}

{\bf Remarks}

\noindent
{\bf (3.6.1)} Since the branch points of $C^{i} \rightarrow H$ map to
the
branch points of $H^{i} \rightarrow {\bf P}^{1}$, we have the
relation
between the genera:
\[ g(H^{i}) = g(C^{i}) - 2\cdot g(H). \]

\noindent
{\bf (3.6.2)} The possible local pictures are exactly the same as in
(3.5.1).
(The use of $C_{0}$, $C_{1}$ in (3.6) is consistent with that of
(2.14).)

\noindent
{\bf (3.6.3)} Another way of proving both lemma (3.5) and proposition
(3.6)
is based on lemma (5.5), which says that the three tetragonal curves
$C, C_{0}, C_{1}$  which are tetragonally related are obtained, via
the
trigonal construction, from one and the same trigonal curve $X$ (with
three distinct double covers). Lemma (3.3) characterizes the possible
curves
$X$, hence proves that the locus of bihyperelliptics is closed under
the
tetragonal construction, lemma (3.5). To complete the proof of
proposition (3.6), one simply needs to characterize the double covers
$\widetilde{X}$ which correspond to Cartesian covers of $C$.

\ \ \ \ For the rest of this section, we specialize to the case
where the hyperelliptic $H$ is an elliptic curve $E$, i.e. $C$ is
bielliptic. First, we write out explicitly the content of
Proposition (3.6) in this case:

\noindent{\bf Corollary 3.7}  The Pryms of double covers $\pi
:\widetilde{C}\to C$ where

\noindent $\bullet$  $C$ is bielliptic,
$C\stackrel{f}{\to}E\stackrel{g}{\to}{\bf P}^1$,

\noindent $\bullet$  $\widetilde{C}\to C$ is Cartesian,
$\widetilde{C}=C^0\times_EC^1, \ \ \ C^0$ is of genus $n$,

\noindent are precisely (via the tetragonal construction) the Pryms
of the following allowable double covers $\widetilde{X}\to X$:

\begin{tabbing}
$n=1$: \ \ \=$X$ is obtained from a hyperelliptic curve by
identifying \\ \>two pairs of points, $X=H/(p\sim q, \ \ r\sim s)$.
\\ \\

$n=2$:\>$X=X_0\cup X_1, \ \ X_0$ rational, $X_1$ hyperelliptic,
\\ \>$\#(X_0\cap X_1)\!=\!4$. \\ \\

$n\ge 3$:\>$X=X_0\cup X_1$, \ each $X_i$ hyperelliptic, \
$g(X_0)=n-2$, \\ \>$g(X_1)=g(C)-n-1, \ \ \#(X_0\cap X_1)=4$,
and both \\ \>hyperelliptic maps are restrictions of the same
tetragonal \\ \>map on $X$ (i.e. they agree on $X_0\cap X_1$).
\end{tabbing}

\ \ \ \ Everything here follows directly from the proposition,
except that for $n=1$ we need to use (twice) the following
observation of Beauville.  Let $\pi:\widetilde{X}\to X$ be an
allowable
double cover where $$X=Y\cup R, \ \ \ \ \ \ \ \ R
\;\mbox{rational}, \ \ \ \ \ \ \ \ \ Y\cap R=\{a,b\}$$
$$\widetilde{X}=\widetilde{Y}\cup\widetilde{R}, \ \
\widetilde{R}=\pi^{-1}(R)\;\mbox{rational}, \ \
\widetilde{Y}\cap\widetilde{R}=\{\widetilde{a},\widetilde{b}\}$$ and
$\pi$ is
ramified at $\widetilde{a}, \widetilde{b}$, which map to $a,b$.
Construct
a new cover $\widetilde{Z}\to Z$ where
$$\widetilde{Z}:=\widetilde{Y}/(\widetilde{a}\sim\widetilde{b})$$
$$Z:=Y/(a\sim
b).$$  Then this is still allowable, and $$P(\widetilde{Z}/Z)\approx
P(\widetilde{X}/X).$$  (Indeed, there are natural isomorphisms of
generalized Jacobians $$J(\widetilde{Z})\approx J(\widetilde{X}),
\;\;J(Z)\approx J(X)$$  commuting with $\pi_*$ and inducing the
desired isomorphisms.)
\begin{flushright} QED \end{flushright}

\ \ \ \ We are left with the Pryms of non-Cartesian double covers
of bielliptic curves.  The result here may be somewhat surprising:

\bigskip

\noindent{\bf Proposition 3.8}  Pryms of \underline{non}-Cartesian
double covers of bielliptic curves are precisely the Pryms of
Cartesian covers (of bielliptic curves) with $n(:=g(C_0))=1.$  (The
isomorphism is obtained through a sequence of 2 tetragonal moves.)

\ \ \ \ The point is that if $X=H/(p\sim q, \ \ r\sim s)$ with $H$
hyperelliptic, and $\widetilde{X}\to X$ is an allowable double cover,
then $P(\widetilde{X}/X)$ is the Prym of a Cartesian cover (with
$n=1$)
of a bielliptic curve, as we've just seen; but $X$ has another
$g^1_4$, and applying the tetragonal construction to it yields a
non-Cartesian double cover of a bielliptic curve.

\ \ \ \ The $g^1_4$ is obtained as follows:  map $H$ to a conic in
${\bf P}^2$ (by the hyperelliptic map), then project the conic to
${\bf P}^1$ from the unique point $x$ in ${\bf P}^2$ (and not on
the conic) on the intersection of the lines $\overline{pq}$ and
$\overline{rs}$.

\vspace{2in}

\ \ \ \ We should now check that the tetragonal construction yields
a non-Cartesian cover of a bielliptic curve, and that all covers
arise this way.  We leave the former to the reader, and do the
latter.

\ \ \ \ Let $\widetilde{C}\to C$ be a non-Cartesian cover of $C$,
which
is bielliptic:  $$ C\stackrel{f}{\to}E\stackrel{g}{\to}{\bf P}^1.$$

\noindent Let $(\widetilde{C}_i,C_i)$, \ \ $i=0,1$, be the
tetragonally
related covers.  By lemma (3.5), $C_i$ is
bihyperelliptic:$$C_i\stackrel{f_i}{\to}H_i\stackrel{g_i}{\to}{\bf
P}^1.$$ By
the
local pictures (2.14), $$B :=\mbox{Branch}(g)=B_0\amalg B_1, \ \
\  B_i:=\mbox{Branch}(g_i).$$
(As we saw in Remark (3.5.1), the possible pictures are (1), (2),
(7),
(3a) and (3b). Of these, (7) and (3b) contribute to $B$, and each
contributes also to one of the $B_{i}$.)
Since $\#B=4$ \ \ ($E$ is elliptic),
and $\# B_i$ is even and $>0$ (non-Cartesian!), we find $$\# B_i=2,
\;\;\; i=0,1,$$ hence $H_i$ is rational and $C_i$ is hyperelliptic.
Again by the local pictures, $C_i$ will have two nodes, at points
lying over $B_{1-i}$.
\begin{flushright} QED \end{flushright}

\ \ \ \ We observe that the last argument works not only for
bielliptics but also for branched double covers of hyperelliptic
curves of
genus 2, since now  $$\#B_0>0, \ \ \#B_1>0, \ \ \#B_0+\#B_1=6, \ \
\#B_i\; \mbox{even}\Rightarrow$$  $$\mbox{either} \ \;\#B_0=2 \
\mbox{or} \ \#B_1=2.$$  However, the resulting hyperelliptic curve
with 4 nodes does not carry other $g^1_4$'s and is not necessarily
related to any other covers.

\ \ \ \ We leave one more corollary of proposition (3.6) to the
reader.

\bigskip

\noindent{\bf Corollary 3.9}  Let $K$ be hyperelliptic,
$\widetilde{K}\to K$ a double cover with 2 branch points.  Then
$P(\widetilde{K}/K)$ is a hyperelliptic Jacobian.

\noindent(Hint:  take both $H$ and $C^0$ in proposition (3.6) to
be rational, show $P(\widetilde{C}_1/C_1)\approx J(C^1)$ and
$C_1=K\cup_{(2\; \mbox{points})}{\bf P}^1, \ K$ hyperelliptic.)
\large
\section{Fibers of ${\cal P}:{\R}_6\to{\A}_5.$}

\subsection{The structure}

\ \ \ \ We recall the main result of [DS]:

\bigskip

\noindent{\bf Theorem 4.1 [DS]}  ${\cal P}:{\R}_6\to{\A}_5$ is
generically finite, of degree 27.

Recall that ${\M}_6^{\rm Tet}$ denotes the moduli space of curves of
genus 6 with a $g^1_4$.  The forgetful map ${\M}_6^{\rm
Tet}\to{\M}_6$
is generically finite, of degree 5 [ACGH].  By base change we get
a corresponding object ${\R}_6^{\rm Tet}$, with map
$${\R}_6^{\rm Tet}\to{\R}_6$$  of degree 5.  The tetragonal
construction gives a triality, or (2,2)-corres\-pon\-den\-ce, on
${\R}_6^{\rm Tet}$.  The image in ${\R}_6$ is then a
(10,10)-correspondence:
\medskip

\noindent{\bf (4.1.1)} \ \ \ \ \ \ \ \ \ \ \ \ \ \ \ \ \ \ \ \ \ \ \
${\rm Tet}\subset{\R}_6\times{\R}_6.$

\bigskip

\noindent{\bf Theorem 4.2}  The correspondence Tet induced by the
tetragonal construction on the fiber ${\cal P}^{-1}(A)$, for
generic $A\in{\A}_4$, is isomorphic to the incidence correspondence
on the 27 lines on a non-singular cubic surface.  The monodromy
group of ${\R}_6$ over ${\A}_5$ (i.e. the Galois group of its
Galois closure) is the Weyl group $WE_6$, the symmetry group of the
incidence of the 27 lines on the cubic surface.

This was conjectured in [DS] and announced in [D1].  The proof will
be given below.  For the symmetry group of the line incidence on a
cubic surface, or other del Pezzo surfaces, we refer to [Dem].

\bigskip

\noindent{\bf (4.3) The blownup map}
Let ${\cal Q}\subset{\M}_6$ denote the moduli space of non-singular
plane
quintic curves,  ${\R\Q}$ its inverse image in ${\R}_6$.  By
Theorem (1.2), it splits:
$${\R\Q}={\R\Q}^+ \ \ \cup \ \ {\R\Q}^-$$ with $(Q,\mu)\in{\R\Q}^+$
(respectively, ${\R\Q}^-$) iff $h^0(\mu\otimes{\cal O}_Q (1)$) is
even (respectively, odd).  The point is that ${\cal O}_Q(1)$ gives
a uniform choice of theta characteristics over ${\Q}$, hence the
spaces of theta characteristics and semiperiods over ${\Q}$ are
identified.

Let ${\cal J}$ be the closure in ${\A}_5$ of the locus of Jacobians
of
curves, and let ${\C}$ denote the moduli space of non-singular
cubic threefolds.  Via the intermediate Jacobian map, we identify
${\C}$ with its image in ${\A}_5$.

The Prym map sends ${\R\Q}^+$ to ${\cal J}$ and ${\R\Q}^-$ to
${\C}$.  Since the fiber dimensions can be positive, it is useful
to consider the blownup Prym map $$\widetilde{\cal
P}:\widetilde{\R}_6\to\widetilde{\A}_5$$  where ${\cal J},{\C}$ on
the
right are blown up to divisors $\widetilde{\cal J},\widetilde{\C}$,
while
on the left we blow up ${\R\Q}^+,{\R\Q}^-$, as well as the locus
${\R}_6^{Trig}$ of double covers of trigonal curves.  The result is
a morphism which is generically finite over $\widetilde{\cal J}$
and
$\widetilde{\C}$.  We recall the geometric description of points of
the
various loci, and
give the map in these geometric terms.  This is taken from [CG],
[T] and [DS].

\bigskip

\noindent{\bf (4.3.1)} \ A point of ${\C}$ is given by a non-singular
cubic threefold $X\subset{\bf P}^4$.  A point of $\widetilde{\C}$
is
given by a pair $(X,H), \ H\in({\bf P}^4)^*$ a hyperplane.

\bigskip

\noindent{\bf (4.3.2)} \ A point of $\widetilde{{\R\Q}}$ is given by
($Q,
\mu, L$), or $(Q, \widetilde{Q}, L)$, where $Q\subset{\bf P}^2$ is
a
plane quintic, $L\in({\bf P}^2)^*$ a line, and $\mu$ a semiperiod
on $Q$ (or $\widetilde{Q}$ the corresponding double cover).

\bigskip

\noindent{\bf (4.3.3)} \ The fiber ${\cal
P}^{-1}(J(X))\subset{\R\Q}^-$ over a cubic threefold $X$ can be
identified with the Fano surface $F(X)$ of lines $\ell\subset
X$.(Projection
from $\ell$ puts a conic bundle structure $\pi:  X
--\!\to{\bf P}^2={\bf P}^4/\ell$ on $X$; the corresponding point of
${\R\Q}^-$ is $(Q,\widetilde{Q})$, where the plane quintic $Q$ is
the
discriminant locus of $\pi$, and its double cover $\widetilde{Q}$
parametrizes lines
$\ell'\in F(X)$ meeting $\ell$.)

\bigskip

\noindent{\bf (4.3.4)} \ The fiber $\widetilde{\cal P}^{-1}(X,H)$
corresponds to the lines $\ell$ in the cubic surface $X\cap H$.For
general
$X,H$, there are 27 of these.  The corresponding
objects are of the form $(Q,\widetilde{Q},L)$ where
$(Q,\widetilde{Q})$ are
as above, and $L\subset{\bf P}^2$ is the projection, $L=\pi(H)$.

\bigskip

\noindent{\bf (4.3.5)} \ A point of ${\R}_6^{Trig}$ is given by a
curve $T\in{\M}_6$ with a trigonal line bundle ${\cal
L}\in\mbox{Pic}^3(T),h^0({\cal L})=2$, and a double cover
$\widetilde{T}\to T$.  The fiber of $\widetilde{\R}_6^{Trig}$ above
it is
given by the linear system $|\omega_T\otimes{\cal L}^{-2}|$, a
${\bf P}^1$.

\bigskip

\noindent{\bf (4.3.6)} \ A point of ${\cal J}$ is given by the
Jacobian of a curve $C\in{\M}_5$.  The canonical curve
$\Phi(C)\subset{\bf P}^4$, for general $C$, is the base locus of a
net of quadrics:  $$A_p\subset{\bf P}^4, \;\;\;\;\;\;p\in{\bf
P}^2={\bf P}^2(C).$$  A point of $\widetilde{\cal J}$ above $C$ is
then
given by a pair $(C,L)$, where $L$ is a line in ${\bf
P}^2(C)$.(Choosing such a
line is the same as choosing a quartic del Pezzo
surface $$S=S_L=\cap_{p\in L}A_p$$ containing $\Phi(C)$.)

\bigskip

\noindent{\bf (4.3.7)} \ Consider the map
$$\alpha:{\M}_5\rightarrow{\R\Q}^+$$ sending $C\in{\M}_5$ to
$\alpha(C)=(Q, \widetilde{Q})$, where: $$Q:=\{p\in{\bf P}^2(C) \ |
\ A_p {\rm \ is \ singular}\}\subset{\bf P}^2(C),$$ and
$\widetilde{Q}$ is the double cover whose fiber over a general
$p\in Q$ corresponds to the two rulings on the rank-4 quadric
$A_p$.  This $\alpha$ is a birational isomorphism; its inverse is
the restriction to ${\R\Q}^+$ of ${\p}$.

The fiber $\widetilde{\p}^{-1}(C, L)$ over generic $(C,
L)\in\widetilde{\cal J}$ is given by the following 27 objects:

\noindent$\bullet$  The quintic
object$(Q,\widetilde{Q},L)\in\widetilde{{\R\Q}}^+$,
where$(Q,\widetilde{Q})=\alpha(C)$ and $L$ is the given line in ${\bf
P}^2(C)$.

\noindent$\bullet$ Ten trigonals $T^\varepsilon_i, \ 1\le i\le 5,
\ \ \varepsilon=0, 1$, each with a double cover
$\widetilde{T}^\varepsilon_i$: each of the 5 points $p_i\in Q\cap
L$ determines two $g^1_4$'s on $C$, cut out by the rulings
$R^\varepsilon_i$ on $A_{p_{i}}$, and the $(T^\varepsilon_i,
\widetilde{T}^\varepsilon_i)$ are associated to these by the
trigonal
construction.

\noindent$\bullet$  Sixteen Wirtinger covers $(X_j,
\widetilde{X}_j)\in\partial^I{\R}_6$:  the quartic del Pezzo
surface
$S_L$ contains 16 lines $\ell_j$ [Dem], each meeting $\Phi(C)$ in
two points, say $p_j, q_j$, and then $$X_j=C/(p_j\sim q_j)$$ and
$\widetilde{X}_j$ is its unique Wirtinger cover (1.9.I).

\bigskip

\noindent{\bf (4.3.8)} \ We observe that the generically finite map
\ \ ${\R}_6^{\mbox{Tet}}\to{\R}_6$ \ \ has \\ 1-dimensional fibers
over both ${\R\Q}$ and ${\R}^{\mbox{Trig}}_6$.  After blowing up
and normalizing, we obtain finite fibers generically over the
exceptional loci.  In the limit:

\noindent$\bullet$  Over $(Q,L)$, the 5 $g^1_4$'s correspond to
projections of the plane quintic $Q$ from one of the 5 points
$p_i\in Q\cap L$.

\noindent$\bullet$  Over $(T,D)$, with $T$ a trigonal curve, ${\cal
L}$ the trigonal bundle, and \\ $D\in|\omega_T\otimes{\cal
L}^{-2}|$, four of the $g^1_4$'s are of the form ${\cal L}(q)$ with
$q\in D$ (i.e. they are the trigonal ${\cal L}$ with base point
$q$); the fifth $g^1_4$ is $\omega_T\otimes{\cal L}^{-2}$.

\noindent$\bullet$  Given $X=C/(p\sim q)\in\partial\bar{\M}_5$,
there is a pencil $L\subset{\bf P}^2(C)$ of quadrics $A_p, \ p\in
L$, which contain both $\Phi(C)$ and its chord $\overline{pq}$.Among
these
there are 5 quadrics $A_{p_{i}}$ which are singular,
generically of rank 4.  Each of these has a single ruling $R_i$
containing a plane containing $\overline{pq}$.  These $R_i$ cut the
5 $g^1_4$'s on $X$.
We conclude that the tetragonal correspondence Tet of (4.1.1) lifts
to $$\widetilde{\rm Tet}\subset\widetilde{\R}_6 \times
\widetilde{\R}_6$$  which is
generically finite, of type (10,10), over each of our special loci.

\bigskip

\noindent {\bf Theorem 4.4 \ Structure of the blownup Prym map.}

Over each of the following loci, the blownup Prym map
$\widetilde{\p}$
has the listed monodromy group, and the lifted tetragonal
correspondence $\widetilde{\rm Tet}$ induces the listed structure.

\begin{list}{{\bf(\arabic{bean})}}{\usecounter{bean}}
\item $\widetilde{\C}$:  The group is $WE_6$, the structure is that
of
lines on a general non-singular cubic surface.
\item $\widetilde{\cal J}$:  The group is $WD_5$, the symmetry
group of
the incidence of lines on a quartic del Pezzo surface, or
stabilizer in $WE_6$ of a line.  The structure is that of lines on
a non-singular cubic surface, one of which is marked.
\item ${\B}$= the locus of intermediate Jacobians of Clemens'
quartic double solids of genus 5 \ [C1]:  The group is $WA_5=S_6$,
the structure is that of lines on a nodal cubic surface.

[Note:  ${\B}$ is contained in the branch locus of ${\cal P}$ [DS,
V.4] and in fact ([D6], and compare also [SV], [I]) equals the
branch locus.  The monodromy along ${\B}$ acting on a nearby,
unramified fiber is $({\Z}/2{\Z})\times S_6$, or the symmetry group
of a double-six, which is a subgroup of $WE_6$.  The group $S_6$
thus occurs as a subquotient of $WE_6$.]
\item (cf. [I]) $\widetilde{\cal P}$ extends naturally to the
boundary $\partial=\partial{\A}_5$; the monodromy is $WE_6$ and the
structure is that of lines on a general cubic surface.
\end{list}

We will prove parts (1), (2) and (3) in \S4.2. In the rest of this
section we
show
that theorems (4.1) and (4.2) follow from (4.4).

\bigskip

\noindent{\bf (4.5) Proofs of Theorem (4.2).}

By Theorem (2.16), Tet commutes with ${\cal P}$, therefore
$\widetilde{\rm Tet}$ commutes with $\widetilde{\cal P}$.  To
identify
this structure over a generic point, it suffices to do so over any
point over which $\widetilde{\cal P}$ and $\widetilde{\rm Tet}$ are
etale.  These conditions hold, e.g., over a generic
$(X,H)\in\widetilde{\C}$, where (4.4.1) identifies the structure.This
implies
that the monodromy is contained in $WE_6$, but we get
all of $WE_6$ already over $\widetilde{\C}$ (by (4.4.1) again), so
we are done.

We can work instead over $\widetilde{\cal J}$:  again,
$\widetilde{\cal P}$ and $\widetilde{\rm Tet}$ are etale over generic
$(C, L)\in\widetilde{\cal J}$, and $\widetilde{\rm Tet}$ has the
right
structure there by (4.4.2).  This shows $$WD_5\subset
\mbox{Monodromy} \subset WE_6.$$  As there are no intermediate
groups, the monodromy must equal $WD_5$ or $WE_6$.But if it were the
former,
$\widetilde{\R}_6$ would be reducible
(since $WD_5$ is the stabilizer in $WE_6$ of one of the 27 lines),
contradiction. \begin{flushright} QED \end{flushright}

\bigskip

\noindent{\bf Remark 4.5.1} \ Along the same lines, we can also
reprove Theorem (4.1).  Let $\widetilde{\rm Tet}^i$ denote the $i$-th
iterate of the correspondence $\widetilde{\rm Tet}$.  On ${\R\Q}^-$
we
have:
$$\widetilde{\rm Tet}^2 \mbox{ \ has degree} \ 27, $$
$$\widetilde{\rm Tet}^i=\widetilde{\rm Tet}^2 \;\;\;\mbox{for}\; i\ge
2.$$Since
$\widetilde{\rm Tet}$ is etale there, these properties persist
generically
on $\widetilde{\R}_6$.  Let $\sim$ be the equivalence relation
generated by $\widetilde{\rm Tet}$.  We conclude that $\sim$ has
degree
27, and that $\widetilde{\cal P}$ factors through a proper
quotient: $${\cal
P}':\widetilde{\R}_6/\sim\longrightarrow\widetilde{\A}_5.$$  We
still need to verify that $\deg ({\cal P}')=1$.  There are several
possibilities:

\noindent$\bullet$  We can work over $\widetilde{\cal J}$; as we will
see in (4.7), the fiber of $\widetilde{\cal P}$ there consists of
aunique
$\sim$-equivalence class;  so we need to check that ${\cal
P}'$ is unramified at that equivalence class.  This reduces to
seeing that $\widetilde{\cal P}$ is unramified at least at one
point of
the fiber; this is trivial at the plane-quintic point.  (This
argument avoids some of the detailed computations of the
codifferential on the boundary, [DS, Ch., IV], but is still very
close in spirit to [DS].)

\noindent$\bullet$  We could instead work over any other point of
$\widetilde{\A}_5$ over which we know the complete fiber, e.g. over
Andreotti-Mayer points, coming from bielliptic
Pryms, as in \S 3.  (This was proposed in [D1], as a way to avoid
the boundary computations.)

\noindent$\bullet$  Izadi [I] applies a similar argument over
boundary points, in $\partial {\A}_5$.  This lets her reduce the
degree computation over ${\A}_5$ to her results on ${\A}_4$, cf.
(4.9).

\subsection{Special Fibers}

In this section we exhibit the cubic surface of theorem (4.2)
explicitly over three special loci in ${\cal A}_{5}$. We do not know
how
to do this at the generic point of ${\cal A}_{5}$.

\noindent{\bf (4.6) Cubic threefolds}
{}From (4.3.4) we have an identification of $\widetilde{\cal
P}^{-1}(X,H)$, where $X\subset{\bf P}^4$ is a cubic threefold and
$H\subset{\bf P}^4$ a hyperplane, with the set of lines $\ell$ on
the cubic surface $X\cap H$.  For Theorem (4.4.1) we need to check
that two of these, say $\ell, \ell'\in F(X)$, intersect each other
if and only if the corresponding objects $(Q,\widetilde{Q},L),
(Q',\widetilde{Q}',L')$ correspond under $\widetilde{\rm Tet}$.  If
the
lines
$\ell,\ell'$ intersect, we are in the situation of (2.15.4), so the
corresponding objects $$(Q,\widetilde{Q},f),
(Q',\widetilde{Q}',f')$$
(notation of (2.15.4)) are tetragonally related.  Since $f,f'$ are
both cut out by hyperplanes through the span $A$ of $\ell,\ell'$,
we find points $$p\in Q\cap L, \ \ \ \ p'\in Q'\cap L'$$ (namely,the
projection
of $A$ from $\ell,\ell'$ respectively) such that
$f,f'$ are the projections of $Q$ from $p$ and of $Q'$ from $p'$,
respectively.  The description of $\widetilde{\R}_6^{\rm Tet}$ in
(4.3.8)
then shows that
$$((Q,\widetilde{Q},L),(Q',\widetilde{Q}',L'))\in\widetilde{\rm
Tet},$$
as required.  Since both the line incidence and $\widetilde{\rm Tet}$
are of bidegree (10,10), and we have an inclusion, it must be an
equality.
This shows that $\widetilde{\rm Tet}$ induces on $\widetilde{\cal
P}^{-1}(X,H)$ the structure of line incidence on the cubic surface
$X\cap H$.  Fix the ambiant ${\bf P}^4$ and the hyperplane $H$, and
let the cubic $X$ vary.  We clearly get all cubic surfaces in $H$
as intersections $X\cap H$; therefore the monodromy group is the
full symmetry group of the line configuration.   This completes the
proof of (4.4.1), hence also of Theorem (4.2).

\bigskip

\pagebreak

\noindent{\bf (4.7) Jacobians}
Start with $(C,L)\in\widetilde{\cal J}$.  The fiber
$\widetilde{\cal P}^{-
1}(C,L)$ consists of the 27 objects listed in (4.3.7).  Each of
these comes with the 5 $g^1_4$'s given in (4.3.8).  These give the
correspondence $\widetilde{\rm Tet}$, which we claim is equivalent to
the
line incidence on a cubic surface.

Let $S=S_L$ be the quartic del Pezzo surface determined by $(C,L)$,
as in (4.3.6).  Let $S'$ be its blowup at a generic point $r\in S$.
Then $S'$ is a cubic surface; its lines correspond to:

 \noindent$\bullet$  $\ell_Q$, the exceptional divisor over $r$.

 \noindent$\bullet$  10 conics through $r$ in $S$; these correspond
naturally to the 10 rulings ${\R}^\varepsilon_i$ (as in
(4.3.7)).[Each
${\R}^\varepsilon_i$ contains a unique plane through
$r$, which meets $S$ in a conic through $r$.]

\noindent$\bullet$  The 16 lines $\ell_j$ in $S$.

There is thus a natural bijection between the lines of $S'$ and
\linebreak
$\widetilde{\cal P}^{-1}(C,L)$.  We need to check that this
correspondence takes incident lines to covers which are
tetragonally related to each other through the $g^1_4$'s of
(4.3.8).  To that end, we list the effects of the tetragonal
constructions on our curves.  The details are straightforward, and
are omitted.

\bigskip

\noindent{\bf (4.7.1)} \ The quintic $(Q,\widetilde{Q})$, with the
$g_4^1: \ \ {\cal O}_Q(1)(-p_i), \ \ p_i\in Q\cap L$, \ \ is taken
to the two trigonals $$(T_i^\varepsilon,
\widetilde{T}^\varepsilon_i), \;\;\;\;\;\varepsilon=0,1,$$  each
with its unique base-point-free $g_4^1, \ \ \omega_T\otimes{\cal
L}^{-2}$.

\bigskip

\noindent{\bf (4.7.2)} \ The trigonal
$(T^\varepsilon_i,\widetilde{T}^\varepsilon_i)$with its
base-point-free $g^1_4$
goes to $(Q,\widetilde{Q})$ with
${\cal O}_Q(1)(-p_i)$, and to $(T^{1-\varepsilon}_i,
\widetilde{T}^{1-
\varepsilon}_i)$ with the base-point-free $g^1_4$.

Consider $(T^\varepsilon_i,\widetilde{T}^\varepsilon_i)$ with the
$g^1_4 \ \ {\cal L}^\varepsilon_i(p).$  The actual 4-sheeted cover
of ${\bf P}^1$ in this case is reducible, consisting of the
trigonal $T^\varepsilon_i$ together with a copy of ${\bf P}^1$
glued to it at $p$.  We are thus precisely in the situation of
Example (2.15.2):  both tetragonally related objects are Wirtinger
covers $(X_j,\widetilde{X}_j)$.

\bigskip

\noindent{\bf (4.7.3)} \ A Wirtinger cover $(X_j,\widetilde{X}_j)$
with
the $g^1_4$ cut out by the ruling ${\R}^\varepsilon_i$ on the
singular quadric $A_{p_{i}}$, is taken to the trigonal
$(T^\varepsilon_i,\widetilde{T}^\varepsilon_i)$ and to another
Wirtinger cover.

\bigskip

\pagebreak[4]

\noindent{\bf (4.8) Quartic double solids and the branch locus of
${\cal P}$}.

The fiber of ${\cal P}$ over the Jacobian $J(X)\in{\cal B}$ of a
quartic double solid $X$ of genus 5 is described in [DS, V.4],
following ideas of Clemens.  It consists of 6 objects
$(C_i,\widetilde{C}_i), \ \ 0\le i\le 5$, each with multiplicity 2,
and
15 objects $(C_{ij},\widetilde{C}_{ij}), \ \ 0\le i<j\le 5$.  The
monodromy group $S_6$ permutes the six values of $i$:  clearly the
two sets $\{C_i\}$ and $\{C_{ij}\}$ must be separately permuted,
and any permutation of the $C_i$ induces a unique permutation of
the $C_{ij}$.  The situation is precisely that of lines on a nodal
cubic surface:  the $C_i$ correspond to lines $\ell_i$ through the
node; and the plane through $\ell_i,\ell_j$ meets the cubic
residually in a line $\ell_{i,j}$.

The best way to see the symmetry is to consider Segre's cubic
threefold $Y\subset{\bf P}^4$, image of ${\bf P}^3$ by the linear
system of quadrics through 5 points $p_i$, $1\le i\le 5$, in
general position in ${\bf P}^3$. (cf. [SR] for the details.)  $Y$
contains six irreducible, two-dimensional families of lines, which
we call the ``rulings" $R_i, \ \ 0\le i\le 5$:  For $1\le i\le 5,
\ \ R_i$ consists of proper transforms of lines through $p_i$;
while $R_0$ parametrizes twisted cubics through $p_1,\cdots,p_5$.$Y$
also
contains 15 planes $\Pi_{ij}, \ \ 0\le i < j\le 5$
\linebreak (= the 5
exceptional divisors and the proper transforms of the 10 planes
$\overline{p_ip_jp_k}$); the ruling $R_i$ is characterized as the
set of lines in ${\bf P}^4$ meeting the 5 planes $\Pi_{ij}, \ \
j\ne i$.  The symmetric group $S_6$ acts linearly on ${\bf P}^4$,
preserving $Y$, permuting the $R_i$ and correspondingly the
$\Pi_{ij}$.

The quartic double solids in question are essentially the double
covers $$\zeta: X\to Y$$  branched along the intersection of $Y$
with a quadric $Q\subset{\bf P}^4$.  The Prym fiber is obtained as
follows:
\begin{tabbing}
$\bullet$  \=$C_i:= \{$ lines $\ell\in R_i$, tangent to $Q\}$ \\
\> $\widetilde{C}_i:= \{$ irreducible curves $\ell'\subset X$
such that $\zeta(\ell')=\ell\in C_i\}$
\end{tabbing}

Thus $(C_i,\widetilde{C}_i)$ is the discriminant of a conic-bundle
structure on $X$ given by $\zeta^{-1}(R_i)$.  The Prym canonical
curve $\Psi(C_i)\subset{\bf P}^4$ is traced by the tangency points
of $\ell$ and $Q$; in particular, $\Psi(C_i)\subset Q$, so
$(C_i,\widetilde{C}_i)$ is a ramification point of ${\cal P}$, by
(1.6).

\noindent$\bullet$  $(C_{ij},\widetilde{C}_{ij})$ is similarly
obtained
as discriminant of a conic bundle structure on $X$ given by
projection from $\Pi_{ij}$, cf. [DS, V4.5].

\bigskip

\noindent{\bf (4.9) Boundary behavior}

In [I], Izadi uses results on the structure of
${\p}:{\R}_5\to{\A}_4$ to find the incidence structure on the
fibers of the compactified map
$\doublebar{\p}:\doublebar{\R}_6\to\bar{\A}_5$
over boundary points of the toroidal compactification $\bar{\A}_5$.
The picture is as follows:

\[
\begin{array}{cc}
\begin{diagram}[AA] \node{\begin{array}{r} \; \doublebar{\cal P}  \;
\; :
\end{array}}
\end{diagram} &
\begin{diagram}[AA]
\node{\doublebar{\cal R}_{6}}
\arrow{e}
\node{\bar{\cal A}_{5} }
\end{diagram} \\
  &
\begin{diagram}[AA]
\node{\cup} \node{\cup}
\end{diagram}
\\
\begin{diagram}[AA]
\node{\begin{array}{r} \partial {\cal P}  \; \; : \end{array}} \\
\node{\begin{array}{r} \; {\cal P}  \; \; :  \end{array} }
\end{diagram} &
\begin{diagram}[AA]
 \node{\partial^{\rm II}\doublebar{\cal R}_{6}} \arrow{e}
\arrow{s,l}{\alpha}
\node{\partial\bar{\cal A}_{5} } \arrow{s,r}{\beta} \\
 \node{{\cal R}_{5}}  \arrow{e}
\node{{\cal A}_{4}}
 \end{diagram}
\end{array}
\]

Over general $A\in{\A}_4$, the fiber $\beta^{-1}(A)$ is isomorphic
to the Kummer variety $A/(\pm 1)$.  Over
$(\widetilde{C},C)\in{\R}_5$, the fiber of $\alpha$ is
$S^2\widetilde{C}/\iota$, and $\partial{\cal P}$ becomes (cf. [D3,
(4.6)]) the map

\[
\begin{array}{l}
x+y \mapsto \psi (x) + \psi (y)
 \\
\begin{diagram}[AA]
\node{S^{2}\widetilde{C}} \arrow{e} \arrow{s} \node{A}  \arrow{s} \\
\node{S^{2}\widetilde{C}/\iota} \arrow{e} \node{A/(\pm 1)}
\end{diagram}
\end{array}
\]

\noindent where $\psi$ is the Abel-Prym map $\widetilde{C}\to A$.All
in all
then, we are considering the map
$$\partial{\p}:\cup_{(\widetilde{C},C)\in{\p}^{-1}(A)}S^2
\widetilde{C}=:E\longrightarrow A.$$  Theorem (4.1) says that its
degree is 27, and Theorem (4.2) predicts an incidence structure on
its fibers, i.e. a way of associating a cubic surface to each point
$a \in A$.

In \S5 we associate to $A\in{\A}_4$ a cubic threefold
$X=\kappa(A)\subset{\bf P}^4$ such that ${\cal P}^{-1}(A)$ is a
double cover of the Fano surface $F(X)$ of lines in $X$.  For
generic $a\in A$, we are looking for a cubic surface; it is
reasonable to hope that this should be of the form $H(a)\cap X$,
where $H(a)$ is an appropriate hyperplane in ${\bf P}^4$.  We thus
want a map $$H:A\to({\bf P}^4)^*$$  such that  $$pr((\partial{\cal
P})^{-1}(a))=\{\mbox{lines in} \ H(a)\cap X\}.$$

\[
\begin{diagram}[A]
\node{E} \arrow[2]{e,t}{\partial {\cal P}} \arrow[2]{s,l}{pr}
\arrow{se}
\node[2]{A}
\arrow[2]{s,r}{H} \\
\node[2]{{\cal P}^{-1}(A)} \arrow{sw} \\
\node{F(X)} \node[2]{({\bf P}^{4})^{*}}
\end{diagram}
\]
Izadi's beautiful observation is that such an $H$ is given by the
linear system $\Gamma_{00}$ (sections of $|2\Theta|$ vanishing to
order $\ge 4$ at 0).  The identification of $\Gamma_{00}$ with the
ambiant ${\bf P}^4$ of $X$ uses a construction of Clemens relating
his double solids to $\Gamma_{00}$, and the interpretaton of (a
cover of) $X$ as parametrizing double
solids with intermediate Jacobians isomorphic to $A$, cf. [D6] or
[I].

\section{Fibers of $P:{\cal R}_5\to A_4.$}

\subsection{The general fiber.}

\ \ \ \ Our main result in this section is:

\bigskip

\noindent{\bf Theorem 5.1}  For generic $A\in A_4$, the fiber
$\overline{\cal P}^{-1}(A)$ is isomorphic to a double cover of
the Fano
surface $F=F(X)$ of lines on some cubic threefold $X$.

Let ${\cal R}{\cal C}$ denote the inverse image in ${\cal R}A_5$
of
the locus ${\cal C}$ of (intermediate Jacobians $J(X)$ of)
cubic threefolds $X$.  We recall from [D4] that it splits into
even
and odd components:
\noindent{\bf (5.1.1)} \ \ \ \ \ \ \ \ \ \ \ \ \ \ \ \ \ \ \ \ \ \
${\cal R}{\cal C}={\cal R}{\cal C}^+\amalg {\cal R}{\cal C}^-,$
\noindent distinguished by a parity funciton.  This follows from
the existence of a natural theta divisor $\Xi\subset J(X)$,
characterized (cf. [CG]) by having a triple point at $0:\Xi$
translates the parity function $q$ of (1.2), on theta
characteristics, to a parity on semiperiods.  More explicitly,
pick
$(Q,\sigma)\in{\p}^{-1}(J(X))\subset{\cal R}{\cal Q}^-$;
Mumford's
exact sequence (Theorem (1.4)(2)) says that any $\delta\in
J_2(X)$ is $\pi^*\nu$ for some $\nu\in(\sigma)^\perp\subset
J_2(Q)$.  The compatibility result, theorem (1.5), then gives
(cf.
[D4], Proposition (5.1)):

\noindent{\bf (5.1.2)} \ \ \ \ \ \ \ \ \ \ \ \ \ \ \ \ \ \ \ \ \ \
$q_X(\delta)=q_Q(\nu)=q_Q(\nu \sigma).$

In case $\delta$ is even, we end up with an isotropic subgroup \\
$(\!\nu,\sigma\!) \ \ \subset \ \ J_2(Q)$, with $\sigma$ odd and
$\nu, \nu\sigma$ even.  The Pryms of the latter are therefore
Jacobians of curves:
\noindent{\bf (5.1.3)} \ \ \ \ \ \ \ \ \ \ \ \ \ \ \ \ \
$P(Q,\nu)\approx J(C), \ \ \ P(Q,\nu\sigma)\approx J(C'),$

\noindent and the image of $\sigma$ gives semiperiods $\mu\in
J_2(C),\mu'\in J_2(C')$.
Reversing direction, we can construct an involution
$$\lambda:{\cal
R}_5\longrightarrow {\cal R}_5$$ and a map $$\kappa: {\cal
R}_5\longrightarrow{\cal R}{\cal C}^+,$$ as follows:  Start with
$(C,\mu)\in {\cal R}_5$, pick the unique \ \ $(Q,\nu)$ \ \ in \\
${\cal P}^{-1}(C)\cap {\cal R}{\cal Q}^+$, and let
$\sigma,\nu\sigma\in J_2(Q)$ map to $\mu\in J_2(C).$  Then
formula (1.3) reads:
\noindent{\bf (5.1.4)} \ \ \ \ \ \ \ \ \ \ \ \ \ $0\equiv3+{\rm \
even}+q(\sigma)+q(\nu\sigma) \ \ \ \ \ {\rm (mod. \ 2)},$

\noindent so after possibly relabeling, we may assume
$$(Q,\sigma)\in {\cal RQ}^-, \ \ \ (Q,\nu\sigma)\in{\cal RQ}^+$$
so
that there is a well-defined curve $C'\in{\cal M}_5$ and a cubic
threefold $X\in{\cal C}$ such that
$$P(Q,\sigma)\approx J(X)$$ {\bf (5.1.5)}
$$P(Q,\nu\sigma)\approx J(C').$$
\noindent We can thus define $\lambda$ and $\kappa$ by:
$$\lambda(C,\mu):=(C',\mu')$$ {\bf (5.1.6)}
$$\kappa(C,\mu):=(X,\delta),$$

\noindent where $\mu'\in J_2(C'), \ \ \ \delta\in J_2(X)$ are the
images of $\nu\in J_2(Q).$  The precise version of our results is
in terms of $\lambda$ and $\kappa$:

\bigskip

\noindent{\bf Theorem 5.2}
\begin{list}{{\rm(\arabic{bean})}}{\usecounter{bean}}
\item $(C,\mu)$ is related to $\lambda(C,\mu)$ by a sequence of
two
tetragonal constructions.  Hence $\lambda$ commutes with the Prym
map:$${\cal P}\circ\lambda={\cal P}, \ \ \
\lambda\circ\lambda=id.$$\item $\kappa$ factors through the Prym map:
$$\kappa:{\cal R}_5\stackrel{{\cal
P}}{\longrightarrow}A_4\stackrel{\chi}{\longrightarrow}{\cal
RC}^+,$$ where $\chi$ is a birational map.
\end{list}

Recall the Abel-Jacobi map [CG], $$AJ:F(X)\longrightarrow
J(X),$$ which is well-defined up to translation in $J(X)$.  (It
can
be identified with the Albanese map of the Fano Surface $F(X)$.)A
point
$\delta\in J_2(X)$ determines a double cover of $J(X)$,
hence of $F(X)$.
\bigskip

\noindent{\bf Theorem 5.3}  For generic $A\in{\cal A}_4$, set
$$(X,\delta):=\chi(A)=\kappa({\cal P}^{-1}(A))\in {\cal R C}^+.$$

Let $F(X)$ be the Fano surface of $X$, $\widetilde{F(X)}$ its double
cover determined by $\delta$ via the Abel-Jacobi map.
\begin{list}{{\rm(\arabic{bean})}}{\usecounter{bean}}
\item  There is a natural isomorphism
$$P^{-1}(A)\approx\widetilde{F(X)}.$$
\item The action of $\lambda$ on the left corresponds to the
sheet
interchange on the right.
\item Two objects $(C,\mu), (C',\mu)\in{\p}^{-1}(A)$ are
tetragonally related if and only of the lines $\ell, \ell'\in
F(X)$
which they determine intersect.
\end{list}

\bigskip

\noindent{\bf Remark 5.4}  Izadi has recently analyzed the
birational map $\chi$, in [I].  In particular, she shows that
$\chi$ is an isomorphism on an explicitly described, large open
subset of ${\cal A}_4$.

\subsection{Isotropic subgroups.}

\begin{tabbing} X \= \kill
\>By isotropic subgroup of rank $r$ on a curve $C$ we mean an
\\$r$-dimensional
${\bf F}_2$-subspace of $J_2(C)$ on
which the intersection pairing \\ $\langle \ , \ \rangle$ is
identically
zero.Choosing an isotropic subgroup of rank 1 is \\ the same as
choosing
a non-zero semiperiod.
\end{tabbing}

Start with a trigonal curve $T\in {\cal M}_{g+1}$, with a rank-2
isotropic subgroup $W\subset J_2(T)$ whose non-zero elements
we denote $\nu_i, \ \ i=0,1,2.$  The trigonal construction
associates to $(T,\nu_i)$ the tetragonal curve $X_i \in {\cal
M}_g$.Mumford's sequence (1.4)(2) sends $W$ to an isotropic subgroup
of
rank $1$ on $X_i$, whose non-zero element we denote $\mu_i$.
\bigskip

\noindent{\bf Lemma 5.5}  The construction above sets up a
bijection between the following data:

\noindent $\bullet$  A trigonal curve $T\in {\cal M}_{g+1}$ with
rank-2 isotropic subgroup.

\noindent $\bullet$  A tetragonally related triple
$(X_i,\mu_i)\in {\cal R}_g, \ \ \ i=0, 1, 2$.

\bigskip

\noindent{\bf Proof.}

We think of $WD_4$ as the group of signed permutations of the 8
objects $\{x^{\pm}_i\}$,  $1\leq i\leq 4$. \ Start with a
tetragonal double cover \linebreak  $\widetilde{X}_0\longrightarrow
X_0\longrightarrow{\bf P}^1.$  It determines a principal
$WD_4$-bundle over ${\bf P}^1\backslash$(Branch).  The original
covers $\widetilde{X}_0, X_0$ are recovered as quotients by the
following subgroups of $WD_4$: $$\widetilde{H}_0:={\rm
Stab}(x^+_1),$$
$$H_0:={\rm Stab}(x^\pm_1),$$
Consider also the subgroup $$G:={\rm Stab}\{\{x^+_1, x^+_2\},
\{x_1^-,x_2^-\}\}.$$  It has index 12 in $WD_4$.  Its normalizer
is:$$N(G)={\rm Stab}\{\{x^\pm_1, x^\pm_2\}, \{x^\pm_3,
x^\pm_4\}\},$$ of index 3.  The quotient is $$N(G)/G\approx({\bf
Z}/2{\bf Z})^2,$$  so there are 3 intermediate groups
$\widetilde{G}_i,\ \ \ i=0,1,2.$  We single out one of
them:$$\widetilde{G}_0:={\rm Stab}\{x^\pm_1, x^\pm_2\}.$$  The three
subgroups $\widetilde{G}_i$ are not conjugate to each other, but can
be
taken to each other by outer automorphisms of $WD_4$.  In fact,
the
action of Out$(WD_4)\approx S_3$ sends $G$, and hence also
$N(G)$,
to conjugate subgroups; it permutes the $\widetilde{G}_i$
transitively,
modulo conjugation; and it also takes $H_0, \widetilde{H}_0$ to
non-conjugate subgroups $H_i, \widetilde{H}_i, \ \ \ i=1,2.$  We
illustrate each of these subgroups as the stabilizer in $WD_4$ of
a corresponding partition of $\left(\begin{array}{llll}x^+_1 &
x^+_2 & x^+_3 & x^+_4 \\ x^-_1 & x^-_2 & x^-_3 & x^-_4
\end{array}
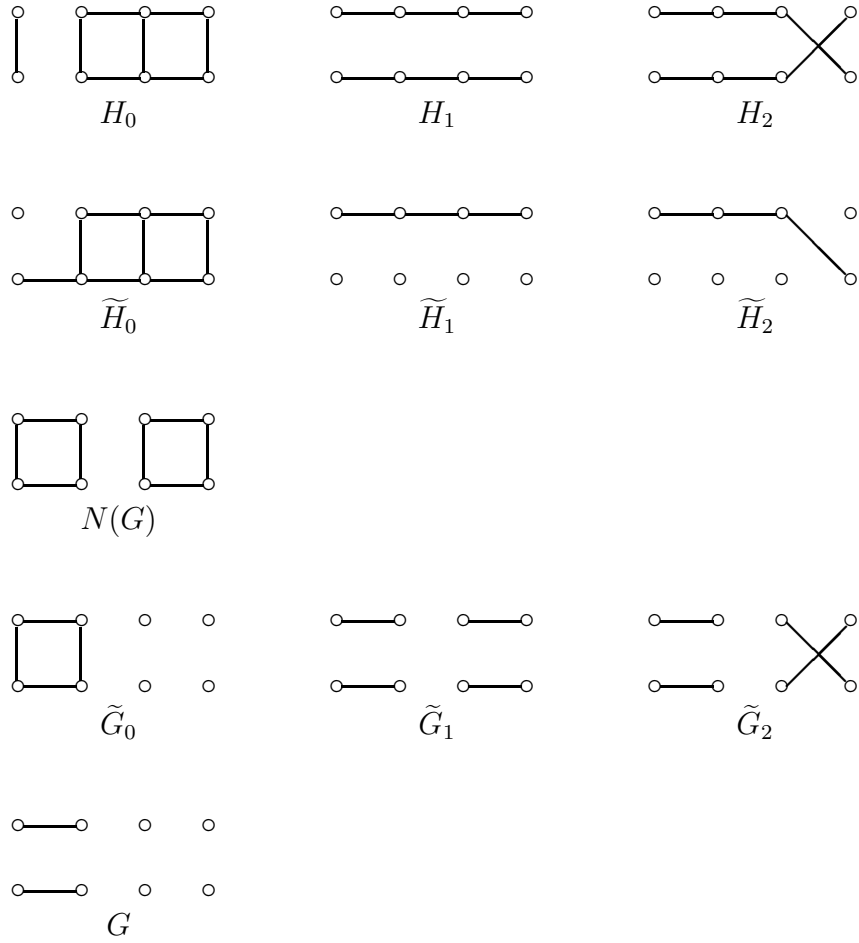
\right)$:

\begin{center}
\begin{tabular}{ccc}
\hspace{1.5in} & \hspace{1.5in} & \hspace{1.5in} \\
\begin{picture}(82,30)(2,1)
\thicklines
\put(2,1){$\circ$} \put(2,25.8){$\circ$}
\put(4.5,6){\line(0,1){20}}
\put(26,1){$\circ$} \put(31,3.5){\line(1,0){20}}
\put(28.5,6){\line(0,1){20}}
\put(26,25.8){${\circ}$} \put(31,28.3){\line(1,0){20}}
\put(50,1){${\circ}$}  \put(55,3.5){\line(1,0){20}}
\put(52.5,6){\line(0,1){20}} \put(50,25.8){${\circ}$}
\put(55,28.3){\line(1,0){20}}
\put(74,1){${\circ}$} \put(76.5,6){\line(0,1){20}}
\put(74,25.8){${\circ}$}
\end{picture} &
\begin{picture}(82,30)(2,1)
\thicklines
\put(2,1){$\circ$} \put(2,25.8){$\circ$} \put(7,28.3){\line(1,0){20}}
\put(7,3.5){\line(1,0){20}}
\put(26,1){$\circ$} \put(31,3.5){\line(1,0){20}}
\put(26,25.8){${\circ}$} \put(31,28.3){\line(1,0){20}}
\put(50,1){${\circ}$}  \put(55,3.5){\line(1,0){20}}
 \put(50,25.8){${\circ}$}
\put(55,28.3){\line(1,0){20}}
\put(74,1){${\circ}$}
\put(74,25.8){${\circ}$}
\end{picture} &
\begin{picture}(82,30)(2,1)
\thicklines
\put(2,1){$\circ$} \put(2,25.8){$\circ$} \put(7,28.3){\line(1,0){20}}
\put(7,3.5){\line(1,0){20}}
\put(26,1){$\circ$} \put(31,3.5){\line(1,0){20}}
\put(26,25.8){${\circ}$} \put(31,28.3){\line(1,0){20}}
\put(50,1){${\circ}$}  \put(55,3.9){\line(1,1){22.4}}
 \put(50,25.8){${\circ}$}
\put(55,27.9){\line(1,-1){22}}
\put(76,1){${\circ}$}
\put(76,25.8){${\circ}$}
\end{picture}  \\
$H_{0}$ & $H_{1}$ & $H_{2}$  \\
  &  &  \\
    &  &  \\
\begin{picture}(82,30)(2,1)
\thicklines
\put(2,1){$\circ$} \put(2,25.8){$\circ$}
\put(7,3.5){\line(1,0){20}}
\put(26,1){$\circ$} \put(31,3.5){\line(1,0){20}}
\put(28.5,6){\line(0,1){20}}
\put(26,25.8){${\circ}$} \put(31,28.3){\line(1,0){20}}
\put(50,1){${\circ}$}  \put(55,3.5){\line(1,0){20}}
\put(52.5,6){\line(0,1){20}} \put(50,25.8){${\circ}$}
\put(55,28.3){\line(1,0){20}}
\put(74,1){${\circ}$} \put(76.5,6){\line(0,1){20}}
\put(74,25.8){${\circ}$}
\end{picture} &
\begin{picture}(82,30)(2,1)
\thicklines
\put(2,1){$\circ$} \put(2,25.8){$\circ$} \put(7,28.3){\line(1,0){20}}
 \put(26,1){$\circ$}
\put(26,25.8){${\circ}$} \put(31,28.3){\line(1,0){20}}
\put(50,1){${\circ}$}
 \put(50,25.8){${\circ}$}
\put(55,28.3){\line(1,0){20}}
\put(74,1){${\circ}$}
\put(74,25.8){${\circ}$}
\end{picture} &
 \begin{picture}(82,30)(2,1)
\thicklines
\put(2,1){$\circ$} \put(2,25.8){$\circ$} \put(7,28.3){\line(1,0){20}}
\put(26,1){$\circ$}
\put(26,25.8){${\circ}$} \put(31,28.3){\line(1,0){20}}
\put(50,1){${\circ}$}
 \put(50,25.8){${\circ}$}
\put(55,27.9){\line(1,-1){22}}
\put(76,1){${\circ}$}
\put(76,25.8){${\circ}$}
\end{picture}  \\
$\widetilde{H}_{0}$ & $\widetilde{H}_{1}$ & $\widetilde{H}_{2}$ \\
  &  &  \\
    &  &  \\
 \begin{picture}(82,30)(2,1)
\thicklines
\put(2,1){$\circ$} \put(2,25.8){$\circ$} \put(7,28.3){\line(1,0){20}}
\put(7,3.5){\line(1,0){20}}   \put(4.5,6){\line(0,1){20}}
\put(26,1){$\circ$}  \put(28.5,6){\line(0,1){20}}
\put(26,25.8){${\circ}$}
\put(50,1){${\circ}$}  \put(55,3.5){\line(1,0){20}}
\put(52.5,6){\line(0,1){20}} \put(50,25.8){${\circ}$}
\put(55,28.3){\line(1,0){20}}
\put(74,1){${\circ}$} \put(76.5,6){\line(0,1){20}}
\put(74,25.8){${\circ}$}
\end{picture} &  & \\
$N(G)$ & & \\
 & & \\
   &  &  \\
  \begin{picture}(82,30)(2,1)
\thicklines
\put(2,1){$\circ$} \put(2,25.8){$\circ$} \put(7,28.3){\line(1,0){20}}
\put(7,3.5){\line(1,0){20}}   \put(4.5,6){\line(0,1){20}}
\put(26,1){$\circ$}  \put(28.5,6){\line(0,1){20}}
\put(26,25.8){${\circ}$}
\put(50,1){${\circ}$}
 \put(50,25.8){${\circ}$}
\put(74,1){${\circ}$}
\put(74,25.8){${\circ}$}
\end{picture} &
  \begin{picture}(82,30)(2,1)
\thicklines
\put(2,1){$\circ$} \put(2,25.8){$\circ$} \put(7,28.3){\line(1,0){20}}
\put(7,3.5){\line(1,0){20}}
\put(26,1){$\circ$}
\put(26,25.8){${\circ}$}
\put(50,1){${\circ}$}  \put(55,3.5){\line(1,0){20}}
 \put(50,25.8){${\circ}$}
\put(55,28.3){\line(1,0){20}}
\put(74,1){${\circ}$}
\put(74,25.8){${\circ}$}
\end{picture} &
 \begin{picture}(82,30)(2,1)
\thicklines
\put(2,1){$\circ$} \put(2,25.8){$\circ$} \put(7,28.3){\line(1,0){20}}
\put(7,3.5){\line(1,0){20}}
\put(26,1){$\circ$}
\put(26,25.8){${\circ}$}
\put(50,1){${\circ}$}  \put(55,3.9){\line(1,1){22.4}}
 \put(50,25.8){${\circ}$}
\put(55,27.9){\line(1,-1){22}}
\put(76,1){${\circ}$}
\put(76,25.8){${\circ}$}
\end{picture}  \\
$\widetilde{G}_{0}$ & $\widetilde{G}_{1}$ & $\widetilde{G}_{2}$ \\
  &  &  \\
    &  &  \\
  \begin{picture}(82,30)(2,1)
\thicklines
\put(2,1){$\circ$} \put(2,25.8){$\circ$} \put(7,28.3){\line(1,0){20}}
\put(7,3.5){\line(1,0){20}}
\put(26,1){$\circ$}
\put(26,25.8){${\circ}$}
\put(50,1){${\circ}$}
 \put(50,25.8){${\circ}$}
\put(74,1){${\circ}$}
\put(74,25.8){${\circ}$}
\end{picture} &   & \\
 $G$&  & \\
   &  & \\
     &  &
\end{tabular}
\end{center}

Let $X_i, \widetilde{X}_i, T, \stackrel{\approx}{T},\widetilde{T}_i \
\
\ (i=0,1,2)$ be the quotients of the principal $WD_4$-bundle by
the
subgroups $H_i, \widetilde{H}_i, N(G), G, \widetilde{G}_i$
respectively,
compactified to branched covers of ${\bf P}^1$.  We see
immediately
that:

\noindent $\bullet$  The trigonal construction takes $X_0\to {\bf
P}^1$ to $\widetilde{T}_0\to T\to{\bf P}^1.$

\noindent $\bullet$  The double cover $\widetilde{X}_0\to X_0$
corresponds via (1.4)(2) to the double cover
$\stackrel{\approx}{T}\to\widetilde{T}_0$.

\noindent $\bullet$  The tetragonal construction acts by outer
antomorphisms,
hence exchanges the three tetragonal double covers
$\widetilde{X}_i\to X_i\to{\bf P}^1.$

Applying the same outer automorphisms, we see that the trigonal
construction also takes $X_i\to{\bf P}^1$ to $\widetilde{T}_i\to T\to
{\bf P}^1, \  i=1,2.$  To a tetragonally related triple
$(\widetilde{X}_i\to X_i \to {\bf P}^1)$ we can thus unambiguously
associate the trigonal \ \ $T\to{\bf P}^1$ \ \ together with the
rank-2, isotropic subgroup corresponding to the covers
$\widetilde{T}_i$.  This inverts the construction predecing the
lemma.
\begin{flushright} QED. \end{flushright}

\bigskip

\noindent{\bf Note 5.5.1}  The basic fact in the above proof is
that the 3 tetragonals $X_i$ yield the same trigonal $T$.  This
can be explained more succinctly:  outer automorphisms
take the natural surjection $\alpha_0:WD_4\to\!\to S_4$ to
homomorphisms $\alpha_1,\alpha_2$ which are not conjugate to it.But
the
composition $\beta\circ\alpha_i:WD_4\to\!\to S_3$, where
$\beta:S_4\to\!\to S_3$ is the Klein map, are conjugate to each
other.

\bigskip

\noindent{\bf Construction 5.6}  Now let $T\in{\cal M}_{g+1}$ be
a trigonal curve, together with an isotropic subgroup of rank 3,
$$V\subset J_2(T).$$  We think of $V$ as a vector space over
${\bf F}_2$; the projective plane ${\bf P}(V)$ is identified with
$V\backslash(0)$.  For each $i\in {\bf P}(V)$, the trigonal
construction gives a tetragonal curve $Y_i\in{\cal M}_g$.Mumford's
sequence
(1.4)(2) gives an isotropic subgroup of rank
2,
$$W_i\subset J_2(Y_i),$$ with a natural identification
$W_i\approx
V/(i).$

Let $U\subset V$ be a rank-2 subgroup, so ${\bf P}(U)\subset{\bf
P}(V)$ is a projective line.  Lemma (5.5) shows that the 3
objects $$(Y_i, U/(i))\in{\bf R}_g, \ \ \ \ \ \ \ \ i\in {\bf
P}(U)$$ are
tetragonally related.  In particular, they have a
common Prym variety $$P_U\approx{\cal P}(Y_i, U/(i))\in {\cal
A}_{g-1}, \ \ \ \ \ \forall i\in {\bf P}(U).$$Applying (1.4) twice,
we see that
the original rank-3 subgroup
$V$
determines a rank-1 subgroup $$V/U\subset(P_{U})_2,$$
so we let $\mu_{U}\in(P_U)_2$ be its non-zero element.Altogether
then, we have a map $${\bf P}(V)^*\longrightarrow{\cal R
A}_{g-1}$$
$$U\longmapsto(P_{U},\mu_{U}).$$

\bigskip

\noindent {\bf (5.6.1)}  Assume now that one of the $Y_i$ happens
to be trigonal.  (This can only happen if $g\le 6.$)  Whenever
$U\ni i$, we find a tetragonal curve $C_U\in {\cal
M}_{g-1}$ such that $P_U\approx J(C_U)$.  Lemma
(5.5), applied to $(Y_i, W_i)$, shows that the 3 objects
$$(C_U,\mu_U)\in{\cal R}_{g-1},  \ \ \ \ \ U\ni i$$ are
tetragonally related, so they have a common Prym variety \\
$A=P_V\in{\cal A}_{g-2}.$

\bigskip

\noindent {\bf (5.6.2)}  Assume instead that $g=6$ and that
$P_U$ happens to be a Jacobian $J(C_U)\in{\cal J}_5$, for some
$U\in{\bf P}(V)^*.$  Of the three $Y_i, \ \ i\in U$, we claim two
are trigonal and the third, a plane quintic.  Indeed, by (4.7),
the tetragonal triples above $J(C_U)$ consist either of a plane
quintic and two trigonals, as claimed, or of a trigonal and two
Wirtingers.  The latter is excluded since the isomorphism
$$J(Y_i)\approx P(T,i)$$ implies that $Y_i$ is non-singular for
each $i\in{\bf P}(V)$.

Assume from now on that $g=6$.  Our data consists of:

\noindent $\bullet$  $T\in {\cal M}_7$, trigonal, with $V\subset
J_2(T)$
isotropic of rank 3.

\noindent $\bullet$  For each $i\in {\bf P}(V)$, a curve $Y_i \in
{\cal
M}_6$, with a rank-2 isotropic subgroup $W_i\subset J_2(Y_i)$.

\noindent $\bullet$  For each $U\in {\bf P}(V)^*$, an object
$(P_U,\mu_U)\in
{\cal RA}_5$

\noindent $\bullet$  An abelian variety $A=P_V\in {\cal A}_4.$

We display ${\bf P}(V)$ as a graph with 7 vertices $i\in {\bf
P}(V)$ and 7 edges \\ $U\in{\bf P}(V)^*$, in
(3,3)-correspondence.We write $T$ (or $Q$) on a vertex corresponding
to a
trigonal (or
quintic) curve, and $C$ on an edge corresponding to a Jacobian.We
restate our
observations:

\bigskip

\noindent{\bf (5.6.3)}:  Edges through a $T$-vertex are $C$-edges.
\bigskip

\noindent{\bf (5.6.4)}:  On a $C$-edge, the vertices are $T,T,Q$.

\bigskip

It follows that only one configuration is possible:

\pagebreak[4]

{ \ }
\vspace{4in}

\centerline{ {\bf Figure 5.7}}

\noindent Thus four of the $Y_i$ are trigonal, the other three
are
quintics, and six of the $P_U$, corresponding to the straight
lines, are Jacobians of curves.  Let $U_0\in{\bf P}(V)^*$
correspond to the circle.  For $i\in U_0$, \ \ $Y_i$ is a quintic
$Q$.  Through $Q$ pass two $C$ edges and $U_0$, and the
semiperiods
corresponding to the $C$-edges are even; by (1.3), the semiperiod
$U_0/(i)$ corresponding to $U_0$ must be \underline{odd}, so
there
is a cubic threefold $X\in{\cal C}$ such that $${\bf
P}_{U_0}\approx J(X).$$ Finally, theorem (1.5), or formula (5.1.2),
shows that the
semiperiod $\delta:=\mu_{U_0}\in J_2(X)$ is \underline{even}.

We observe that the three tetragonally related quintics
correspond
to 3 lines on the cubic threefold which meet each other and thus
form the intersection of $X$ with a (tritangent) plane.  We are
thus exactly in the situation of (2.15.4).

\subsection{Proofs.}

\noindent{\bf (5.8)}  Theorems (5.1),(5.2) and (5.3) all follow
from the following statements:
\begin{list}{{\rm(\arabic{bean})}}{\usecounter{bean}}
\item $(C,\mu)$ is related to $\lambda(C,\mu)$ by a sequence of
two tetragonal constructions.
\item $\kappa$ is invariant under the tetragonal construction
\item For $(X,\delta)\in {\cal R}{\cal C}^+, \ \ \
\kappa^{-1}(X,\delta)\approx\widetilde{F(X)}$, the isomorphism takes
$\lambda$ to the involution on $\widetilde{F(X)}$ over $F(X)$,
and two objects on the left are tetragonally related iff the
corresponding lines intersect.
\item Any two objects in ${\cal P}^{-1}(A)$, generic $A\in{\cal
A}_4$, are connected by a sequence of (two) tetragonal
constructions.
\end{list}

Indeed, (1) is (5.2)(1); \ (2) and (4) imply the existence of \\
$\chi:{\cal A}_4\longrightarrow{\cal R}{\cal C}^+$ such that
$\kappa=\chi\circ{\cal P}$, while (3) shows that any two objects
in a $\kappa$-fiber are also connected by a sequence of two
tetragonal constructions, so $\chi$ must be birational, giving
(5.2)(2).  This gives an isomorphism ${\cal
P}^{-1}(A)\approx\kappa^{-1}(X,\delta)$, so (5.3) follows.
\bigskip

\noindent{\bf (5.9)}  We let ${\cal R}^2{\cal Q}^+, {\cal R}^2{\cal
Q}^-$ denote the moduli spaces of plane quintic curves $Q$
together with:

 \noindent $\bullet$ A rank-2, isotropic subgroup $W\subset
J_2(Q)$, containing one odd and two even semiperiods, and

\noindent $\bullet$  a marked even (respectively odd) semiperiod in
$W\backslash(0).$

Exchanging the two even semiperiods gives an involution on ${\cal
R}^2{\cal Q}^+$, with quotient ${\cal R}^2{\cal Q}^-$.  The
birational map $$\alpha:{\cal M}_5\widetilde{\longrightarrow}{\cal
R}{\cal Q}^+,$$ of (4.3.7), lifts to a birational map
\noindent{\bf (5.9.1)} \ \ \ \ \ \ \ \ \ \ \ \ \ \ \ \ \ \ \ \ \ \
\
${\cal R}\alpha:{\cal R}_5\widetilde{\longrightarrow}{\cal R}^2{\cal
Q}^+.$

\noindent From the construction of $\lambda$ in (5.1.6) it
follows that the involution on the right hand side corresponds to
$\lambda$ on the left, so we have a commutative diagram:

\begin{equation}
\renewcommand{\theequation}{\bf
{\arabic{section}}.{\theau}.{\arabic{equation}}}
\setcounter{au}{9}
\setcounter{equation}{2}
\begin{diagram}[AA]
\node{{\cal R}_{5}} \arrow{s} \arrow{e,tb}{{\cal R}\alpha}{\sim}
\node{{\cal R}^{2}Q^{+}} \arrow{s,lr}{\pi}{2:1} \\
\node{{\cal R}_{5}/\lambda} \arrow{e,b}{\sim} \node{{\cal
R}^{2}Q^{-}}
\end{diagram}
\end{equation}

Start with $(C,\mu)\in{\cal R}_5$ and any $g^1_4$ on $C$.  The
trigonal construction produces a trigonal $Y\in{\cal M}_6$ with
rank-2, isotropic subgroup $W_Y$.  On $Y$ we have a natural
$g^1_4$, namely $w_Y\otimes L^{-2}$, where $L$ is the trigonal
bundle; so we bootstrap again, to a trigonal $T\in{\cal M}_7$
with
rank-3 isotropic subgroup $V$.  Applying construction (4.6) we
obtain a diagram like (5.7), including an edge for $(C,\mu)$ and
on
it a vertex for $(Q,W_Q):=$ \\ $\pi{\cal R}\alpha(C,\mu)$.  But then
$\lambda(C,\mu)$ and $\kappa(C,\mu)$ also appear in the same
diagram, as the two other edges (the line, respectively the
circle) through $Q$! Statement (5.8.1) now follows, since any two
edges of (5.7) which meet in a trigonal vertex are tetragonally
related.  (5.8.2) also follows, since any $(C',\mu')$
tetragonally related to
$(C,\mu)$ will appear in the same diagram with $(C,\mu)$ (for the
obvious initial choice of $g^1_4$ on $C$), so they have the same
$\kappa$.

{}From the restriction to ${\cal R}{\cal Q}^-$ of the Prym map we
obtain, by base change:

\begin{equation}
\renewcommand{\theequation}{\bf
{\arabic{section}}.{\theau}.{\arabic{equation}}}
\setcounter{au}{9}
\setcounter{equation}{3}
\begin{diagram}[AA]
\node{{\cal R}^{2}Q^{-}} \arrow{s,l}{\cal RP} \arrow{e}
\node{{\cal R}Q^{-}} \arrow{s,r}{\cal P} \\
\node{{\cal RC}^{+}} \arrow{e} \node{{\cal C}}
\end{diagram}
\end{equation}

Combining with (5.8)(1),(2) and (5.9.2), we find that $\kappa$
factors

\begin{equation}
\renewcommand{\theequation}{\bf
{\arabic{section}}.{\theau}.{\arabic{equation}}}
\setcounter{au}{9}
\setcounter{equation}{4}
\begin{diagram}[AA]
\node{{\cal R}_{5}} \arrow{s} \arrow{e,tb}{{\cal R}\alpha}{\sim}
\node{{\cal R}^{2}Q^{+}} \arrow{s,r}{\pi} \\
\node{{\cal R}_{5}/\lambda} \arrow{e,b}{\sim} \node{{\cal
R}^{2}Q^{-}}
\arrow{s,r}{\cal RP} \\
\node[2]{{\cal RC}^{+}}
\end{diagram}
\end{equation}

We know ${\cal P}^{-1}(X)$ from (4.6), so by (5.9.3):

\noindent{\bf (5.9.5)} \ \ \ \ \ \ \ \ \ \ \ \ \ \ \ \ \ \ \ ${\cal
RP}^{-1}(X,\delta)\approx{\cal P}^{-1}(X)\approx F(X),$

\noindent and $\kappa^{-1}(X,\delta)$ is a double cover, which by
the following lemma is identified with $\widetilde{F(X)}$.  (The
compatibility with $\lambda$ follows from (5.9.4); line incidence
in $F(X)$ corresponds by (4.6) to the tetragonal relation among
the quintics, which by figure (5.7) corresponds, in turn, to
the
tetragonal relation in ${\cal R}_5$, so the proof of (5.8)(3) is
complete.)

\bigskip

\noindent{\bf Lemma 5.10}   The Albanese double cover
$\widetilde{F(X)}$ determined by \\ $\delta\in J_2(X)$  is
isomorphic to $\pi^{-1}{\cal R P}^{-1}(X,\delta)$ (notation of
(5.9.4)).

\bigskip

\noindent{\bf Proof.}

The second isomorphism in (5.9.5) sends a line $\ell\in
F(X)$ to the object $(\widetilde{Q}_{\ell},Q_{\ell})\in{\cal
P}^{-1}(X)$, where the curves $\widetilde{Q}_{\ell},Q_{\ell}$
parametrize ordered (respectively, unordered) pairs $\ell',\ell''
\in F(X)$ satisfying: $$\ell + \ell' + \ell''= 0 \ \ \ \ \ \ \ \
\
\ {\rm (sum \ in \ } \ \ \ J(X)).$$
 We may of course think of $\widetilde{Q}_{\ell}$ as sitting in
$F(X)$, since $\ell'$ uniquely determines $\ell''$:
$\widetilde{Q}_{\ell}$ is the closure in $F(X)$ of

\noindent{\bf (5.10.1)} \ \ \ \ \ \ \ \ \ \ \ \ \ \ \ \ \
$\{ \ell'\in F(X) \ \ | \ \ \ell'\cap\ell\neq\phi  \, , \, \ell' \neq
\ell \}.
$

\medskip

\noindent
The corresponding object of ${\cal R}{\cal P}^{-1}(X,\delta)$ is
$(\stackrel{\approx}{Q}_\ell, \widetilde{Q}_\ell, Q_\ell)$, where
$\stackrel{\approx}{Q}_\ell$ is the inverse image in
$\widetilde{F(X)}$ of $\widetilde{Q}_\ell$ embedded in $F(X)$ via
(5.10.1).  Now to specify a point in $\pi^{-1}{\cal R}{\cal
P}^{-1}(X,\delta)$ we need, additionally, a double cover
$\widetilde{Q}_\ell' \to Q_\ell$ satisfying:

\noindent{\bf (5.10.2)} \ \ \ \ \ \ \ \ \ \ \ \ \ \ \ \ \ \ \ \ \
$\widetilde{Q}_\ell\times_{Q_\ell}\widetilde{Q}_\ell'\approx
\stackrel{\approx}{Q}_\ell \ .$

We need to show that a choice of $\widetilde{\ell}
\in\widetilde{F(X)}$
over $\ell\in F(X)$ determines such a $\widetilde{Q}_\ell'$.Recall
that
$\widetilde{F(X)} \to F(X)$ is obtained by base change,
via
the Albanese map, from the double cover $\widetilde{J(X)} \to
J(X)$ determined by $\delta$. \ $\widetilde{Q}_\ell'$ can
thus be taken to parametrize unordered pairs $\widetilde{\ell}',
\widetilde{\ell}''\in\widetilde{F(X)}$ satisfying: $$\widetilde{\ell}
+ \widetilde{\ell}' + \widetilde{\ell}'' = 0 \ \ \ \ \ \ \ {\rm (sum
\
in} \ \ \ \widetilde{J(X)})\ .$$ The fiber product in (5.10.2)
then parametrizes such ordered pairs, so the required
isomorphism to $\doubletilde{Q}_\ell$ simply sends
$$(\widetilde{\ell}',\widetilde{\ell}'') \mapsto \widetilde{\ell}'.$$
\begin{flushright} Q.E.D. \end{flushright}

Finally, we prove (5.8)(4).  Let $\overline{{\cal P}}:
\overline{{\cal R}}_5 \to {\cal A}_4$ be the proper Prym map.  By
(5.8)(3) it factors $$\overline{\cal P} = \iota \circ \kappa$$
where $\iota: {\cal R}{\cal C}^+ \to {\cal A}_4$ is a rational
map,
which we are trying to show is birational.  It suffices to find
some $A \in{\cal A}_4$ such that:

\begin{list}{{\rm(\arabic{bean})}}{\usecounter{bean}}
\item Any two objects in $\overline{\cal P}^{-1}(A)$ can be
related by a sequence of tetragonal constructions.
\item The differential $d{\cal P}$ is surjective over $A$.
\end{list}

In \S5.4 we see that $(1)$ holds for various examples,
including generic Jacobians $\in{\cal J}_4$: for generic $C \in
{\cal M}_4$, $\overline{\cal P}^{-1}(J(C))$ consists of
Wirtinger covers $\widetilde{C} \to C'$ (with normalilzation $C$) and
of trigonals $T$, and the two types are exchanged by $\lambda$.It is
easier to
check surjectivity of $d{\cal P}$ at the
Wirtingers: by theorem (1.6), this amounts to showing that the
Prym-canonical curve $\Psi(X) \subset {\bf P}^3$ is contained in
no quadrics.  By [DS] IV, Propo. 3.4.1, $\Psi(X)$ consists of the
canonical curve $\Phi(C)$ together with an (arbitrarily chosen)
chord.  Since $\Phi(C)$ is contained in a unique quadric $Q$,
which does not contain the generic chord, we are done.  [Another
argument: it suffices to show that no one quadric contains
$\Psi(T)$ for all trigonal $T$ in ${\cal P}^{-1}(J(C))$.By [DS], III
2.3 we
have $$\cup_T \Psi(T) \ \supset \ \Phi(C) ,$$
so the only possible quadric would be $Q$.  Consider the $g^1_4$
on
$C$ given by $\omega_C$(-$p$-$q$), where $p,q \in C$ are such
that the
chord $\overline{\Phi(p),\Phi(q)}$ is not in $Q$.  Let $T$ be the
trigonal curve associated to $(C,\omega_C$(-$p$-$q$)), and choose
a plane $A \subset {\bf P}^3$ through $\Phi(p),\Phi(q)$, meeting
$Q$ and $\Phi(C)$ transversally, say $$A \cap \Phi(C) =
\Phi(p+q+\sum^{4}_{i=1}x_i),$$ then by [DS],III 2.1, \ \
$\Psi(T)$
contains the point $$\overline{\Phi(x_1),\Phi(x_2)} \cap
\overline{\Phi(x_3), \Phi(x_4)}$$ which cannot be in
$Q$.]\begin{flushright}
Q.E.D. \end{flushright}

\subsection{Special fibers.}

\begin{tabbing} X \= \kill
\>We want to illustrate the behavior of the Prym map over some \\
special loci in $\overline{\cal A}_4$.  The common feature to all
of these examples is that \\ the cubic threefold $X$ given in
Theorem (5.1) acquires a node.  We thus \\ begin with a review of
some
results, mostly from [CG], on nodal \\ cubics.
\end{tabbing}
\bigskip

\noindent{\bf (5.11) Nodal cubic threefolds}

There in a natural correspondence between nodal cubic threefolds
$X\subset{\bf P}^4$ and nonhyperelliptic curves $B$ of genus 4.
Either object can
be described by a pair of homogeneous polynomials $F_2, \ F_3$,
of degrees 2 and 3 respectively, in 4 variables $x_1, ..., x_4:X$
has homogeneous equation $0=F_3+x_0F_2$ \ (in ${\bf P}^4$), and
the canonical curve $\Phi(B)$ has equations $F_2=F_3=0$ in ${\bf
P}^3$.

More geometrically, we express the Fano surface $F(X)$ in terms
of $B$.  Assume the two $g^1_3$'s on $B$, \ ${\cal L}'$ and
${\cal L}''$, are distinct.  They give maps $$\tau', \
\tau'':B\hookrightarrow S^2B$$ sending $r\in B$ to $p+q$ if
$p+q+r$ is a trigonal divisor in $|{\cal L}'|$, \ $|{\cal L}''|$
respectively.  We then have the identification
\noindent{\bf (5.11.1)} \ \ \ \ \ \ \ \ \ \ \ \ \ \ \ $F(X)\approx
S^2B/(\tau'(B)\sim\tau''(B))$.

Indeed, we have an embedding $$\tau:B\hookrightarrow F(X),$$
identifying $B$ with the family of lines through the node
$n = (1,0,0,0,0)$. This gives a map $S^2B\rightarrow F(X)$
sending a pair $\ell_{1}, \ell_{2}$ of lines through $n$ to the
residual
intersection with $X$ of the plane $(\ell_1, \ell_2)$.this map
identifies
$\tau'(B)$ with $\tau''(B)$, and induces the
isomorphism (5.11.1).

\bigskip

\noindent{\bf (5.11.2)}  A line $\ell\in F(X)$ determines a pair
$(Q,
\widetilde{Q})\in\overline{\cal RQ}^-$, which must be in
$\partial^{\rm II}{\cal RQ}^-$, i.e. for generic $\ell$ we obtain a
nodal quintic $Q$ with \'{e}tale double cover $\widetilde{Q}$.  We
can
interpret (5.11.1) in terms of these nodal quintics:  Start with
a
divisor $p+q\in S^2B$.  Then $\omega_B(-p-q)$ is a $g^1_4$ on
$B$,
so the trigonal construction produces a double cover
$\widetilde{T}\rightarrow T$, where $T\in{\cal M}_5$ comes with a
trigonal bundle ${\cal L}$.  The linear system
$|\omega_T\otimes{\cal L}^{-1}|$ maps $T$ to a plane quintic $Q$,
with a single node given by the divisor $|\omega_T\otimes{\cal
L}^{-2}|$ on $T$.

\bigskip

\noindent{\bf (5.11.3)}  In the special case that there exists
$r\in B$ such that \\ $p+q=\tau''(r)$, i.e. $p+q+r\in|{\cal
L}''|$ is a trigonal divisor, our $g^1_4$ acquires a base point:
$$\omega_B(-p-q)\approx{\cal L}'(r).$$ As seen in (2.10.ii), the
trigonal construction produces the nodal trigonal curve
$$T:=B/(p'\sim q')$$ with its Wirtinger double cover $\widetilde{T}$,
where $p', q'\in B$ are determined by: $$p'+q'+r\in|{\cal L}'|,$$
i.e. $p'+q'=\tau'(r)$.  In this case, the quintic $Q$ is the
projection of $\Phi(B)$ from $\Phi(r)$, with 2 nodes $p\sim q, \
\
p'\sim q'$, and $\widetilde{Q}$ is the reducible double cover with
crossings over both nodes.

\bigskip

\noindent{\bf (5.12) Degenerations in ${\cal RC}^+$.}

We fix our notation as in \S 5.1.  Thus we have:

\begin{tabbing}
\=$X\in{\cal C} \ \ \ \ \ \ \ (X, \delta)\in{\cal RC}^+$ \\
\>$(Q,\sigma)\in{\cal RQ}^-, \ \ \ \ \ \ (Q, \nu), \ (Q,
\nu\sigma)\in{\cal RQ}^+$ \\ \>$(C, \mu), \ (C', \mu')\in{\cal
R}_5$ \\ \>$A\in{\cal A}_4$
\end{tabbing}

\noindent and these objects satisfy:

\begin{tabbing}
\= ${\cal P}(Q, \sigma)$ \= $=$ \= $J(X)$ \= , \ \ \ \ \ \ \ \ \
\
\ \ \ \ \ \ \= $\nu, \nu\!$ \= $\!\sigma\mapsto\delta$ \\
\>${\cal
P}(Q, \nu)$ \>= \>$J(C)$ \>, \>$\sigma, \nu\!$
\>$\!\sigma\mapsto\mu$ \\ \>${\cal P}(Q, \nu\sigma)$ \>=
\>$J(C')$
\>, \>$\nu,$ \>$\!\sigma\mapsto\mu'$ \end{tabbing}

\bigskip

\begin{tabbing}
\= ${\cal P}(C, \mu)={\cal P}\!$ \= $\!(C'\!$ \= $, \mu'\!$ \=
$\!)=A$ \\ \>$\lambda(C, \mu)=$ \>$\!(C'\!$ \>$, \mu'\!$ \>$\!)$
\\
\>$\kappa(C, \mu)=$ \>$\!(X$ \>$, \delta\!$ \>$\!).$
\end{tabbing}

Now let $X$ degenerate, acquiring a node, with
$\bar{\varepsilon}\in J_2(X)\backslash(0)$ the vanishing cycle
mod. 2. \ From (5.11) we see that $Q$ also degenerates, with a
vanishing cycle $\varepsilon$ which maps (via. (1.4)) to
$\bar{\varepsilon}$. Lemma (5.9) of [D4] shows that
$\varepsilon$, hence also $\bar{\varepsilon}$, must be even.

There are 3 types of degenerations of $(X, \delta)$,
distinguished
as in (1.7) by the relationship of $\delta, \bar{\varepsilon}$.
(A fourth type, where $Q$ degenerates but $X$ does not, is
explained in (5.13).)The possibilities are summarized below:

\bigskip

\begin{list}{{\rm(\Roman{butter})}}{\usecounter{butter}}
\item If $\bar{\varepsilon}=\delta$ then either
$\varepsilon=\nu$ or $\varepsilon=\nu\sigma$, which gives the same
picture with $C, C'$ exchanged.  In case $\varepsilon = \nu$, $(Q,
\nu)$
undergoes
a $\partial^{\rm I}$ degeneration, while $(Q, \nu\sigma)$ is
$\partial^{\rm II}$. (The notation is that of (1.7).) Thus $A$ is a
Jacobian.

\ \ \ The double cover $\widetilde{F(X)}$ is itself a $\partial^I$
cover.  In terms of the curve $B$ of (5.11), we have
$$\widetilde{F(X)}=(S^2B)_0\amalg(S^2B)_1 \ \ / \ \
(\tau'(B)_0\sim\tau''(B)_1, \ \ \ \tau''(B)_0\sim\tau'(B)_1).$$This
is clear,
either from the definition of $\widetilde{F(X)}$ viathe Albanese map,
or by
considering the restriction to ${\cal
RP}^{-1}(X, \delta)$ of the double cover $${\cal R}^2{\cal
Q}^+\stackrel{\pi}{\rightarrow}{\cal R}^2{\cal Q}^-$$ of (5.9).
One of the components parametrizes the trigonal objects $(C,
\mu)$, the other parametrizes the nodals $(C', \mu')$.

\item $\sigma$ is always perpendicular to $\varepsilon,
\nu$, and the condition $\langle\bar{\varepsilon}, \delta\rangle=0$
implies
$\langle\varepsilon, \nu\rangle=0$ by (1.4.3).  Both $(Q, \nu)$ and
$(Q,
\nu\sigma)$ then give $\partial^{\rm II}$-covers, so $C, C'$ are
nodal. Again by (1.4.3), both $(C, \mu)$ and $(C', \mu')$ are
$\partial^{\rm II}$, so their common Prym $A$ is in
$\partial\bar{\cal
A}_4$.

\ \ \ \ From the Albanese map we see that $\widetilde{F(X)}$ is an
etale cover of $F(X)$.  Indeed, $\delta$ comes from a semiperiod
$\delta'$ on $B$, giving a double cover $\widetilde{B}$ with
involution $\iota$; the normalization of $\widetilde{F(X)}$ is then
$S^2\widetilde{B}/\iota$, and $\widetilde{F(X)}$ is obtained by
glueing
above $\tau(B)$.

\item In this case both $(Q, \nu)$ and $(Q, \nu\sigma)$ are
$\partial^{\rm III}$, so $C, C'$ are nonsingular.  The node of $Q$
represents a quadric of rank 3 through $\Phi(C)$, so ${\cal L}$
is
cut out by the unique ruling.  By the Schottky-Jung relations
[M2],
the vanishing theta null on $C$ descends to one on $A$.

The double cover $\widetilde{F(X)}$ is again a
$\partial^{\rm III}$-cover, in the sense that its normalization is
ramified over $\tau'(B), \tau''(B)$, the sheets being glued.Each of
the
quintics in (5.11.2) gives two points of
$\widetilde{F(X)}$, while the two-nodal quintics (5.11.3) land in the
branch locus of $\pi$ \ (5.9.4).

\end{list}

\begin{center}
\begin{tabular}{|l|l|l|l|l|}
\hline
Degeneration & Degeneration &  & & \\
type of  & type of & & & \\
\multicolumn{1}{|c|}{$(X,\delta)$} &
\multicolumn{1}{|c|}{$(Q,\sigma,\nu,\nu\sigma)$}
& \multicolumn{1}{|c|}{$(C,\mu)$} & \multicolumn{1}{|c|}{$(C',\mu')$}
&
\multicolumn{1}{|c|}{$A$} \\ \hline \hline
 & & & & \\
I : $\bar{\varepsilon} = \delta$ &$\varepsilon = \nu$ &
nonsingular  & nodal, $\partial^{\rm I}$ & ${\cal J}_{4}$ \\
 & & trigonal &  &  \\ \hline
  & & & &  \\
 II :  $\bar{\varepsilon} \neq \delta$, & $(\varepsilon, \sigma,
\nu)$  &
 nodal, $\partial^{\rm II}$& nodal, $\partial^{\rm II}$& $\partial
\bar{\cal
A}_{4}$ \\
 \multicolumn{1}{|c|}{$\langle \bar{\varepsilon} , \delta \rangle =
0$} & rank
3 &
  &   &
  \\
   & isotropic & & & \\
    & subgroup & & & \\ \hline
    & & & & \\
  III : $\langle \bar{\varepsilon} , \delta \rangle \neq 0$ &
 $\langle \varepsilon , \sigma \rangle = 0$ & nonsingular, &
  nonsingular, & $\theta_{\rm null}$ \\
   &$\langle \varepsilon , \nu \rangle \neq 0$ & has vanishing &
   has vanishing & \\
    & & thetanull ${\cal L}$, & thetanull ${\cal L}'$, & \\
     & & ${\cal L}(\mu)$ even & ${\cal L}'(\mu')$ even  & \\ \hline
   IV : & $\langle \varepsilon , \nu \rangle = 0$ &
    & nonsingular, &  \\
   nonsingular & $\langle \varepsilon , \sigma \rangle \neq 0$ &
nodal,
   $\partial^{\rm II}$ &
   has vanishing & ${\cal A}_{4}$ \\
    & & & thetanull ${\cal L}'$, & \\
     & & & ${\cal L}'(\mu')$ odd & \\ \hline
      \end{tabular}
 \end{center}

\bigskip

\pagebreak[4]

\noindent{\bf (5.13) Degenerations in ${\cal R}^2{\cal Q}^+$.}

We have just described the universe as seen by a degenerating
cubic threefold.  From the point of view of a degenerating plane
quintic, there are a few more possibilities though they lead to
no new components.  We retain the notation:  $Q, \nu, \sigma,
\varepsilon$ etc.

\noindent 0. \ \ $\varepsilon$ cannot equal $\sigma$, since
$\varepsilon$ is even, $\sigma$ odd.

\noindent I. \ \ $\varepsilon=\nu$ reproduces case I of (5.12), as
does:

\noindent I$'$. \ $\varepsilon=\nu\sigma$.

\begin{tabbing}
\noindent II. \ \=Excluding the above, $\nu, \sigma,
\varepsilon$ generate a subgroup of rank 3.  If \\ \>it is
isotropic, we are in case II above.
\end{tabbing}

\noindent III. If $\langle\varepsilon, \sigma\rangle=0$ but
$\langle\varepsilon,
\nu\rangle=\langle\varepsilon, \nu\sigma\rangle\neq0$, we're in case
III.

The only new cases are thus:
\noindent IV. \ \ \ \ \ \ \ \ \ \ \ \ \ $\langle\varepsilon,
\nu \ \ \rangle=0\neq\langle\varepsilon, \sigma\rangle, \ \ {\rm
or}:$

\noindent IV.$'$ \ \ \ \ \ \ \ \ \ \ \ \ \ $\langle\varepsilon,
\nu\sigma\rangle=0\neq\langle\varepsilon, \sigma\rangle,$
which is the same as IV after exchanging $C, C'$.

In case IV, we find:

\noindent $\bullet$ $X$ is non-singular, in fact any $X$ can
arise.  What is special is the line $\ell\in F(X)$ corresponding
to $Q:$ it is contained in a plane which is tangent to $X$ along
another line, $\ell'$.

\noindent $\bullet$ $(Q, \nu)$ is a $\partial^{\rm II}$ degeneration,
so
$C$ is nodal, and $(C, \mu)$ is a $\partial^{\rm II}$ degeneration.

\noindent $\bullet$ On the other hand, $(Q, \nu\sigma)$ is
$\partial^{\rm III}$, so $C'$ is non-singular, and has a vanishing
theta null ${\cal L}'$ (corresponding, as before, to the node of
$Q$).

\noindent $\bullet$ This time though, ${\cal L}'(\mu')$ is odd,
so
$A\in{\cal A}_4$ does not inherit a vanishing theta null.  In
fact,
any $A\in{\cal A}_4$ arises from a singular quintic with
degeneration of type IV.

So far, we found three loci in $\bar{\cal A}_4$ which are related
to nodal cubics:

\[
\begin{array}{lcl}
{\cal P}\circ \kappa^{-1}(\partial^{\rm I} {\cal RC}^{+}) & \subset &
{\cal
J}_{4} \\
 \bar{\cal P}\circ \kappa^{-1}(\partial^{\rm II} {\cal RC}^{+}) &
\subset &
 \partial \bar{A}_{4} \\
 {\cal P}\circ \kappa^{-1}(\partial^{\rm III} {\cal RC}^{+}) &
\subset &
\theta_{\rm
 null}
 \end{array}
 \]

 We are now going to study, one at a time, the fibers of ${\cal P}$
above
 generic points in these three loci. We note that related results
have
 recently been obtained by Izadi. In a sense, her results are more
precise:
 she knows (cf. Remark 5.4) that $\chi$ is an isomorphism on the open
 complement ${\cal U}$ of a certain 6-dimensional locus in ${\cal
A}_{4}$.
 In [I] she shows that for $A \in {\cal U}$, $\chi (A)$ is singular
if and only
if
 \[ A \in {\cal J}_{4} \cup \theta_{\rm null}. \]
 Her description of the cubic threefold corresponding to $A \in {\cal
J}_{4}$
 complements the one we give below. In general her techiques, based
on
 $\Gamma_{00}$, are very different than our degeneration arguments.
\bigskip

\noindent{\bf (5.14) Jacobians}

\bigskip
\noindent{\bf Theorem 5.14}  Let $B\in{\M}_4$ be a general curve
of genus 4, and let $(X,\delta)=\chi(J(B)).$
\begin{list}{{\rm(\arabic{bean})}}{\usecounter{bean}}
\item $X$ is the nodal cubic threefold corresponding to $B$
(5.11).
\item $(X,\delta)\in\partial^I$, so $\widetilde{F(X)}$ is reducible,
each component is isomorphic to $S^2B.$
\item Let $(Q,\sigma,\nu)$ be the plane quintic with rank-2
isotropic subgroup corresponding to some $\ell\in F(X)$.  Then
$Q$
is nodal, with trigonal normalization $T, \ \nu$ is the vanishing
cycle, and \\ $(Q,\sigma)=(Q,\nu\sigma)\in \partial^{\rm II}$.
\item$\bar{\p}^{-1}(J(B))$ is isomorphic to $\widetilde{F(X)}$.The
component
corresponding to $\nu$ (respectively $\nu\sigma$)
consists of trigonal curves $T_{p,q}$ (respectively Wirtinger
covers of singular curves $S_{p,q}$), \ \ \ $(p,q)\in S^2B.$
\item The tetragonal construction takes both $S_{p,q}$ and
$T_{p,q}$ to $S_{r,s}$ and $T_{r,s}$ if and only if $p+q+r+s$ is
a
special divisor on $B$.  The involution $\lambda$ exchanges
$S_{p,q}, T_{p,q}$.
\item Any two objects in $\bar{\p}^{-1}(J(B))$ can be connected
by a sequence of two tetragonal moves (generally, in 10 ways).
\end{list}

\bigskip
\noindent{\bf Proof}

Since at least some of these results are needed for the proof of
(5.8)(4), we do not use Theorem (5.3).  For $(p,q)\in S^2B$, we
consider:

\noindent $\bullet$  $\widetilde{T}_{p,q}\to T_{p,q}$, the trigonal
double cover associated by the trigonal construction to $B$ with
the $g^1_4$ given by $\omega_B$(-$p$-$q$).

\noindent $\bullet$  $\widetilde{S}_{p,q}\to S_{p,q}$, the Wirtinger
cover of $S_{p,q}:=B/(p\sim q)$.  (When $p=q$, this specializes
to
$B\cup_pR$, where $R$ is a nodal rational curve in which $p$
is a non singular point.)

These objects are clearly in $\bar{\p}^{-1}(J(B))$.  Beauville's
list ([B1], (4.10)) shows that they exhaust the fiber.  This
proves
part (4).  Now clearly $\kappa$, as defined in (5.1.6), takes any
of these objects to our $(X,\delta)$; so the analysis in
(5.12)(I)
applies, proving (1)-(3).  (Note:  this already suffices to
complete the proof of (5.8)(4)!)

Let $r+s+t+u$ be an arbitrary divisor in
$|\omega_B(-p-q)|$.Projection of
$\Phi(B)$ from the chord
$\overline{\Phi(t),\Phi(u)}$
gives (the general) $g^1_4$ on $S_{p,q}$.  The tetragonal
construction takes this to the curves $T_{t,u}$ and $S_{r,s}$.(The
situation is
that of (2.15.2).)

On $T_{p,q}$ there are two types of $g^1_4$'s, of the form ${\cal
L}(x)$ and $\omega\otimes{\cal L}^{-1}(-x)$, where ${\cal L}$ is
the trigonal bundle and $x\in T_{p,q}$.  Now $x$ corresponds to a
(2,2) partition, say $\{\{r,s\},\{t,u\}\}$, of some divisor in
$|\omega_B(-p-q)|$.  The tetragonal construction,
applied to ${\cal L}(x)$, yields the curves $S_{r,s}$ and
$S_{t,u}$; while when applied to $\omega\otimes{\cal
L}^{-1}(-x)$,
it gives $T_{r,s}$ and $T_{t,u}$.  Altogether, this proves (5).We
conclude
with:

\bigskip

\noindent{\bf Lemma 5.14.7}  Given any $p,q,r,s\in B$, there are
points $t,u\in B$ (in general, 5 such pairs) such that both
$p+q+t+u, r+s+t+u$ are special.

\bigskip

\noindent{\bf Proof}

Let $\alpha, \beta$ be the maps of degree 4 from $B$ to ${\bf
P}^1$ given by \\ $|\omega_B(-p-q)|, |\omega_B(-r-s)|$.
Then $$\alpha\times\beta:B\to{\bf P}^1\times{\bf P}^1$$
exhibits $B$ as a curve of type (4,4) on a non-singular quadric
surface, hence the image has arithmetic genus
$(4-1)^2=9\rangle4=g(B)$,
so there must be (in general, 5) singular points; these give the
desired pairs $(t,u)$. \begin{flushright} QED \end{flushright}

\bigskip

\noindent{\bf  (5.15) The Boundary.}

The results in this case were obtained by Clemens [C2].
A general point $A$ of the boundary $\partial\bar{\A}_4$ of
a toroidal compactification $\bar{\A}_4$ is a ${\bf
C}^*$-extension
of some $A_0\in{\A}_3$.  The extension data is given by a point
$a$
in the Kummer variety $A_0/(\pm 1)$.

Given $a\in A_0$, consider the curve $$\widetilde{B}=
\widetilde{B}_a:=\Theta\cap\Theta_a\subset A_0$$ (where
$x\in\Theta_a\Leftrightarrow x+a\in\Theta$), and its quotient
$B=B_a$ by the involution $x\mapsto -a-x$.  We have
$$(B,\widetilde{B})\in{\R}_4$$ and $${\p}(B,\widetilde{B})\approx
A_0.$$  The pair $(B,\widetilde{B})$ does not change (up to
isomorphism) when $a$ is replaced by $-a$.

\bigskip

\noindent{\bf Theorem 5.15 ([C2])}  Let
$A\in\partial\bar{\A}_4$ be the ${\bf C}^*$-extension of
\\ $A_0\in{\A}_3$, a generic  $P\!P\!A\!V$, determined by $\pm
a\in A_0$.  Let $(X, \delta)=\chi(A)$.
\begin{list}{{\bf(\arabic{bean})}}{\usecounter{bean}}
\item $X$ is the nodal cubic threefold corresponding to $B=B_a$.
\item $(X,\delta)\in\partial^{\rm II}$, so $\widetilde{F(X)}$ is the
etale
double cover of $F(X)$ with normalization $S^2\widetilde{B}/\iota$,
as
in (5.12.II).
\item The corresponding quintics $Q$ are nodal; all three of
$\sigma,\nu,\nu\sigma$ are of type $\partial^{\rm II}$.
\item $\doublebar{\p}^{-1}(A)$ is isomorphic to $\widetilde{F(X)}$,
and
consists of $\partial^{\rm II}$-covers $(C,\widetilde{C})$ whose
normalizations (at one point) are of the form $(B_b,\widetilde{B}_b)$
for $b=b_1-b_2, \ \ b_1, b_2\in\psi(\widetilde{B})$.
\end{list}

\bigskip

\noindent{\bf Proof}

Clearly $\doublebar{\p}^{-1}(A)\subset\partial^{\rm
II}\doublebar{\R}_5$, so
consider a pair $(C,\widetilde{C})\in\partial^{\rm II}$, say
$$C=N/(p\sim
q),\;\;\;\widetilde{C}=\widetilde{N}/(p'\sim q', p"\sim q")$$  with
$(N,\widetilde{N})\in\bar{\R}_4$.  Then
$\doublebar{\p}(C,\widetilde{C})$ is a
${\bf C}^*$-extension of $P(N,\widetilde{N})$, with extension data
$$\pm(\psi(p')-\psi(q'))\in{\p}(N, \widetilde{N})/(\pm 1).$$  We see
that $\doublebar{\p}(C,\widetilde{C})=A$ if and only if
\noindent{\bf (5.15.5)} \ \ \ \ \ \ \ \ \ \ \ \ \ \ \ \ \ \ \
$(N,\widetilde{N})\in\bar{\p}^{-1}(A_0),$

\noindent and:

\noindent{\bf (5.15.6)} \ \ \ \ \ \ \ \ \ \ \ \ \ \ \
$\psi(p')-\psi(q')=a, \ \ \ p',q'\in\widetilde{N}.$

\medskip

Now, (5.15.5) says that $(N,\widetilde{N})$ is taken, by its
Abel-Prym map $\psi$, to $(B_b,\widetilde{B}_b)$ for some $b\in A_0$,
and then (5.15.6) translates to: $$a=a_1-a_2, \ \ \ \ \ \ \ \ \ \
\ a_1,a_2\in\Theta\cap\Theta_b$$ which is equivalent to$$b=b_1-b_2, \
\ \ \ \ \
\ \ \ b_1,b_2\in\Theta\cap\Theta_a$$
(take $b_1=a_2+b, \ b_2=a_2)$.  This proves (4), and everything
else follows from what we have already seen. \begin{flushright}
QED \end{flushright}

\bigskip

\noindent{\bf (5.16) Theta nulls}

Let $A\in{\A}_4$ be a generic $PPAV$ with vanishing thetanull,
and $(C,\widetilde{C})$ a generic element of ${\p}^{-1}(A)$.  By
[B1], Proposition (7.3), $C$ has a vanishing thetanull.  This
implies that the plane quintic $Q$ parametrizing singular
quadrics through $\Phi(C)$ has a node, corresponding to the
thetanull.  The corresponding cubic threefold $X$ is thus also
nodal, and we are again in the situation of (5.11.III).  I do not
see, however, a more direct way of describing the curve $B$ (or
the cubic $X$) in terms of $A$.

\bigskip

\noindent{\bf (5.17) Pentagons and wheels}.

In [V], Varley exhibits a two dimensional family of double covers
$(C,\widetilde{C})\in{\R}_5$ whose Prym is the unique
non-hyperelliptic
$P\!P\!A\!V$ \linebreak $A\in{\A}_4$ with 10 vanishing thetanulls.
The
curves
$C$ involved are Humbert curves, and each of these comes with a
distinguished double cover $\widetilde{C}$. As an illustration of our
technique, we work out the fiber of $\doublebar{\p}$ over $A$ and the
tetragonal moves on this fiber.  This is, of course, a very
special
case of (5.12)(III) or (5.16).

We recall the construction of Humbert curves and their double
covers.  Start by marking 5 points $p_1,\cdots,p_5\in{\bf P}^1$.Take
5 copies
$L_i$ of ${\bf P}^1$, and let $E_i$ be the double
cover of $L_i$ branched at the 4 points $p_j, \ \ j\ne i$.  Let
\noindent{\bf (5.17.1)} \ \ \ \ \ \ \ \ \ \ \ \ \ \ \ \ \ \ \ \ \ \
\ \ \ $A:=\coprod^5_{i=1}L_i, \;\;\;\;\; B:=\coprod^5_{i=1}E_i.$

\noindent The \underline{pentagonal} construction applied to

\noindent{\bf (5.17.2)} \ \ \ \ \ \ \ \ \ \ \ \ \ \ \ \ \ \ \ \ \ \
\  \ \ \ \ \ \ \ $B\stackrel{g}{\to}A\stackrel{f}{\to}{\bf P}^1$

\noindent($f$ is the forgetful map, of degree 5), yields a
32-sheeted branched cover $f_*B\to{\bf P}^1$ which splits, by
(2.1.1), into 2 copies of the Humbert curve $C$, of degree 16
over ${\bf P}^1$.

Let ${\beta}_I$, $I\subset S:=\{1,\cdots,5\}$, be the
involution of (5.17.2) which fixes $A$ and acts non-trivially on
$E_i, \ \ i\in I$.  It induces an involution $\alpha_I$ on
$f_*B$, hence on its quotient $C$.  Let $$G:=\{\alpha_I \ \ | \ \
I\subset S\} \ / \ (\alpha_S).$$  Then $C$ is Galois over ${\bf
P}^1$, with group $G\approx({\Z}/2{\Z})^4$.  Let $G_i, \ \ \ 1\le
i\le 5$, be the image in $G$ of $$\{\alpha_I|i  \not\in I, \ \
\#(I)=\;\mbox{even}\}.$$Then $$C/\alpha_i\approx C/G_i\approx E_i,$$

and the
quotient map
$$E_i\approx C/\alpha_i\to C/G_i\approx E_i$$ becomes
multiplication by 2 on $E_i$.  In particular, the Humbert curve
$C$ has 5 bielliptic maps $h_i:C\to E_i$.  The branch locus of
$h_i$
consists of the 8 points $x\in E_i$ satisfying $g(2x)=p_i$.

For ease of notation, set $E:=E_5, \ \ \ p=p_5\in{\bf P}^1$,
$$C\stackrel{h}{\to}E\stackrel{g}{\to}{\bf P}^1,$$ and
$$\{p^0,p^1\}:=g^{-1}(p)\subset E.$$  Then for $j=0,1, \ \ E$ has
a natural double cover $C^{j}$, branched at the four points
${1\over 2}p^j$ and given by the line bundle ${\cal O}_E(2p^j)$.The
fiber
product
\noindent{\bf (5.17.3)} \ \ \ \ \ \ \ \ \ \ \ \ \ \ \ \ \ \ \ \ \ \
\ \ \ \ \ \ \ \ \ \ \ $\widetilde{C}:=C^0\times_EC^1$

\noindent gives a Cartesian double cover of $C$.

Replacing $E_5$ by another $E_i$, we get an isomorphic double
cover
$\widetilde{C}$.  Here is an invariant description of this cover:

Let $p_{i,j}:=L_i\cap f^{-1}(p_j)\in A$, and consider the curve
$$Q:=A/(p_{i,j}\sim p_{j,i}, \ \ \ i\ne j).$$  Then $Q$ can be
embedded in ${\bf P}^2$ as a pentagon, or completely reducible
plane quintic curve: embed ${\bf P}^1$ as a non-singular conic,
and take $L_i$ to be the tangent line of the conic at $p_i$.  We
have two natural branched double covers of $Q$:
\noindent{\bf (5.17.4)} \ \ \ \ \ \ \ \ \ \ \ \ \ \ \
$\widetilde{Q}_\sigma:=(\coprod^{ \: 5 \ \ \ 1} _{i=1,\varepsilon=0}
\ \ L^\epsilon_i)/(p^0_{i,j}\sim p^1_{j,i}, \ \ \ i\ne j)$

\noindent{\bf (5.17.5)} \ \ \ \ \ \ \ \ \ \ \ \ \ \ \ \ \ \ \ \ \ \
\ \
$\widetilde{Q}_\nu:=B/(\widetilde{p}_{i,j}\sim\widetilde{p}_{j,i}, \
\
\ i\ne j),$

\noindent where $\widetilde{p}_{i,j}\in E_i$ is the unique
(ramification) point above $p_{i,j}\in L_i$.  We may think of
$\widetilde{Q}_\sigma$ as a "totally $\partial^I$" degeneration, and
of
$\widetilde{Q}_\nu$ as a "totally $\partial^{\rm III}$" degeneration.
We
then find:

\begin{tabbing}
\noindent{\bf (5.17.6)} \=$(Q,\widetilde{Q}_\nu)\in\overline{\R\Q}^+$
is the quintic double cover corresponding to the \\ \>Humbert
curve $C\in{\M}_5$.
\end{tabbing}

\begin{tabbing}
\noindent{\bf (5.17.7)}  \=The double cover $\widetilde{Q}_\sigma$ of
$Q$ corresponds, via (1.4.2), to the double \\ \>cover
$\widetilde{C}$ of $C$.
\end{tabbing}

We note that $\widetilde{Q}_\sigma$ is itself an odd cover, so it
corresponds to some (singular) cubic threefold.  A moment's
reflection shows that this must be Segre's cubic threefold $Y$ which
 we have
already met in (4.8). Indeed, the Fano surface $F(Y)$ consists of the
six rulings $R_i, \ \ \
0\le i\le 5$, plus the 15 dual planes $\Pi_{i,j}^*$ of lines in
$\Pi_{i,j}$ (notation of (4.8)).  We see that:

\begin{tabbing}
\noindent{\bf (5.17.8)}  \=The discriminant of projection of $Y$
from a line $\ell\in R_i$ is a plane \\ \>pentagon $Q$, with its
double cover $\widetilde{Q}_\sigma$ as above.
\end{tabbing}

The other covers, $\widetilde{Q}_\sigma$, fit together to determine a
point $(Y,\delta)\in\overline{\R\C}^+$:

\noindent{\bf (5.17.9)} \ \ \ \ \ \ \ \ \ \ \ \ \ \ \ \ \ \ \ \ \ \
\ \ $(Y,\delta)=\kappa(C, \widetilde{C}),$

\noindent for any Humbert cover $(C,\widetilde{C})$.  The
tetragonal construction takes any $(Q,\widetilde{Q}_\sigma)$ to any
other (in two steps), so we recover Varley's theorem:

\noindent{\bf (5.17.10)} $A:={\p}(C, \widetilde{C})\in{\A}_4$ is
independent of the Humbert cover $(C,\widetilde{C})$.

But this is not the complete fiber:  we have only used one of the
two component types of $F(Y)$.  We note:

\noindent{\bf (5.17.11)}  The discriminant of projection of $Y$
from a line $\ell\subset\Pi_{ij}$ consists of a conic plus
three lines meeting at a point; the double cover is split.

\vspace{2in}

\begin{center}
\begin{tabular}{cc}
\hspace{2.3in} & \hspace{2.3in} \\
pentagon & wheel
\end{tabular}
\end{center}

Consider a tritangent plane, meeting $Y$ in lines $\ell_i\in R_i,
\ \ \ \ell_j\in R_j$, and $\ell_{ij}\in\Pi_{ij}^*$.  It
corresponds to a tetragonal construction involving two pentagons
and a wheel.  The other kind of tritangent plane intersects $Y$
in
lines $\ell_{ij}\in\Pi_{ij}^*, \ \ell_{kl}\in\Pi_{kl}^*,
\\ell_{mn}\in\Pi_{mn}^*$, where $\{i,j,k,l,m,n\}=\{0,1,2,3,4,5\}$;
the tetragonal construction then relates three wheels.

\noindent{\bf Theorem 5.18}  Let $A\in{\A}_4$ be the
non-hyperelliptic $P\!P\!A\!V$ with 10 vanishing thetanulls.
\begin{list}{{\bf(\arabic{bean})}}{\usecounter{bean}}
\item $\chi(A)$ consists of the Segre cubic threefold $Y$, with
its
degenerate semi-period $\delta$ (5.17.9).
\item The corresponding curve $B\in\bar{\M}_4$ (5.11) consists of
six ${\bf P}^1$'s:

{ \ }

\vspace{2in}

\item The Fano surface $F(Y)$ consists of the 6 rulings
$R_i \ \ (0\le i\le 5)$ and the 15 dual planes $\Pi_{i,j}^*$.The
plane quintics are pentagons, for $\ell\in R_i$, and wheels, for
$\ell\in\Pi_{ij}^*$, \ all with split covers $\sigma$ \ (5.17.4,
5.17.11).  (The $\nu$ covers are branched over all the double
points.)
\item The fiber $\doublebar{\p}^{-1}(A)$ is contained in the fixed
locus
of the involution $\lambda:\doublebar{\R}_5\to\doublebar{\R}_5$
(5.1.6), so
it is a quotient of $F(Y).$
\item $\doublebar{\p}^{-1}(A)$ consists of two components:

\begin{itemize}
\item Humbert double covers $\widetilde{C} \rightarrow C$ (5.17.3).
\item Allowable covers $\widetilde{X}_{0} \cup \widetilde{X}_{1}
\rightarrow
X_{0} \cup X_{1}$, where $X_{0}$, $X_{1}$ are elliptic, meeting at
their
4 points of order 2.
\end{itemize}\end{list}

All of this follows from our previous analysis, except (5).  The
new, allowable, covers are obtained by applying Corollary (3.7),
with $n=3$, to the Cartesian cover $\widetilde{C}\to C$ in (5.17.3).
It is also easy to see that the plane quintic parametrizing
singular quadrics through the canonical curve $\Phi(X_0\cup X_1)$
is a wheel, and vice versa, that the generalized Prym of any
wheel

(with its $\partial^{\rm III}$-cover) is the generalized Jacobian
$J(X_0\cup X_1)$ of such a curve.  Thus every line in $F(Y)$ is
accounted for, so we have the complete fiber
$\doublebar{\p}^{-1}(A)$.
\begin{flushright} QED \end{flushright}
\large
\section{Other genera}

\ \ \ \ For $g\le 4$, it is relatively easy to describe the fibers
of ${\p}:\bar{\R}_g\to{\A}_{g-1}$.  Indeed, every curve in
${\M}_g$ is trigonal, and every $A\in{\A}_{g-1}$ is a Jacobian (of a
possibly reducible curve), so the situation is completely
controlled by Recillas' trigonal construction.  Similar results can
be obtained, for \linebreak $g\le 3$, by using Masiewicki's criterion
[Ma].

\bigskip

\noindent{\bf (6.1)} \ \underline{$g=1$}.  Here $\bar{\p}$ sends
$\bar{\R}_1\approx{\bf P}^1$ to ${\A}_0$ (= a point).  Thefibers
of $\bar{\p},{\p}$ are then ${\bf P}^1, {\bf C}^*$ respectively.

\bigskip

\noindent{\bf (6.2)} \ \underline{$g=2$}.  All curves of genus 2 are
hyperelliptic, and all covers are Cartesian (3.2).  An element of
${\R}_2$ is thus given by 6 points in ${\bf P}^1$, with 4 of them
marked, modulo ${\bf P}GL(2)$; an element $E$ of ${\A}_1$ is given
by 4 points of ${\bf P}^1$ modulo ${\bf P}GL(2)$; and ${\p}$
forgets the 2 unmarked points.  The fiber of ${\p}$ is thus
rational; it can be described as $S/G$ where
\[ S := S^{2}\left( {\bf P}^{1} \setminus ({\rm 4 \; points}) \right)
\setminus
({\rm diagonal})\]
and $G \approx ({\bf Z}/2{\bf Z})^{2}$ is the Klein group, whose
action on $S$ is induced from its action on ${\bf P}^{1}$ permuting
the
4 marked points.

We note that $S$ is ${\bf P}^{2}$ minus a conic $C$ and four lines
$L_{i}$ tangent to it. To compactify it  we add:

$\bullet$  a $\partial^{\rm I}$ cover for each point of
$C\backslash\cup L_i$,

$\bullet$  a $\partial^{\rm III}$ cover for each point of
$L_i\backslash C$, and

$\bullet$  an "elliptic tail" cover [DS, IV 1.3] for each point in
the exceptional divisor obtained by blowing up one of the points
$L_i\cap C$.  (The
limiting double cover obtained is
\[{(E_0\amalg E_1)/\approx} \; \longrightarrow \; {E/\sim}\]
where $\sim$ places a cusp at one of the four marked points
$p_i$ on $E$ and $\approx$ places a
tacnode above it.  These curves are unstable, and the family of
elliptic-tail covers gives their stable models, each elliptic tail
being blown down to the cusp.)

The resulting $\overline{S}$ is ${\bf P}^{2}$ with 4 points in
general position blown up, and the compactified fiber is
$\overline{S}/G$,
or ${\bf P}^{2}/G$ with one point blown up.

\bigskip

\noindent{\bf (6.3)} \ \underline{$g=3$}.  Fix $A\in{\A}_2$.  The
Abel-Prym map sends pairs $(C,\tilde{C})\in{\cal P}^{-1}(A)$ to
curves $\psi(\tilde{C})$ in the linear system $|2\Theta|$ on $A$,
uniquely defined modulo translation by the group
$G=A_2\approx({\Z}/2{\Z})^4$.  The fiber is therefore, birationally,
the quotient ${\bf P}^3/G$.  Since some curves in $|2\Theta|$ are
not stable, some blowing up is required to obtain the biregular
model of $\bar{\p}^{-1}(A)$.  This is carried out in [Ve].  The
quotient ${\bf P}^3/G$ is identified with Siegel's modularquartic
threefold, or the minimal compactification $\bar{\A}_2^{(2)}$ of
the moduli space of $P\!P\!A\!V$'s with level-2 structure.  To
obtain $\bar{\p}^{-1}(A)$, Verra shows that we need to blow
$\bar{\A}_2^{(2)}$ up at a point $A'$, corresponding to a level-2
structure on $A$ itself, and along a rational curve.  The 2
exceptional divisors then parametrize hyperelliptic and
elliptic-tail covers, respectively.

\bigskip

\noindent{\bf (6.4)} \ \underline{$g=4$}.

As we noted in (5.15), the fiber ${\p}^{-1}(A), \ \ A\in{\A}_3$,
consists of
covers $(B_a,\tilde{B}_a), \ \ a\in A/(\pm 1):$
$$\tilde{B}_a=\Theta\cap\Theta_a,\;\;\;B_a=\tilde{B}_a/(x\sim(-
a-x)).$$  The fiber is thus (birationally) the Kummer variety
$A/(\pm 1)$.

\bigskip

\noindent{\bf (6.5)}  \underline{$g\ge 7$}.

In this case, it was proved in [FS], [K], and [W], that ${\p}$ is
generically injective.  The results in \S3 show that it is never
injective: on the hyperelliptic loci there are positive-dimensional
fibers, and various coincidences occur on the bielliptic loci. In
[D1] we conjectured:
\bigskip

\noindent{\bf Conjecture 6.5.1}  Any two objects in a fiber of
${\p}$ are connected by a sequence of tetragonal constructions.

We state this for ${\cal P}$, rather than $\overline{\cal P}$, since
      various other phenomena can contribute to non-trivial fibers at
      the boundary. For example, all fibers of $\overline{\cal P}$ on
      $\partial^{\rm I}$ are two-dimensional. On the other hand, from
       the local pictures (2.14) it is clear that the tetragonal
construction can take a nonsingular curve to a singular one. In fact
proposition (3.8) shows that it is possible for two objects in
${\R}_g$ to be tetragonally related through an intermediate object of
$\partial{\R}_g$, so some care must be taken in clarifying which
class of tetragonal covers should be allowed.  The conjecture is
consistent with our results for $g\le 6$.  For $g\ge 13$, Debarre
[Deb2] proved it for curves which are neither hyperelliptic,
trigonal, or bielliptic.  Naranjo [N] extended this to generic
bielliptics, $g\ge 10$.  The following result was communicated to
me by Radionov:

\noindent{\bf Theorem 6.5.2} [Ra] \ \ For $g\ge 7,{\R}_g^{\rm Tet}$
is
an irreducible component of the noninjectivity locus of the Prym
map, and for generic $(C,\tilde{C})\in{\R}_g^{\rm Tet}$,
${\p}^{-1}({\p}(C,\tilde{C}))$
consists precisely of three tetragonally related objects.

\pagebreak[4]

\centerline{\bf REFERENCES}

\medskip

\def\bib#1{\noindent\hbox to50pt{[#1]\hfil}\hang}
\def\bibline#1{\bib{#1}\vrule height.1pt width0.75in depth.1pt
\/,}

\vskip 25pt
\parindent=50pt
\frenchspacing

\bib{ACGH}E. Arbarello, M . Cornalba, P. Griffiths, J. Harris,
{\it Geometry of algebraic curves}, Vol. I, Springer-Verlag, New
York (1985).

\bib{B1}A. Beauville, {\it Prym varieties and the Schottky
problem}, Inv. Math. 41 (1977), 149-196.

\bib{B2}A. Beauville, {\it Sous-vari\'{e}t\'{e}s sp\'{e}ciales
des vari\'{e}t\'{e}s de Prym}, Compos. Math. 45, 357-383 (1982).

\bib{C1}H. Clemens, {\it Double Solids}, Advances in Math. 47
(1983) pp. 107-230.

\bib{C2}H. Clemens, {\it The fiber of the Prym map and the
period map for double solids, as given by Ron Donagi},U. of Utah,
preprint.

\bib{CG}H. Clemens, P. Griffiths, {\it The intermediate Jacobian
of the cubic threefold}, Ann. Math. 95 (1972), 281-356.

\bib{D1}R. Donagi, {\it The tetragonal construction}, AMS Bull.
4 (1981), 181-185.

\bib{D2}R. Donagi, {\it The unirationality of ${\cal A}_5$},
Ann. Math. 119 (1984), 269-307.

\bib{D3}R. Donagi, {\it Big Schottky}, Inv. Math. 89 (1987),
569-599.

\bib{D4}R. Donagi, {\it Non-Jacobians in the Schottky loci},
Ann. of Math. 126 (1987), 193-217.

\bib{D5}R. Donagi, {\it The Schottky problem, in: Theory of
Moduli}, LNM 1337, Springer-Verlag (1988), 84-137.

\bib{D6}R. Donagi, {\it On the period map for Clemens' double
solids}, preprint.

\bib{DS}R. Donagi, R. Smith, {\it The structure of the Prym
map}, Acta Math. 146 (1981) 25-102.

\bib{Deb1}O. Debarre, {\it  Vari\'{e}t\'{e}s de Prym,
conjecture de la tris\'{e}cante et ensembles d'Andreotti et
Mayer},   Univ. Paris Sud,Thesis, \linebreak Orsay (1987).

\bib{Deb2}O. Debarre, {\it Sur les vari\'{e}t\'{e}s de Prym des
courbes t\'{e}tragonales}, Ann. Sci. E.N.S. 21 (1988), 545-559.

\bib{Dem}M. Demazure, {\it Seminaire sur les singularit\'{e}s des
surfaces.}  LNM  777 , Springer-Verlag (1980),  23-69.

\bib{FS}R. Friedman, R. Smith, {\it The generic Torelli Theorem
for the Prym map}, Inv. Math. 67 (1982), 473-490.

\bib{vG}B. van Geemen, {\it Siegel modular forms vanishing on the
moduli space of curves}, Inv. Math. 78 (1984), 329-349.

\bib{vGvdG}B. van Geemen, G. van der Geer, {\it Kummer varieties
and the moduli space of curves}, Am. J. of Math. 108 (1986),
615-642.

\bib{vGP}B. van Geemen, E. Previato, {\it Prym varieties and the
Verlinde formula}, MSRI preprint, May 1991.

\bib{I}E. Izadi, {\it On the moduli space of four dimensional
principally polarized abelian varieties}, Univ. of Utah
Thesis, June 1991.

\bib{K}V. Kanev, {\it The global Torelli theorem for Prym
varieties at a generic point}, Math. USSR Izvestija 20 (1983),
235-258.

\bib{M1}D. Mumford, {\it Theta characteristics on an algebraic
curve}, Ann. Sci. E.N.S. 4 (1971), 181-192.

\bib{M2}D. Mumford, {\it Prym varieties I. Contributions to
Analysis}, 325-350, New York, Acad. Press, 1974.

\bib{Ma}L. Masiewicki, {\it Universal properties of Prym
varieties
with an application to algebraic curves of genus five}, Trans.
Amer. Math. Soc. 222 (1976), 221-240.

\bib{N}J. C. Naranjo, {\it Prym varieties of bi-elliptic curves},
Univ. de Barcelona preprint no. 65, June 1989.

\bib{P}S. Pantazis, {\it Prym varieties and the geodesic flow on
$SO(n)$}, Math. Ann. 273 (1986), 297-315.

\bib{R}S. Recillas, {\it Jacobians of curves with a $g^1_4$ are
Prym varieties of trigonal curves}, Bol. Soc. Math. Mexicana 19
(1974), 9-13.

\bib{Ra}D. Radionov, {\it letter}.

\bib{SR}J.G. Semple, L. Roth, {\it Introduction to Algebraic
Geometry}, Oxford U. Press, 1949.

\bib{SV}R. Smith, R. Varley, {\it Components of the locus of
singular theta divisors of genus 5}, LNM 1124, Springer-Verlag
(1983),
338-416.

\bib{T}A. Tjurin, {\it Five lectures on three dimensional
varieties}, Russ. Math. Surv. 27 (1972).

\bib{V}R. Varley, {\it Weddle's surfaces, Humbert's curves, and a
certain 4-dimensional abelian variety}, Amer. J. Math. 108
(1986), 931-952.

\bib{Ve}A. Verra, {\it The fibre of the Prym map in genus
three}, Math. Ann. 276 (1987), 433-448.

\bib{W}G. Welters, {\it Recovering the curve data from a general
Prym variety}, Amer. J. of Math 109 (1987), 165-182.
\end{document}